\documentclass[twocolumn]{autart}   
\usepackage[british]{babel}
\usepackage[T1]{fontenc}
\usepackage{amsmath,amssymb,amsfonts}
\usepackage{graphicx}          
\usepackage{mathrsfs}
\usepackage{booktabs}          
\usepackage{multirow}          
\usepackage{array}             
\usepackage{float}             
\usepackage{enumitem}          
\usepackage{multicol}         
\usepackage{adjustbox}         
\usepackage{algorithm}         
\usepackage{algpseudocode}     
\usepackage{natbib}

\usepackage{tikz}
\usetikzlibrary{decorations.pathreplacing}
\usepackage{overpic}

\makeatletter
\def\ps@copyright{%
 \let\@oddhead\@empty
 \let\@evenhead\@empty
 \let\@oddfoot\@empty
 \let\@evenfoot\@empty
}
\makeatother

\begin{document}

\begin{frontmatter}

\title{Predictive-Switching Control of Stochastic Gene Regulatory Networks: A Contractive PIDE Framework}

\author[Valencia]{Christian Fern\'andez}
\author[Vigo]{Manuel P\'ajaro}
\author[Budapest]{G\'abor Szederk\'enyi}
\author[Valencia]{Irene Otero-Muras}

\address[Valencia]{Institute for Integrative Systems Biology, Spanish National Research Council, 46015 Valencia, Spain}
\address[Vigo]{Universidade de Vigo, Department of Mathematics, Campus Ourense 32004; CITMAga, Santiago de Compostela 15782, Spain}
\address[Budapest]{Systems and Control Laboratory, ELKH Institute for Computer Science and Control (SZTAKI), Budapest 1111, Hungary}

\begin{abstract}
This paper develops a predictive switching control algorithm for stochastic gene regulatory networks described by a Partial Integro-Differential Equation (PIDE) model, which enables direct shape control of the probability density function. Control inputs are selected from a finite candidate set to minimize a prescribed cost functional. A hybrid framework is proposed for scalability in higher-dimensional systems, using neural networks to approximate the control policy. A central theoretical contribution is a contraction-based analysis of the closed-loop PIDE dynamics. The paper establishes  $L^
1$-contractivity under the proposed control scheme, yielding formal stability guarantees and showing that the evolution of the probability density becomes progressively independent of the initial condition. Moreover, under strictly positive leakage terms, exponential convergence is obtained. The effectiveness and flexibility of the approach, together with the theoretical contractivity results, are illustrated through numerical simulations on three representative examples of increasing dimensionality.
\end{abstract}

\end{frontmatter}

\section{Introduction}

The dynamics of a stochastic gene regulatory network (GRN) are governed by a chemical master equation (CME) which is untractable in control applications of practical interest \citep{SAKURAI2022110647}. First-moment control of stochastic biomolecular systems has been successfully implemented \citep{11312970}. However, many relevant biological behaviors exhibit non-Gaussian features such as bimodality or oscillations, which require full probability density  shaping rather than mean regulation alone. Alternative approaches, including reinforcement learning \citep{Brancato2023} and adaptive pulsatile control \citep{Guarino2020Balancing}, provide flexibility but typically rely on extensive training or problem-specific tuning. Deep learning-enabled feedback has also shown empirical success in controlling gene expression \citep{Lugagne2024}. Nevertheless, there remains a need for model-based frameworks that provide both expressiveness and formal theoretical guarantees.

Control of non-Gaussian stochastic dynamical systems can be significantly enhanced by extending the objective from moment regulation to full output probability density shaping, thereby improving overall performance \citep{SUN2025112101}.
To address shape control of the full probability distribution of gene regulatory networks, we adopt Partial Integro-Differential Equation (PIDE) models \citep{Pajaro17}. These distributed-parameter descriptions approximate the CME, capturing the time evolution of the full probability distribution in closed form while remaining numerically tractable via semi-Lagrangian \citep{Pajaro18} or finite-volume methods \citep{Vaghy2024}. PIDE-based control has previously been explored using proportional–integral and adaptive schemes \citep{Vaghy2024, Fernandez2025}, and more recently via model predictive control (MPC) \citep{Faquir2025}. However, MPC formulations remain computationally demanding due to repeated optimization over high-dimensional stochastic dynamics. In contrast, switching or ON–OFF strategies restrict inputs to a finite set, offering low computational cost \citep{Menolascina2011, Oduola2017}, but typically lack predictive capability, limiting their ability to achieve complex distributional objectives. This motivates control strategies that combine the efficiency of switching policies with the foresight of model-based prediction.

This paper proposes a predictive-switching control (PSC) method for stochastic GRNs. PSC selects, at each decision step, the optimal action from a finite set by minimizing a cost functional evaluated over a prediction horizon. The approach is related to finite-control-set predictive methods in power electronics \citep{Herrera2024_OSS} and switched systems \citep{FCSMPC2024, Mosca2005}. A key feature is that decisions are driven directly by the PIDE-governed evolution of the probability density function (PDF), enabling control of distributional properties such as multimodality and mass allocation. To improve scalability, a hybrid architecture is introduced in which a neural network proposes candidate actions that are validated within the model-based optimization loop, in line with recent efforts to accelerate MPC using machine learning \citep{Nielsen20, Zhang25_arxiv, Chen22}.

Building on contraction analysis \citep{Lohmiller1998, FB-CTDS, FIORE2016279}, a central contribution of this paper is the proof of contractivity for controlled PIDE dynamics. A previous study by \cite{Canizo2018} established convergence for the PIDE model without external time-dependent inputs using relative entropy methods. The present work provides a more general analysis based on stochastic semigroup theory, explicitly incorporating control inputs and enabling time-dependent actuation. Therefore, $L^1$-contractivity of the controlled PIDE dynamics is established, together with exponential convergence under strictly positive leakage conditions. As a consequence, all trajectories converge to a unique distribution determined by the control input, ensuring robustness with respect to uncertainty in the initial conditions.

The remainder of the paper is structured as follows: in section \ref{sec:pide} the PIDE model is introduced, section \ref{sec:psc} develops a Predictive-Switching Control (PSC) framework for stochastic GRNs  incorporating a hybrid implementation that utilizes neural networks to facilitate scalability in higher-dimensional systems. In section \ref{sec:contractivity} the $L^1$-contractivity of the controlled PIDE dynamics is established under biologically grounded assumptions, providing a theoretical foundation for the stability of the closed-loop density evolution. Finally, in section \ref{sec:examples} the proposed framework is illustrated through three case studies of increasing dimensionality, demonstrating its effectiveness in achieving diverse control objectives.

\section{PIDE Model}
\label{sec:pide}

\begin{figure}[!htbp]
\centering
\includegraphics[width=0.4\textwidth]{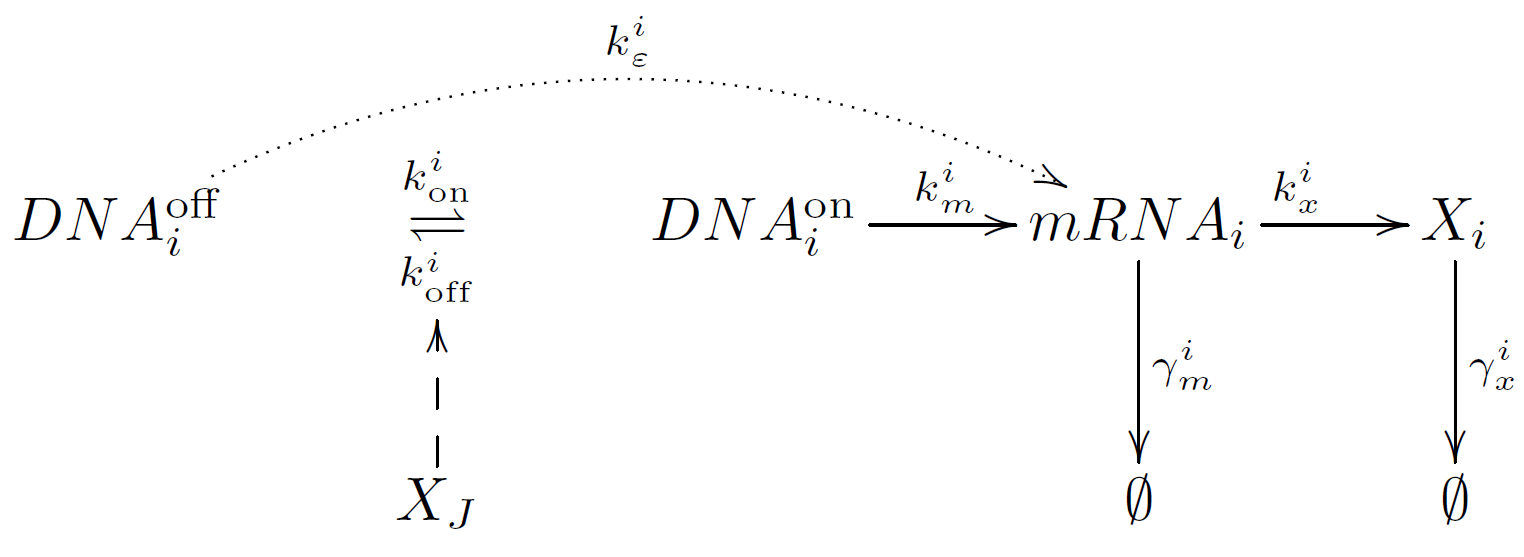} 
\caption{Schematic of the gene regulatory network. $DNA_i^{\mathrm{on/off}}$ denote active/inactive promoter states with transition rates $k_{\mathrm{on/off}}^{i}$ governed by protein $X_j$. The rates $k_m^{i}$, $k_x^{i}$, and $\gamma_{m/x}^{i}$ represent transcription, translation, and first-order degradation, respectively, while $k_{\varepsilon}^{i}$ denotes transcriptional leakage.}
\label{fig:tr_tr_scheme}
\end{figure}

Consider a gene regulatory network (GRN) consisting of $n$ genes, where $X_i$ denotes the concentration of the $i$-th protein. The system architecture, including transcription, translation, and regulatory feedback, is illustrated in Fig. \ref{fig:tr_tr_scheme}. 
Under the assumption of transcriptional bursting, where messenger RNA (mRNA) degradation is significantly faster than protein degradation ($\gamma_m^i \gg \gamma_x^i$), the dynamics of the protein concentrations $\mathbf{x} \in \mathbb{R}^n_+$ are governed by a multidimensional partial integro-differential equation (PIDE) \citep{Pajaro17}:
\begin{equation}\label{eq:PIDE_general}
\begin{split}
&\frac{\partial p(t,\mathbf{x})}{\partial t} = \sum_{i=1}^{n}\frac{\partial}{\partial x_i}\left[\gamma_x^i x_i p(t,\mathbf{x})\right] \\
&\quad + \sum_{i=1}^{n} k_m^i \int_0^{x_i} \beta_i(x_i-y_i)c_{i,\mathbf{u}}(\mathbf{y}_i)p(t,\mathbf{y}_i) \, \mathrm{d}y_i,
\end{split}
\end{equation}
\noindent where $p(t, \mathbf{x})$ is the PDF of the system state at time $t$. The first term on the right-hand side represents the drift due to first-order protein degradation at rate $\gamma_x^i$. The second term describes the stochastic production of proteins in bursts, where the kernel $\beta_i$ defines the jump size distribution as:
\begin{equation}\label{eq_burst_delta}
\beta_i(x_i - y_i) = \frac{1}{b_i}\exp \left(-\frac{x_i-y_i}{b_i} \right) -\delta(x_i - y_i),
\end{equation}
with $b_i = k_x^i/\gamma_m^i$ denotes the mean burst size and $\delta(\cdot)$ is the Dirac delta function. The vector $\mathbf{y}_i$ is defined such that $(\mathbf{y}_i)_j = x_j$ for $j \neq i$ and $(\mathbf{y}_i)_i = y_i$.

The regulatory function $c_{i,\mathbf{u}}(\mathbf{x}) \in [0, 1]$ represents the probability of gene $i$ being active, conditioned on the state $\mathbf{x}$ and the external control input $\mathbf{u}(t) = [I_1(t), \dots, I_n(t)]$, where $I_j$ are inducer concentrations. For Hill-type kinetics, the regulatory contribution of protein $X_j$ to gene $i$ is described by:
\begin{equation}
\rho_{ij}(x_j) = \frac{x_j^{H_{ij}}}{x_j^{H_{ij}}+K_{ij}^{H_{ij}}},
\end{equation}
where $H_{ij}$ and $K_{ij}$ are the cooperativity coefficient and Hill constant, respectively. Control is exerted by modulating repressive interactions via a scaling function $F_i(I_i)$:
\begin{equation}\label{eq:ind_function}
F_i(I_i) = \left[1 + \left( \frac{I_i}{\theta_{I_i}} \right)^{\mu_i}\right]^{-1},
\end{equation}
where $\theta_{I_i}$ is the half-saturation constant and $\mu_i$ the Hill coefficient. The modulated repressive contribution is given by:
\begin{equation}\label{eq:rho_modulated}
\bar{\rho}_{ij,\mathbf{u}}(x_j) = \frac{K_{ij}^{H_{ij}}}{K_{ij}^{H_{ij}} + x_j^{H_{ij}} F_i(I_i)}.
\end{equation}
To account for basal expression, a dimensionless leakage factor $\varepsilon_i = k_{\varepsilon}^i / k_m^i$ is introduced, such that the effective regulatory function becomes:
\begin{equation}\label{eq:c_with_leakage}
c_{i,\mathbf{u}}(\mathbf{x}) = \bar{\rho}_{ij,\mathbf{u}}(x_j) + \varepsilon_i (1 - \bar{\rho}_{ij,\mathbf{u}}(x_j)).
\end{equation}
The PIDE \eqref{eq:PIDE_general} is solved numerically using a semi-Lagrangian scheme on a discretized temporal grid $\{t_k\}_{k=0}^K$ over the interval $[0, T_{\text{end}}]$ with step $\Delta t_k = T_{\text{end}}/K$ and on a bounded spatial domain partitioned into $N_i$ intervals per dimension \citep{Pajaro18}.  For each protein $x_i$ is bounded by $x_{i,\max}$ and partitioned into $N_i$ intervals of size $\Delta x_i = x_{i,\max}/N_i$. The stationary solution of this model is unique and independent of the initial distribution $p(0, \mathbf{x})$ \citep{Canizo2018}.

\section{Predictive-Switching Control}
\label{sec:psc}

Predictive-Switching Control (PSC) selects inputs from a finite set of
admissible configurations at discrete time instants. Control actions are
applied jointly across all inducers, operating over the full input
configuration space. At each switching instant, the predicted evolution of
the probability distribution is evaluated under all admissible input
configurations, and the configuration that optimally satisfies a prescribed
cost functional is selected.

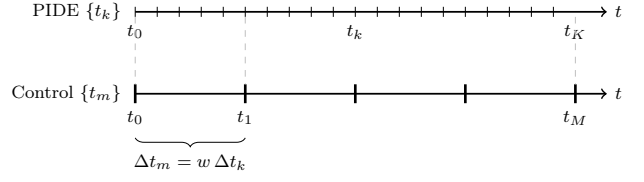
\begin{figure}[!htbp]
\centering
\resizebox{\columnwidth}{!}{
\begin{tikzpicture}[node distance=2cm, font=\scriptsize]

    \draw[->, thick] (0,1.3) -- (7.5,1.3) node[right] {$t$};
    \node[anchor=east] at (-0.1, 1.3) {PIDE $\{t_k\}$};

    \foreach \x in {0,0.35,...,7}
        \draw (\x, 1.38) -- (\x, 1.22);

    \node[below] at (0,   1.22) {$t_0$};
    \node[below] at (3.5, 1.22) {$t_k$};
    \node[below] at (7,   1.22) {$t_K$};

    \draw[->, thick] (0,0) -- (7.5,0) node[right] {$t$};
    \node[anchor=east] at (-0.1, 0) {Control $\{t_m\}$};

    \foreach \x in {0,1.75,3.5,5.25,7}
        \draw[very thick] (\x, 0.15) -- (\x, -0.15);

    \node[below] at (0,    -0.15) {$t_0$};
    \node[below] at (1.75, -0.15) {$t_1$};
    \node[below] at (7,    -0.15) {$t_M$};

    \draw[dashed, gray!50] (0,1.22)    -- (0,0.15);
    \draw[dashed, gray!50] (1.75,1.22) -- (1.75,0.15);
    \draw[dashed, gray!50] (7,1.22)    -- (7,0.15);

    \draw [decorate, decoration={brace, amplitude=4pt, mirror}]
        (0,-0.7) -- (1.75,-0.7)
        node [midway, below=4pt] {$\Delta t_m = w\,\Delta t_k$};

\end{tikzpicture}
}
\caption{Relationship between the fine numerical grid $\{t_k\}$ and the coarse switching grid $\{t_m\}$.}
\label{fig:discretization}
\end{figure}

\textbf{Control Structure and Discretization.}
The control action is held constant over an actuation window of fixed
duration $\Delta t_m = w \Delta t_k$, where $w \in \mathbb{N}_{\geq 1}$
governs the number of fine integration steps per switching interval. The
control grid $\{t_m\}_{m=0}^{M}$ constitutes a coarser partition of the
operation interval $\mathcal{T}$:
\begin{equation}
    \mathcal{T} = \bigcup_{m=0}^{M-1} (t_m,\, t_{m+1}],
\end{equation}
where $t_m = m\,\Delta t_m$. By construction,
$\{t_m\}_{m=0}^{M} \subseteq \{t_k\}_{k=0}^{K}$, and each actuation
window spans exactly $w$ integration steps of the PIDE solver
(see Fig.~\ref{fig:discretization}).

The parameter $w$ determines the sampling rate of the control action
relative to the integration step. Its selection is guided by the
characteristic time scales of the gene regulatory network, such as the
protein degradation rates $\gamma_x^i$: smaller values of $w$ increase the
temporal resolution of the switching policy, whereas larger values reduce
the number of PIDE evaluations per unit time and improve computational
efficiency.

\textbf{Admissible Input Set.}
Each external inducer $I_i$ is restricted to a binary switching logic,
taking either an inactive (OFF) state or a saturated active (ON) state.
For each inducer $I_i$, a minimal saturation level $\kappa_i$ is determined
such that the modulated repression function~\eqref{eq:rho_modulated} reaches
a prescribed suppression target at the upper boundary of the protein domain
$x_{j_{\max}}$. Specifically, for a given tolerance $\alpha_i \in (0,1)$,
the constant $\kappa_i$ is obtained by solving:
\begin{equation}\label{eq:kappa_def}
    \bar{\rho}_{ij,\mathbf{u}}(x_{j_{\max}},\, \kappa_i) = 1 - \alpha_i.
\end{equation}
The saturation vector $\boldsymbol{\kappa} = [\kappa_1, \dots, \kappa_n]$
defines the inducer magnitudes required to achieve full transcriptional
modulation across the network. Since each of the $n$ inducers independently assumes one of two states,
the admissible input space comprises $2^n$ configurations. These are encoded
in a binary structural matrix $S \in \{0,1\}^{2^n \times n}$, where each
row $S_r$ represents a distinct configuration. The element $S_{ri} = 1$
indicates that inducer $I_i$ is active at its saturation level $\kappa_i$,
and $S_{ri} = 0$ indicates that it is inactive. The corresponding physical
inducer vector for configuration $r$ is given by the Hadamard product 
\begin{equation}
    \mathbf{u}^{(r)} = \boldsymbol{\kappa} \odot S_r
    = [\kappa_1 S_{r1},\, \dots,\, \kappa_n S_{rn}], \quad r = 1, \dots, 2^n,
\end{equation}
and the finite admissible set is
$\mathcal{U} = \left\{\mathbf{u}^{(r)}\right\}_{r=1}^{2^n}$.

\textbf{Optimal Switching Policy.}
At each switching instant $t_m$, the PIDE model~\eqref{eq:PIDE_general} is
integrated over the horizon $(t_m, t_{m+1}]$ for every candidate
$\mathbf{u}^{(r)} \in \mathcal{U}$, starting from the current distribution
$p(t_m, \mathbf{x})$. This yields a set of predicted distributions
$\left\{p_{\mathbf{u}^{(r)}}(t_{m+1}, \mathbf{x})\right\}_{r=1}^{2^n}$. A cost
functional $J(\cdot)$ is evaluated for each predicted outcome, and the
optimal configuration is selected as:
\begin{equation}\label{eq:psc_opt}
    \mathbf{u}_m = \arg\max_{\mathbf{u}^{(r)} \in \mathcal{U}}\;
    J\!\left(p_{\mathbf{u}^{(r)}}(t_{m+1},\mathbf{x})\right).
\end{equation}
For minimization objectives, $\arg\max$ is replaced by $\arg\min$. The
selected configuration is applied over the entire interval
$(t_m, t_{m+1}]$, and the resulting distribution
$p_{\mathbf{u}_m}(t_{m+1}, \mathbf{x})$ serves as the initial condition
for the subsequent window. The concatenated sequence of decisions: $
    \mathbf{u} = (\mathbf{u}_0,\, \mathbf{u}_1,\, \dots,\, \mathbf{u}_{M-1})$, constitutes a piecewise-constant control trajectory over $\mathcal{T}$.
The full procedure is summarized in
Algorithm~\ref{alg:optimal_inductor_selection}.

\begin{algorithm}[H]
\caption{Predictive-Switching Control (PSC)}
\label{alg:optimal_inductor_selection}
\begin{algorithmic}[1]
    \renewcommand{\algorithmicrequire}{\textbf{Initialization:}}
    \Require Initial distribution $p(t_0, \mathbf{x})$, tolerance $\alpha$, and spatial bounds $x_{\max}$.
    \State Compute $\boldsymbol{\kappa}$ and define the matrix $S$.
    
    \For{$m = 0 : M-1$}
        \For{$r = 1 : 2^n$}
            \State Define the inducer combination $\mathbf{u}^{(r)} = \boldsymbol{\kappa} \odot S_r$.
            \State Solve the PIDE model using $\mathbf{u}^{(r)}$ over $(t_m, t_{m+1}]$ starting from $p(t_m, \mathbf{x})$.
            \State Evaluate the cost function $J_r = J(p_{\mathbf{u}^{(r)}}(t_{m+1}, \mathbf{x}))$.
            \State Store the candidate distribution $p_{\mathbf{u}^{(r)}}(t_{m+1}, \mathbf{x})$.
        \EndFor
        
        \State Select the optimal index $r^*$ and configuration $\mathbf{u}_m$:
        \[
            r^* = \arg \max_{r} J_r, \quad \mathbf{u}_m = \mathbf{u}^{(r^*)}.
        \]
        
        \State Use the stored solution as the initial distribution for the next interval $(t_{m+1}, t_{m+2}]$:
        \[
            p(t_{m+1}, \mathbf{x}) \gets p_{\mathbf{u}^{(r^*)}}(t_{m+1}, \mathbf{x}).
        \]
    \EndFor
    
    \State \textbf{Output:} Optimized inducer configurations $\mathbf{u} = (\mathbf{u}_0, \dots, \mathbf{u}_{M-1})$.
\end{algorithmic}
\end{algorithm}

\textbf{Computational Acceleration for High Dimensional Networks.}
For networks with $n \geq 3$ genes, the exhaustive evaluation of $2^n$
PIDE solutions per switching instant becomes computationally expensive,
compounded by the cost of multidimensional spatial discretization. To
address this, a data-driven proposal mechanism is integrated into the
control loop, replacing exhaustive search with a neural network
approximation of the optimal configuration, subject to a model-based
validation step.

A feedforward neural network $\mathcal{N}_\theta$ is trained offline on
datasets generated by exhaustive PSC simulations. At each switching instant
$t_m$, the network maps a compact feature vector $\mathbf{z}(t_m)$,
encoding the current probability distribution and control history, to a
continuous candidate output $\hat{\mathbf{s}} =
\mathcal{N}_\theta(\mathbf{z}(t_m))$. This output is projected onto the
discrete configuration set by rounding and identifying the corresponding
row index $r^*$ in $S$: $ r^* = \arg\min_{r}\; \|S_r - \mathrm{round}(\hat{\mathbf{s}})\|_2$,
    $\mathbf{u}^{(r^*)} = \boldsymbol{\kappa} \odot S_{r^*}$.
The network architecture and hyperparameters are selected by
cross-validation on the offline dataset. Details of the architecture and
training procedure are provided in Appendix~\ref{app:nn_psc}; the
implementation is available in the accompanying repository.

To preserve performance guarantees, the proposed configuration is validated
within a hybrid decision scheme. The PIDE model is integrated once over
$(t_m, t_{m+1}]$ using $\mathbf{u}^{(r^*)}$, and the resulting cost
$J_{r^*}$ is compared against the cost $J_{t_m}$ obtained at the previous
switching instant. If $J_{r^*} \geq J_{t_m}$, the network proposal is
accepted. Otherwise, an exhaustive search over
$\mathcal{U} \setminus \{\mathbf{u}^{(r^*)}\}$ is performed to recover the
optimal decision. The procedure is summarized in
Algorithm~\ref{alg:psc_accelerated}.

\begin{algorithm}[H]
\caption{Accelerated Predictive-Switching Control}
\label{alg:psc_accelerated}
\begin{algorithmic}[1]
    \renewcommand{\algorithmicrequire}{\textbf{Initialization:}}
    \Require Initial probability distribution \( p(t_0,\mathbf{x}) \).
    \State Load pre-trained neural network \( \mathcal{N}_\theta \).
    \State Compute \( \boldsymbol{\kappa} \) and define the configuration matrix \( S \).
    \For{\( m = 0 : M-1 \)}
        \State Construct state feature vector \( \mathbf{z}(t_m) \).
        \State Predict candidate configuration: \( \mathbf{\hat{s}} = \mathcal{N}_\theta(\mathbf{z}(t_m)) \)
        \State Identify index \( r^* \) such that \( \mathbf{S}_{r^*} = \text{round}(\mathbf{\hat{s}}) \).
        \State Map to physical inducer vector: $\mathbf{u}^{(r^*)} = \boldsymbol{\kappa} \odot \mathbf{S}_{r^*}$
        \State Solve PIDE using \( \mathbf{u}^{(r^*)} \) over \( (t_m, t_{m+1}] \).
        \State Evaluate and store candidate cost: \( J_{r^*} = J(p_{\mathbf{u}^{(r^*)}}(t_{m+1},\mathbf{x})) \). \label{step:save_j}
        \State Store the distribution $p_{\mathbf{u}^{(r^*)}}(t_{m+1}, \mathbf{x})$.
        
        \If{ \( J_{r^*} \ge J_{t_{m}} \) } 
        \Comment{Accept neural network proposal}
            \State Set $\mathbf{u}_m = \mathbf{u}^{(r^*)}$.
        \Else   \Comment{Model-based fallback to exhaustive search}
            \For{\( r = 1 : 2^n \)}
                \If{ \( r = r^* \) } \textbf{continue} \Comment{Skip redundant evaluation}
                \EndIf
                \State Define candidate \( \mathbf{u}^{(r)} =\boldsymbol{\kappa} \odot \mathbf{S}_r\).
                \State Solve PIDE using \( \mathbf{u}^{(r)} \) over \( (t_m, t_{m+1}] \).
                \State Evaluate and store candidate cost \( J_r = J(p_{\mathbf{u}^{(r)}}(t_{m+1},\mathbf{x})) \).
                \State Store the candidate distribution $p_{\mathbf{u}^{(r)}}(t_{m+1}, \mathbf{x})$.
            \EndFor
            \State Select optimal action: \( \mathbf{u}_m = \arg\max_{\mathbf{u}^{(r)}} J_r \).
        \EndIf
        \State Update initial distribution for the next window $(t_{m+1}, t_{m+2}]$ with optimal result.
    \EndFor
    
    \State \textbf{Output:} Optimized inducer trajectory $\mathbf{u} = (\mathbf{u}_0, \dots, \mathbf{u}_{M-1})$.
\end{algorithmic}
\end{algorithm}

\section{Contractivity analysis of the PIDE dynamics}
\label{sec:contractivity}
A central property of the PIDE model~\eqref{eq:PIDE_general} is $L^1$-contractivity, which provides fundamental stability guarantees for the evolution of the probability distribution. For a fixed control input, the evolution operator is non-expansive: the $L^1$ distance between any two solutions is non-increasing over time. Under the standing hypotheses of strictly positive degradation rates and strictly positive leakage, this property strengthens to exponential contractivity, establishing the robustness of the dynamics against uncertainty in initial conditions. This ensures that the long-time behaviour of the system is uniquely determined by the applied control input, regardless of the initial distribution. The formal statements and proofs are given below, relying on results from semigroup theory and Markov process theory \citep{Pazy83, EngelNagel, MeynTweedie93, Davis93, Benaim2015}.

\textbf{Notation and standing hypotheses.}
Let \(\mathbf{x}=(x_1,\dots,x_n)\in\mathbb{R}_+^n\). For a fixed inducer profile \(\mathbf{u}\) denote by \(c_{i,\mathbf{u}}(\mathbf{x})\) the corresponding input functions and by \(U_\mathbf{u}(t,\sigma)\) the propagator mapping a density at time \(\sigma\) to the density at time \(t\ge \sigma\). We assume:
\begin{itemize}\setlength{\itemindent}{1em}
  \item[\bf (H1)] Constant positive degradation rates \(\gamma^i>0\) for \(i=1,\dots,n\).
  \item[\bf (H2)] Burst kernels \(\omega_i\ge0\) with \( \displaystyle \int_0^\infty\omega_i(s)\,\mathrm{d}s=1\)
    and finite mean \(\displaystyle b_i:=\int_0^\infty s\,\omega_i(s)\,\mathrm{d}s<\infty\).
  \item[\bf (H3)] \(0\le c_{i,\mathbf{u}}(\mathbf{x})\le 1\) for all \(\mathbf{x}\); for geometric contractivity we further assume strict leakage \(c_{i,\mathbf{u}}(\mathbf{x})\ge\varepsilon_i>0\) for all \(\mathbf{x}\) and \(i\),
    and set \( \Lambda_{\min}:=\sum_{i=1}^n k_m^i\varepsilon_i>0\).
  \item[\bf (H4)] Weak solution regularity: we consider non-negative, mass-preserving weak
    solutions \(p\in C([0,\infty);L^1(\mathbb{R}_+^n))\).
\end{itemize}

\textbf{Generator and its adjoint.}
Let \(L\) denote the infinitesimal generator of the Markov process on observables; for smooth test functions \(\varphi \in C_c^\infty(\mathbb{R}_+^n)\) we write
\begin{equation*}
    L\varphi(\mathbf{x}) = A\varphi(\mathbf{x}) + B\varphi(\mathbf{x}),
\end{equation*}
with the transport term
\begin{equation*}
  A\varphi(\mathbf{x}) = -\sum_{i=1}^n \gamma^i x_i\,\partial_{x_i}\varphi(\mathbf{x}),
\end{equation*}
and the jump term
\begin{equation*}
  B\varphi(\mathbf{x}) = \sum_{i=1}^n k_m^i c_{i,\mathbf{u}}(\mathbf{x})\int_0^\infty\big(\varphi(\mathbf{x}+s \mathbf{e}_i)-\varphi(\mathbf{x})\big)\,\omega_i(s)\,\mathrm{d}s,
\end{equation*}
where \(\mathbf{e}_i\in\mathbb{R}^n\) denotes the $i$-th canonical basis vector and $s\ge 0$ is the burst size.

The adjoint \(L^\dagger\) acting on densities \(p\) (the forward Kolmogorov / generalized Fokker--Planck form) is obtained by the duality relation 
\begin{equation}\label{eq:dual_rel}
  \int_{\mathbb{R}_+^n} (L\varphi)(\mathbf{x})\,p(\mathbf{x})\,\mathrm{d}\mathbf{x}
  \;=\;
  \int_{\mathbb{R}_+^n} \varphi(\mathbf{x})\,(L^\dagger p)(\mathbf{x})\,\mathrm{d}\mathbf{x},
\end{equation}
and is given by
\begin{equation}
\label{eq:adjoint_a}
  A^\dagger p(\mathbf{x}) = \sum_{i=1}^n \partial_{x_i}\big(\gamma^i x_i p(\mathbf{x})\big), \\
  \end{equation}
and
\begin{equation}
  \label{eq:adjoint_b}
      B^\dagger p(\mathbf{x}) = \sum_{i=1}^n k_m^i \Big( (\omega_i \ast_i (c_{i,\mathbf{u}} p))(\mathbf{x})  - c_{i,\mathbf{u}}(\mathbf{x})p(\mathbf{x}) \Big),
  \end{equation}
where the convolution in the $i$-th coordinate is defined by letting $\mathbf{y}_i = (x_1, \dots, y_i, \dots, x_n)$ such that:
\begin{equation*} 
\big(\omega_i \ast_i (c_{i,\mathbf{u}} p)\big)(\mathbf{x}) = \int_0^{x_i} \omega_i(x_i-y_i) c_{i,\mathbf{u}}(\mathbf{y}_i) p(\mathbf{y}_i) \, \mathrm{d}y_i.
\end{equation*}
\textbf{Derivation of the adjoint.}
The adjoint operator $L^\dagger$ is defined through the duality relation \eqref{eq:dual_rel} for all test functions $\varphi \in C_c^\infty(\mathbb{R}_+^n)$.

\textit{Transport part ($A^\dagger$)}.
We integrate $A\varphi$ against $p$ over $\mathbb{R}_+^n = [0,\infty)^n$ and perform integration by parts in the coordinate $x_i$:
\begin{equation*}
    \int_{\mathbb{R}_+^n} A\varphi(\mathbf{x})\,p(\mathbf{x})\,\mathrm{d}\mathbf{x}
  = \int_{\mathbb{R}_+^n} -\sum_{i=1}^n \gamma^i x_i \,[\partial_{x_i}\varphi(\mathbf{x})]\,p(\mathbf{x})\,\mathrm{d}\mathbf{x}.
\end{equation*}
We consider the $i$-th summand in the right hand side in last expression and by integrating by parts we obtain
\begin{equation*}
\begin{split}
    &\int_{\mathbb{R}_+^n} -\gamma^i x_i \,[\partial_{x_i}\varphi(\mathbf{x})]\,p(\mathbf{x})\,\mathrm{d}\mathbf{x}  \\
    & =\int_{\mathbb{R}_+^{n-1}} \! \Big[ -\gamma^i x_i \varphi p \Big]_{x_i=0}^{x_i=\infty} \mathrm{d}\mathbf{x}_{i'} 
    + \int_{\mathbb{R}_+^n} \! \varphi \partial_{x_i}[\gamma^i x_i p]\,\mathrm{d}\mathbf{x},
\end{split}
\end{equation*}  
where $d\mathbf{x}_{i'}$ denotes integration over all variables except $x_i$.
The boundary term evaluates to:
\begin{equation*}
\begin{split}
    \Big[ -\gamma^i x_i \varphi(\mathbf{x}) p(\mathbf{x}) \Big]_{x_i=0}^{x_i=\infty} &= 
    \lim_{x_i\to\infty} \big(-\gamma^i x_i \varphi(\mathbf{x}) p(\mathbf{x})\big) \\
    &\quad - \lim_{x_i\to 0^+} \big(-\gamma^i x_i \varphi(\mathbf{x}) p(\mathbf{x})\big).
\end{split}
\end{equation*}
The integration by parts is understood in the distributional sense. Since
$\varphi\in C_c^\infty(\mathbb{R}_+^n)$, it is bounded and vanishes at infinity.
Moreover, the boundary flux at $x_i=0$ vanishes: for any $\varepsilon>0$,
\begin{equation*}
\int_{\mathbb R_+^{n-1}}\int_0^\varepsilon x_i\,|\varphi(\mathbf x)p(\mathbf x)|\,\mathrm{d}x_i\,\mathrm{d}\mathbf x_{i'}
\le \varepsilon \|\varphi\|_\infty \|p\|_{L^1(\mathbb R_+^n)},
\end{equation*}
which tends to zero as $\varepsilon\to0$. Hence the boundary contribution is zero, and we obtain
\begin{equation*}
    \int_{\mathbb{R}_+^n} -\gamma^i x_i \,\partial_{x_i}\varphi(\mathbf{x})\,p(\mathbf{x})\,\mathrm{d}\mathbf{x}
  = \int_{\mathbb{R}_+^n} \varphi(\mathbf{x})\,\partial_{x_i}[\gamma^i x_i p(\mathbf{x})]\,\mathrm{d}\mathbf{x}.
\end{equation*}
Summing over $i=1,\dots,n$ we prove \eqref{eq:adjoint_a}.

\textit{Jump part ($B^\dagger$)}.
We now integrate $B\varphi$ against $p$:
\begin{equation*}
\begin{split}
    &\int_{\mathbb{R}_+^n} B\varphi(\mathbf{x})\,p(\mathbf{x})\,\mathrm{d}\mathbf{x}  \\
    &=\int_{\mathbb{R}_+^n} \sum_{i=1}^n k_m^i c_{i,\mathbf{u}} \int_{0}^{\infty} \!\! \big(\varphi(\mathbf{x}+s\mathbf{e}_i)-\varphi(\mathbf{x})\big)\omega_i(s)\,\mathrm{d}s\, p\,\mathrm{d}\mathbf{x}.
\end{split}
\end{equation*}
Exchanging integrals and sums,
\begin{equation*}
\begin{split}
    &\int_{\mathbb{R}_+^n} B\varphi(\mathbf{x})\,p(\mathbf{x})\,\mathrm{d}\mathbf{x} \\
    &\quad = \sum_{i=1}^n k_m^i \Biggl( \int_{\mathbb{R}_+^n}\!\int_{0}^{\infty} c_{i,\mathbf{u}}(\mathbf{x})\varphi(\mathbf{x}+s\mathbf{e}_i)\omega_i(s)p(\mathbf{x}) \,\mathrm{d}s\mathrm{d}\mathbf{x} \\
    &\quad - \int_{\mathbb{R}_+^n} c_{i,\mathbf{u}}(\mathbf{x})\varphi(\mathbf{x})p(\mathbf{x})\,\mathrm{d}\mathbf{x} \Biggr).
\end{split}
\end{equation*}
Focus on the gain term
\begin{equation*}
    \mathcal{I}_i = \int_{\mathbb{R}_+^n}\!\int_{0}^{\infty}
        c_{i,\mathbf{u}}(\mathbf{x})\,\varphi(\mathbf{x}+s\mathbf{e}_i)\,\omega_i(s)\,p(\mathbf{x})\,\mathrm{d}s\,\mathrm{d}\mathbf{x}.
\end{equation*}
Apply the change of variables $\mathbf{y} = \mathbf{x} + s\mathbf{e}_i$.  
For each fixed $s$, this is a translation in the $i$-th coordinate, hence $\partial \mathbf{y}/\partial \mathbf{x}= \mathcal{I}_n$, $\det(\partial \mathbf{y}/\partial \mathbf{x}) = 1$, and $\mathrm{d}\mathbf{y} = \mathrm{d}\mathbf{x}$. The inverse relation is $\mathbf{x} = \mathbf{y} - s\mathbf{e}_i$, and the constraint  
$x_i = y_i - s > 0$ implies that the outer integral becomes restricted to \( s \in (0,\,y_i)\).

Substituting and applying Fubini's Theorem (see Appendix \ref{standard_theorems}),
\begin{equation*}
    \mathcal{I}_i = \int_{\mathbb{R}_+^n} \! \varphi(\mathbf{y}) \bigg[ \int_{0}^{y_i} \! \omega_i(s) c_{i,\mathbf{u}}(\mathbf{y}-s\mathbf{e}_i) p(\mathbf{y}-s\mathbf{e}_i) \mathrm{d}s \bigg] \mathrm{d}\mathbf{y}.
\end{equation*}
  Collecting terms, we find
\begin{equation*}
     \int_{\mathbb{R}_+^n} (B\varphi)(\mathbf{x})\,p(\mathbf{x})\,\mathrm{d}\mathbf{x}
  = \int_{\mathbb{R}_+^n} \varphi(\mathbf{y})\,\big(B^\dagger p\big)(\mathbf{y})\,\mathrm{d}\mathbf{y},
\end{equation*}
with
\begin{equation*}
\begin{split}
    &B^\dagger p(\mathbf{y})\\ 
   &= \sum_{i=1}^n k_m^i \Biggl[ \int_{0}^{y_i} \omega_i(s) c_{i,\mathbf{u}}(\mathbf{y}-s\mathbf{e}_i)  p(\mathbf{y}-s\mathbf{e}_i)\,\mathrm{d}s \\
    & - c_{i,\mathbf{u}}(\mathbf{y})p(\mathbf{y}) \Biggr].
\end{split}
\end{equation*}
Using the convolution notation introduced above, this expression coincides with the form stated in \eqref{eq:adjoint_b}.

\subsection*{(I) Non-expansivity in \(L^1\) for PSC}
\label{par:nonexpansive_psc}

In the implementation of the Predictive–Switching Control (PSC) we consider a temporal mesh \( t_0<t_1<\cdots<t_M\). On each interval \((t_m,t_{m+1}]\) the controller selects a mode \(\mathbf{u_m}\in\{0,1\}^n\) (a binary vector indicating which inducers are ON or OFF) and keeps that mode fixed during the whole interval.
Accordingly, we distinguish two kinds of propagators:

\begin{itemize}
\item For a fixed (time-independent) mode \(\mathbf{u}_m\) we denote by \(U_{\mathbf{u}_m}(t_{m+1},t_m)\) the evolution operator (or semigroup) that maps a density at time \(t_m\) to the density at time \(t_{m+1}\) when the mode \(\mathbf{u}_m\) is held constant on \((t_m,t_{m+1}]\).
\item For a given partition $\mathcal{P} = \{t_m\}_{m=0}^M$ with $t_0 < t_1 < \dots < t_M$, consider a switching realization $\mathbf{u}=(u_0, u_1, \dots, u_{M-1})$. We denote by $U_{\mathbf{u}}(t,\sigma)$ the composition
  \[
    U_{\mathbf{u}}(t_{M},t_0) \;=\; U_{u_{M-1}}(t_M,t_{M-1})\circ\cdots\circ U_{u_0}(t_1,t_0),
  \]
  i.e. the propagator obtained by applying successively the block-operators corresponding to each PSC decision.
\end{itemize}

The following proposition states that the evolution operators are non-expansive in $L^1$, both for a single fixed mode and for sequences of PSC switches.

\textbf{Proposition (Non-expansivity).}
Under the standing hypotheses, for every fixed mode \(\mathbf{u}_m\) the operator norm inequality
\begin{equation*}
  \|U_{\mathbf{u}_m}(t_{m+1},t_m) f\|_{L^1} \le \|f\|_{L^1}, \qquad \forall f\in L^1(\mathbb{R}_+^n),
\end{equation*}
holds. Consequently, for any PSC realization \(\mathbf{u}\) and any pair of initial densities \(p_0,q_0\),
\begin{equation*}
  \|U_{\mathbf{u}}(t_M,t_0)p_0 - U_{\mathbf{u}}(t_M,t_0)q_0\|_{L^1} \le \|p_0-q_0\|_{L^1}.
\end{equation*}

\textbf{Proof.} The proof proceeds in several detailed steps.

\medskip\noindent\textsc{Step 1. Transport semigroup \(S(t)\).}  
The transport semigroup $S(t)$ is associated with the operator $A^\dagger$. The generator $A$ acts on observables $\varphi$ as:
\begin{equation*}
    A\varphi(\mathbf{x}) := -\sum_{i=1}^n \gamma^i x_i\,\partial_{x_i}\varphi(\mathbf{x}).
\end{equation*}
The flow $\Phi_t(\mathbf{x})$ is defined by the system of ODEs $\mathrm{d}\mathbf{x}/\mathrm{d}t = \mathbf{F}(\mathbf{x})$, where the vector field $\mathbf{F}(\mathbf{x})$ is read directly from the coefficients of $A$. In this case, the field is $\mathbf{F}(\mathbf{x}) = (-\gamma^1 x_1, \dots, -\gamma^n x_n)$.

The solution to these separated ODEs $\mathrm{d}x_i/\mathrm{d}t = -\gamma^i x_i$ with the initial condition $\mathbf{x}(0) = \mathbf{x}$ is the characteristic flow:
\begin{equation*}
    \Phi_t(\mathbf{x})=(e^{-\gamma^1 t}x_1,\dots,e^{-\gamma^n t}x_n).
\end{equation*}

For any initial density \(p_0\in L^1\), define the push-forward  
\begin{equation*}
    (S(t)p_0)(\mathbf{x}) := p_0(\Phi_{-t}(\mathbf{x}))\,\big|\det D\Phi_{-t}(\mathbf{x})\big|.
\end{equation*}
To check that \(S(t)\) preserves mass, consider the integral
\begin{equation*}
    \int_{\mathbb{R}_+^n} (S(t)p_0)(\mathbf{x})\,\mathrm{d}\mathbf{x} = \int_{\mathbb{R}_+^n} p_0(\Phi_{-t}(\mathbf{x}))\,\big|\det D\Phi_{-t}(\mathbf{x})\big|\,\mathrm{d}\mathbf{x}.
\end{equation*}
Here, $\det D\Phi_{-t}(\mathbf{x})$ denotes the Jacobian determinant of the inverse flow. Perform the change of variables \(\mathbf{y}=\Phi_{-t}(x)\), so that \(\mathrm{d}\mathbf{y} = |\det D\Phi_{-t}(\mathbf{x})|\,\mathrm{d}\mathbf{x}\). Then the integral becomes
\begin{equation*}
    \int_{\mathbb{R}_+^n} p_0(\mathbf{y})\,\mathrm{d}\mathbf{y} = \int_{\mathbb{R}_+^n} | p_0(\mathbf{y})| \,\mathrm{d}\mathbf{y} = \|p_0\|_{L^1}.
\end{equation*}
Hence, \(S(t)\) preserves mass. Moreover, since \(S(t)\) is defined as the push-forward under the flow \(\Phi_t\), it also preserves positivity: if \(p_0\ge 0\) then \((S(t)p_0)(\mathbf{x})\ge 0\) for all \(\mathbf{x}\). By linearity, mass preservation, and positivity, it follows that for all \(p_0 \in L^1\)  
\begin{equation*}
    \|S(t)p_0\|_{L^1} \le \|p_0\|_{L^1},
\end{equation*}
i.e., \(S(t)\) is a contraction in \(L^1\).

Finally, standard semigroup theory asserts that \(A^\dagger\) generates a strongly continuous semigroup \(S(t)\) on a suitable domain of functions where the derivatives and boundary conditions are well-defined, with no-flux or decay at the boundaries \(x_i\to0^+\) and \(x_i\to\infty\), \citep{Pazy83,EngelNagel} (Appendix \ref{standard_theorems}).

\medskip\noindent\textsc{Step 2. Boundedness and mass-preservation of \(B^\dagger\).}
For each \(i\) and \(p\in L^1\), using \(0\le c_{i,\mathbf{u}}\le 1\), we first consider the case \(p\ge 0\). Tonelli's theorem (see Appendix \ref{standard_theorems}) yields
\begin{equation*}
\begin{aligned}
  &\|(\omega_i\ast_i (c_{i,\mathbf{u}} p))\|_{L^1}
  \le \int_{\mathbb{R}_+^n} \int_0^{x_i} \omega_i(x_i-y_i)\,p(\mathbf{y}_i)\,\mathrm{d}y_i\, \mathrm{d}\mathbf{x} \\
  &\qquad= \int_{\mathbb{R}_+^{n-1}}
  \int_0^\infty
  p(\mathbf{y}_i)
  \left(
    \int_{y_i}^\infty \omega_i(x_i-y_i)\,\mathrm{d}x_i
  \right)
  \mathrm{d}y_i\,\mathrm{d}\mathbf{x}_{i'}\\
   &\qquad= \|p\|_{L^1},
\end{aligned}
\end{equation*}
where we used \( \displaystyle \int_0^\infty \omega_i(s)\,\mathrm{d}s = 1\).

For general \(p\in L^1\), observe that
\begin{equation*}
    |(\omega_i\ast_i (c_{i,\mathbf{u}} p))(\mathbf{x})| \le (\omega_i \ast_i |c_{i,\mathbf{u}} p|)(\mathbf{x}) \le (\omega_i \ast_i |p|)(\mathbf{x}),
\end{equation*}
so that applying the previous calculation to \(|p|\) gives
\begin{equation*}
    \|(\omega_i \ast_i (c_{i,\mathbf{u}} p))\|_{L^1} \le \|p\|_{L^1}.
\end{equation*}
Using also \(\|c_{i,\mathbf{u}} p\|_{L^1} \le \|p\|_{L^1}\), we obtain:
\begin{equation*}
    \begin{split}
    \|B^\dagger p\|_{L^1} &\le \sum_{i=1}^n k_m^i \big( \|(\omega_i \ast_i (c_{i,\mathbf{u}} p))\|_{L^1} + \|c_{i,\mathbf{u}} p\|_{L^1} \big) \\
    &\le 2 \sum_{i=1}^n k_m^i \|p\|_{L^1}.
    \end{split}
\end{equation*}
Hence \(B^\dagger \in \mathcal{B}(L^1)\) is bounded.
 
To check mass-preservation, integrate \(B^\dagger p\) over \(\mathbb{R}_+^n\):
\begin{equation*}
    \begin{split}
    \int_{\mathbb{R}_+^n} (B^\dagger p)(\mathbf{x})\,\mathrm{d}\mathbf{x}
    &= \sum_{i=1}^n k_m^i \Biggl( \int_{\mathbb{R}_+^n} (\omega_i\ast_i (c_{i,\mathbf{u}} p))(\mathbf{x})\,\mathrm{d}\mathbf{x} \\
    &\quad - \int_{\mathbb{R}_+^n} c_{i,\mathbf{u}}(\mathbf{x})p(\mathbf{x})\,\mathrm{d}\mathbf{x} \Biggr).
    \end{split}
\end{equation*}
  
Consider the first term in parentheses. Using Fubini's theorem and \(\displaystyle  \int_0^\infty \omega_i(s)\,\mathrm{d}s=1\),
  \begin{equation*}
\begin{split}
    &\int_{\mathbb{R}_+^n} (\omega_i\ast_i (c_{i,\mathbf{u}} p))(\mathbf{x})\,\mathrm{d}\mathbf{x} \\
    & =\int_{\mathbb{R}_+^{n-1}} \int_0^\infty
c_{i,\mathbf{u}}(\mathbf{y}_i)p(\mathbf{y}_i)
\int_{y_i}^\infty \omega_i(x_i-y_i)\,\mathrm{d}x_i
\,\mathrm{d}y_i\,\mathrm{d}\mathbf{x}_{i'} \\
    &= \int_{\mathbb{R}_+^n} c_{i,\mathbf{u}}(\mathbf{x})p(\mathbf{x})\,\mathrm{d}\mathbf{x}.
\end{split}
\end{equation*}
Hence, for each \(i\),
\begin{equation*}
    \int_{\mathbb{R}_+^n} (\omega_i\ast_i (c_{i,\mathbf{u}} p))(\mathbf{x})\,\mathrm{d}\mathbf{x} - \int_{\mathbb{R}_+^n} c_{i,\mathbf{u}}(\mathbf{x})p(\mathbf{x})\,\mathrm{d}\mathbf{x} = 0.
\end{equation*}
Summing over \(i=1,\dots,n\) gives
\begin{equation}\label{eq:Bp0}
    \int_{\mathbb{R}_+^n} B^\dagger p(\mathbf{x})\,\mathrm{d}\mathbf{x} = 0, \qquad p\in L^1.
\end{equation}
Therefore \(B^\dagger\) preserves mass.

\medskip\noindent\textsc{Step 3. Generation of the full semigroup.}
We combine the transport generator \(A^\dagger\) and the jump operator \(B^\dagger\) to construct the full forward semigroup \(T(t)\) solving \(\partial_t p = L^\dagger p\), \(L^\dagger = A^\dagger + B^\dagger\). The generation and positivity arguments rely on the bounded perturbation theorem and the Miyadera--Voigt criterion stated in Appendix~\ref{standard_theorems} \citep{Pazy83, EngelNagel, Miyadera66, Voigt77}.

From Step 1 we have that \(A^\dagger\) generates positive $C_0$-semigroup \(S(t)\) on \(L^1(\mathbb{R}_+^n)\). In fact \(S(t)\) is the push-forward along the deterministic flow \(\Phi_t\), hence it preserves positivity and mass and satisfies the isometry identity
\begin{equation*}
    \|S(t)p\|_{L^1}=\|p\|_{L^1}\qquad\text{for all }p\in L^1,
\end{equation*}
so in particular \(\|S(t)\|_{L^1\to L^1}=1\). From Step 2 we have \(B^\dagger\in\mathcal B(L^1)\) and the operator norm bound
\begin{equation*}
     \|B^\dagger\|_{L^1\to L^1} \le 2\sum_{i=1}^n k_m^i.
\end{equation*}
Therefore the hypotheses of the bounded-perturbation theorem apply: since \(B^\dagger\) is bounded, the operator \(L^\dagger:=A^\dagger+B^\dagger\) (with domain \(D(A^\dagger)\)) is closable and its closure generates a \(C_0\)-semigroup \(T(t)\) on \(L^1\); moreover the short-time expansion is the Dyson--Phillips series above. Choosing any \(\hat{t}_0>0\) with
\begin{equation*}
    \|B^\dagger\|\,\hat{t}_0 \le 2\sum_{i=1}^n k_m^i \hat{t}_0 < 1,
\end{equation*}
ensures convergence of the series in operator norm on \([0,\hat{t}_0]\).

For the Miyadera--Voigt criterion, we note that $(S(t))$ is positive and for $p \ge 0$,
\begin{equation*}
    \|B^\dagger S(s) p\|_{L^1} 
\le  \|B^\dagger\| \, \|p\|_{L^1} 
\le 2 \sum_{i=1}^n k_m^i \, \|p\|_{L^1}.
\end{equation*}
Integrating over $s\in[0,\hat{t}_0]$ gives
\begin{equation*}
    \int_0^{\hat{t}_0} \|B^\dagger S(s)p\|_{L^1} \, \mathrm{d}s
\le 2 \sum_{i=1}^n k_m^i\hat{t}_0 \, \|p\|_{L^1}.
\end{equation*}
Hence, choosing $\hat{t}_0>0$ such that 
\begin{equation*}
    q := 2 \sum_{i=1}^n k_m^i \hat{t}_0 < 1,
\end{equation*}
the Miyadera--Voigt bound holds. 

It follows that the Dyson--Phillips series converges on the positive cone to a positive $C_0$-semigroup $T(t)$, and mass preservation follows immediately from \eqref{eq:Bp0}.

To control the Dyson--Phillips expansion, we rely on the standard convergence estimate. Fix \(p\in L^1\) and \(t\in[0,\hat{t}_0]\). Define \( \displaystyle m_k(t) := \sup_{0\le s\le t}\|T_k(s)p\|_{L^1}\). Using \(\|S(s)\|_{L^1\to L^1}=1\) and the recurrence
\begin{equation}
\begin{split}
  \|T_{k+1}(s)p\|_{L^1}
  &\le \int_0^s \|S(s-\sigma)\|\,\|B^\dagger\|\,\|T_k(\sigma)p\|_{L^1}\,\mathrm{d}\sigma \\
  &\le \|B^\dagger\| \int_0^s m_k(t)\,\mathrm{d}\sigma \le \|B^\dagger\|\,t\,m_k(t),
\end{split}
\end{equation}
we obtain \(m_{k+1}(t)\le(\|B^\dagger\|\,t)\,m_k(t)\). By induction
\(m_k(t)\le(\|B^\dagger\|\,t)^k\|p\|_{L^1}\), so for \(t\le \hat{t}_0\) with
\(\|B^\dagger\|\,\hat{t}_0<1\) the series \(\sum_{k\ge0}T_k(t)p\) converges absolutely
in \(L^1\), uniformly on \([0,\hat{t}_0]\). This gives existence of \(T(t)p\) and
continuity in \(t\) for \(t\in[0,\hat{t}_0]\); the semigroup property then extends
\(T(t)\) to all \(t\ge0\).

Finally, we check that the semigroup preserves total mass.
Write
\begin{equation*}
    (B^\dagger q)(\mathbf{x})=\sum_{i=1}^n k_m^i\big((\omega_i\ast_i(c_{i,\mathbf{u}} q))(\mathbf{x})-c_{i,\mathbf{u}}(\mathbf{x})q(\mathbf{x})\big).
\end{equation*}
From \eqref{eq:Bp0} we have \( \displaystyle \int_{\mathbb{R}_+^n}(B^\dagger q)(\mathbf{x})\,\mathrm{d}\mathbf{x}=0\) for all \(q\in L^1\). We prove by induction on \(k\ge0\) that for every \(p\in L^1\) and every \(t\ge0\),
\begin{equation*}
    \int_{\mathbb{R}_+^n} T_k(t)p(\mathbf{x})\,\mathrm{d}\mathbf{x} =
  \begin{cases}
    \displaystyle \int_{\mathbb{R}_+^n} p(\mathbf{x})\,\mathrm{d}\mathbf{x}, & k=0,\\[6pt]
    0, & k\ge1.
  \end{cases}
\end{equation*}
The base case \(k=0\) holds because \(S(t)\) preserves mass. Assume the claim for \(k\). Then, using Fubini and mass-preservation of \(S(\cdot)\),
\begin{equation*}
\begin{aligned}
  &\int_{\mathbb{R}_+^n} T_{k+1}(t)p(\mathbf{x})\,\mathrm{d}\mathbf{x}\\
  &= \int_{\mathbb{R}_+^n} \left(\int_0^t S(t-s)\big(B^\dagger T_k(s)p\big)(\mathbf{x})\,\mathrm{d}s\right) \mathrm{d}\mathbf{x} \\
  &= \int_0^t \left(\int_{\mathbb{R}_+^n} S(t-s)\big(B^\dagger T_k(s)p\big)(\mathbf{x})\,\mathrm{d}\mathbf{x}\right) \mathrm{d}s
    \\
    &= \int_0^t \left(\int_{\mathbb{R}_+^n} (B^\dagger T_k(s)p)(\mathbf{x})\,\mathrm{d}\mathbf{x}\right) \mathrm{d}s = \int_0^t 0\,\mathrm{d}s = 0,
\end{aligned}
\end{equation*}
where the penultimate equality uses \eqref{eq:Bp0}. Therefore, all \(T_k\) with \(k\ge1\) have zero total integral, and summing the series termwise yields
\begin{equation*}
\begin{split}
  &\int_{\mathbb{R}_+^n} T(t)p(\mathbf{x})\,\mathrm{d}\mathbf{x} \\
  &= \int_{\mathbb{R}_+^n} T_0(t)p(\mathbf{x})  \,\mathrm{d}\mathbf{x} + \sum_{k\ge1}\int_{\mathbb{R}_+^n} T_k(t)p(\mathbf{x}) \,\mathrm{d}\mathbf{x} \\
  &= \int_{\mathbb{R}_+^n} p(\mathbf{x})\,\mathrm{d}\mathbf{x},
  \end{split}
\end{equation*}
so \(T(t)\) preserves total mass for all \(t\ge0\).

\medskip\noindent\textsc{Step 4. Mass preservation and non-expansivity.}
Let \(p\in L^1\) with \(p\ge0\). By the previous observation \(\displaystyle  \int_{\mathbb{R}_+^n} T(t)p(\mathbf{x})\,\mathrm{d}\mathbf{x} = \int_{\mathbb{R}_+^n} p(\mathbf{x})\, \mathrm{d}\mathbf{x}\) for all \(t\ge0\). In particular, for non-negative functions \(T(t)\) preserves the \(L^1\) norm:
\begin{equation*}
      \|T(t)p\|_{L^1} = \int_{\mathbb{R}_+^n} T(t)p(\mathbf{x})\,\mathrm{d}\mathbf{x} = \int_{\mathbb{R}_+^n} p(\mathbf{x})\,\mathrm{d}\mathbf{x} = \|p\|_{L^1}.
\end{equation*}
Now let \(f\in L^1\) be arbitrary and decompose it into positive and negative parts \(f=f_+ - f_-\) (Jordan decomposition), with \(f_\pm\ge0\) and \(\|f\|_{L^1}=\|f_+\|_{L^1}+\|f_-\|_{L^1}\). Applying \(T(t)\) and using linearity:
\begin{equation*}
    T(t)f = T(t)f_+ - T(t)f_-.
\end{equation*}
By the triangle inequality in \(L^1\) and norm preservation for non-negative functions:
\begin{equation*}
\begin{aligned}
  \|T(t)f\|_{L^1}
  &\le \|T(t)f_+\|_{L^1} + \|T(t)f_-\|_{L^1} \\
  &= \|f_+\|_{L^1} + \|f_-\|_{L^1} = \|f\|_{L^1}.
\end{aligned}
\end{equation*}
Hence \(T(t)\) is contractive in \(L^1\) for all \(t\ge0\), i.e., \(\|T(t)\|_{L^1\to L^1}\le 1\).

\medskip\noindent\textsc{Step 5. From fixed-mode operators to PSC evolution.}
For each block \(  (t_m,t_{m+1}]\) with fixed mode \(\mathbf{u}_m\) the evolution is given by the operator \(U_{\mathbf{u}_m}(t_{m+1},t_m)=T_{\mathbf{u}_m}(\Delta t_m)\), which is contractive by the previous step. Then the total evolution generated by the PSC policy \(\mathbf{u}=(u_0,\ldots,u_{M-1})\) over the partition \(t_0<t_1<\cdots<t_M\) is the composition
\begin{equation*}
    U_{\mathbf{u}}(t_M,t_0)=T_{\mathbf{u}_{M-1}}(\Delta t_{M-1})\circ\cdots\circ T_{\mathbf{u}_0}(\Delta t_0).
\end{equation*}
Since a composition of contractive operators is contractive, for any initial densities \(p_0,q_0\) we have
\begin{equation*}
\begin{split}
    &\|U_{\mathbf{u}}(t_M,t_0)p_0 - U_{\mathbf{u}}(t_M,t_0)q_0\|_{L^1}\\
  &= \|U_{\mathbf{u}}(t_M,t_0)(p_0-q_0)\|_{L^1}
  \le \|p_0-q_0\|_{L^1}.
  \end{split}
\end{equation*}
This inequality expresses the trajectory form of non-expansivity: the \(L^1\) distance between two solutions evolving under the \emph{same} PSC realization cannot increase over time.

This completes the proof of non-expansivity in the PSC context.
\hfill $\Box$

\subsection*{(II) Geometric contractivity under PSC profiles}
\label{par:geomPSC_full}

\textbf{Theorem (Geometric contractivity for PSC).}
Assume the standing hypotheses and, in addition, that $c_{i,\mathbf{u}}(\mathbf{x})\ge \varepsilon_i>0$ for all $\mathbf{x}$ and $i$ (so $\Lambda_{\min}>0$). Let $\mathbf{u}(t)$ be a piecewise constant inducer profile, i.e.\ a sequence of constant controls on intervals $(t_m,t_{m+1}]$. Then there exist constants $K\geq 1$ and $\phi>0$ (depending only on the model parameters, but independent of the switching sequence) such that for any densities $p_0,q_0\in L^1(\mathbb{R}_+^n)$ and all $t\ge0$,
\begin{equation*}
  \|U_\mathbf{u}(t,0)p_0 - U_\mathbf{u}(t,0)q_0\|_{L^1} \;\le\; K\,e^{-\phi t}\,\|p_0-q_0\|_{L^1}.
\end{equation*}

\textbf{Proof.} 
The proof proceeds through a sequence of steps that explicitly construct the drift and minorization properties and then combine them to obtain a global geometric contractivity estimate. The strategy is based on Harris' classical ergodic theorem \citep{Harris1956} and its modern formulation by Meyn--Tweedie \citep{MeynTweedie93}, which states that for a Markov semigroup on a measurable space, the existence of a Lyapunov function $V$ and a small set $C$ satisfying a uniform minorization condition guarantees geometric ergodicity. Specifically, the evolution is contractive in the weighted norm $\|\cdot\|_V$ defined by the Lyapunov function (see Appendix \ref{standard_theorems}). Since $V(\mathbf{x}) \ge 1$ for all $\mathbf{x} \in \mathbb{R}_+^n$, this weighted contractivity directly implies the geometric convergence in the $L^1$ norm stated in the Theorem.

The present construction of drift (Step~1) and minorization (Step~2) for piecewise-constant PSC profiles directly verifies the conditions of the theorem, allowing us to conclude the global geometric contractivity stated in the theorem.

\medskip\noindent\textsc{Step 1. Lyapunov inequality.}
Define the Lyapunov function
\begin{equation}
     V(\mathbf{x}) := 1 + \sum_{i=1}^n x_i \ge 1.
\end{equation}
For each coordinate $i$, the infinitesimal generator $L$ acting on $x_i$ reads
\begin{equation}
    L x_i = -\gamma^i x_i + k_m^i c_{i,\mathbf{u}}(\mathbf{x})\,b_i,
\end{equation}
where $b_i = \displaystyle \int_0^\infty s\,\omega_i(s)\,ds < \infty$ is the expected jump size.
By the linearity of the generator, applying $L$ to $V(\mathbf{x})$ and using the fact that $c_{i,\mathbf{u}}(\mathbf{x}) \le 1$ for all $i$, we have

\begin{equation*}
    \begin{split}
        L V(\mathbf{x}) &= L\left(1 + \sum_{i=1}^n x_i\right) = \sum_{i=1}^n L x_i \\
        &= -\sum_{i=1}^n \gamma^i x_i + \sum_{i=1}^n k_m^i c_{i,\mathbf{u}}(\mathbf{x})b_i \\
        &\le -\gamma_{\min} \sum_{i=1}^n x_i + \sum_{i=1}^n k_m^i b_i \\
        &=-\gamma_{\min} (V(\mathbf{x}) - 1) + \sum_{i=1}^n k_m^i b_i \\
        &= -\gamma_{\min} V(\mathbf{x}) + \left( \gamma_{\min} + \sum_{i=1}^n k_m^i b_i \right),
    \end{split}
\end{equation*}

where $\displaystyle \gamma_{\min} := \min_i \gamma^i > 0$ and we have used the fact that $c_{i,\mathbf{u}}(\mathbf{x}) \le 1$. By defining the constants $a := \gamma_{\min}$ and $\displaystyle b := \gamma_{\min} + \sum_{i=1}^n k_m^i b_i < \infty$, the inequality simplifies to the standard drift form:
\begin{equation}
    L V(\mathbf{x}) \le -a V(\mathbf{x}) + b.
\end{equation}
These constants are independent of the inducer profile $\mathbf{u}$. This establishes a uniform drift inequality, ensuring that the process is pulled back towards compact subsets at rate~$a$, independently of the control profile.

\medskip\noindent\textsc{Step 2. Minorization condition on a small set.}
Fix $R>0$ and define the compact set $\displaystyle C := \{ \mathbf{x}\in\mathbb{R}_+^n : \sum_{i=1}^n x_i \le R \}$. We now construct a uniform minorization condition on $C$. 

To ensure the existence of a transition density, we assume (as is standard in gene expression models with continuous protein levels) that at least one jump kernel $\omega_{i_0}$ is lower-bounded by a strictly positive constant $\eta > 0$ on a subinterval $[s_0, s_0+\delta]$. Let $t_{\mathrm{min}}>0$ be a fixed time horizon. For any $\mathbf{x} \in C$, consider the set of sample paths that experience exactly one jump of type $i_0$ at some time $\tau \in [\tau_1, \tau_2] \subset (0, t_{\mathrm{min}})$, and no other jumps in $[0, t_{\min}]$.

The state at time $t_{\min}$ following such a path is given by the flow:
\begin{equation*}
    \mathbf{X}_{t_{\min}} = \Phi_{t_{\min}-\tau}\left( \Phi_{\tau}(\mathbf{x}) + s \mathbf{e}_{i_0} \right),
\end{equation*}
where $\mathbf{e}_{i_0}$ is the $i_0$-th unit vector and $s$ is the jump size. Since the jump size $s$ possesses a density $\omega_{i_0}(s) \ge \eta > 0$, and the mapping $s \mapsto \mathbf{X}_{t_{\min}}$ is a diffeomorphism for a fixed $\tau$, the transition probability $U_\mathbf{u}(t_{\mathrm{min}},0)(\mathbf{x}, d\mathbf{y})$ possesses a partial density component with respect to the Lebesgue measure. 

By choosing $s_0$ sufficiently large, the deterministic contraction $\Phi_{t_{\min}}$ is compensated, ensuring that the reachable states from any $\mathbf{x} \in C$ cover a fixed hyperrectangle $A \subset \mathbb{R}_+^n$ independent of $\mathbf{x}$. Under the standing hypotheses, the jump rates 
\begin{equation*}
    \lambda(\mathbf{y})=\sum_{i=1}^n k_m^i c_i(\mathbf{y})
\end{equation*}
are bounded ($0 < \Lambda_{\min}=\sum_{i=1}^n k_m^i\varepsilon_i \le \lambda(\mathbf{y}) \le \Lambda_{\max}=\sum_{i=1}^n k_m^i< \infty$). Combining the lower bounds for the survival probability 
$e^{-\Lambda_{\max} t_{\min}}$, the jump rate 
$k_m^{i_0} \varepsilon_{i_0}$, and the jump-size density 
$\omega_{i_0}(s) \ge \eta$ on $[s_0,s_0+\delta]$, and integrating over the admissible jump time interval $[\tau_1,\tau_2]$, the restriction of the transition kernel to the event of a single $i_0$-type jump induces a component that is absolutely continuous with respect to the Lebesgue measure on a fixed reachable set $A \subset \mathbb{R}_+^n$.

Moreover, the mapping from jump size to final state defines a smooth, non-degenerate transformation with uniformly bounded Jacobian on $C \times [\tau_1,\tau_2]$, ensuring that the pushforward of the jump-size density induces a component that is uniformly bounded below on $A$.

Therefore, there exists a constant $\tilde{\varepsilon}_0 > 0$ such that
\begin{equation*}
    U_\mathbf{u}(t_{\mathrm{min}},0)(\mathbf{x}, B) 
    \ge \tilde{\varepsilon}_0 \operatorname{Leb}(B \cap A), \quad \forall \mathbf{x} \in C,
\end{equation*}
with $ B \in \mathcal{B}(\mathbb{R}_+^n)$. Defining the probability measure $\nu(B) = \operatorname{Leb}(B \cap A) / \operatorname{Leb}(A)$ and the minorization constant $\varepsilon_0 = \tilde{\varepsilon}_0 \operatorname{Leb}(A) > 0$, we obtain the uniform minorization condition:
\begin{equation*}
    U_\mathbf{u}(t_{\mathrm{min}},0)(\mathbf{x}, \cdot) \ge \varepsilon_0 \nu(\cdot), \quad \forall \mathbf{x} \in C,
\end{equation*}
where $\varepsilon_0$ and $\nu$ depend only on the system parameters and the chosen $t_{\min}$, but are independent of the specific control realization $\mathbf{u}_m$.

\medskip\noindent\textsc{Step 3. Contraction over a single regeneration interval.}
Consider two probability densities \(p\) and \(q\), and denote by
\(P := U_{\mathbf{u}}(t_{\min},0)\) the transition kernel over one time interval of length \(t_{\min}\).

From Step~2, we have the minorization condition on the small set \(C\):
\begin{equation*}
    P(\mathbf{x},\cdot)\ge \varepsilon_0\,\nu(\cdot), 
    \qquad \forall \mathbf{x}\in C.
\end{equation*}
This implies that for all \(\mathbf{x}\in C\), the kernel admits the decomposition
\begin{equation*}
    P(\mathbf{x},\cdot)
    = \varepsilon_0\,\nu(\cdot) + (1-\varepsilon_0)\,R(\mathbf{x},\cdot),
\end{equation*}
where \(R(\mathbf{x},\cdot)\) is a probability kernel.
\noindent
Let \(\mu,\eta\) be probability measures supported on \(C\). Then, by linearity of the action of \(P\),
\begin{equation*}
\begin{aligned}
\mu P - \eta P
&= (1-\varepsilon_0)(\mu R - \eta R),
\end{aligned}
\end{equation*}
and therefore, using that \(R\) is a Markov kernel,
\begin{equation*}
    \|P\mu - P\eta\|_{TV}
    \le (1-\varepsilon_0)\,\|\mu - \eta\|_{TV}.
\end{equation*}
\noindent
When \(\mu,\eta\) admit densities \(p,q\) supported on \(C\), using
\(\|\mu-\eta\|_{TV}=\tfrac{1}{2}\|p-q\|_{L^1}\), the above inequality becomes
\begin{equation*}
    \|Pp - Pq\|_{L^1}
    \le (1-\varepsilon_0)\,\|p-q\|_{L^1}.
\end{equation*}
\noindent
Recalling that \(P := U_{\mathbf{u}}(t_{\min},0)\), we finally obtain
\begin{equation*}
    \|U_{\mathbf{u}}(t_{\min},0)p - U_{\mathbf{u}}(t_{\min},0)q\|_{L^1}
    \le (1-\varepsilon_0)\,\|p-q\|_{L^1}.
\end{equation*}
\noindent
This contraction holds only at regeneration times, i.e. for states that are in (or have returned to) the small set \(C\). In general, trajectories evolve outside \(C\), and the drift condition from Step~1 ensures return to \(C\) in finite time with controlled moments.

\medskip
\noindent
Combining the drift inequality \(LV(\mathbf{x}) \le -aV(\mathbf{x}) + b\) with the above minorization condition, the pair \((V,C)\) satisfies the classical Harris ergodicity framework, as developed in its modern form by Meyn and Tweedie \cite{MeynTweedie93}, building on the foundational work of Harris \cite{Harris1956}. This yields geometric ergodicity of the semigroup and guarantees exponential convergence in a weighted norm, from which the \(L^1\)-contraction in the statement follows via \(V(\mathbf{x})\ge 1\).

\medskip\noindent\textsc{Step 4. Continuous-time estimate and concatenation.} 
Let $t\ge 0$ and write $t = N t_{\mathrm{min}} + r$ with $0\le r < t_{\mathrm{min}}$. Using interval contraction for $N$ intervals and the non-expansivity of $U_\mathbf{u}(r,0)$ on the residual interval, we have
\begin{equation}
\begin{split}
&\|U_\mathbf{u}(t,0)p_0 - U_\mathbf{u}(t,0)q_0\|_{L^1}\\
&= \|U_\mathbf{u}(r,0)\,(U_\mathbf{u}(t_{\min},0))^N p_0\\
  & - U_\mathbf{u}(r,0)\,(U_\mathbf{u}(t_{\min},0))^N q_0\|_{L^1} \\
&\le \|(U_\mathbf{u}(t_{\min},0))^N p_0
     - (U_\mathbf{u}(t_{\min},0))^N q_0\|_{L^1} \\
&\le (1-\varepsilon_0)^N \|p_0 - q_0\|_{L^1}.
\end{split}
\end{equation}

Since
\begin{equation}
(1-\varepsilon_0)^N
= \exp\!\left(N \log(1-\varepsilon_0)\right)
= \exp\!\left(-\phi N t_{\min}\right),
\end{equation}

with $\phi := -\frac{1}{t_{\min}}\log(1-\varepsilon_0) > 0,$
and using $N t_{\min} \ge t - t_{\min}$, we deduce
\begin{equation}
(1-\varepsilon_0)^N \le e^{\phi t_{\min}} e^{-\phi t}.
\end{equation}

Hence, defining $K := e^{\phi t_{\min}}\geq 1$, we obtain the continuous-time estimate
\begin{equation}
\|U_\mathbf{u}(t,0)p_0 - U_\mathbf{u}(t,0)q_0\|_{L^1}
\le K e^{-\phi t}\,\|p_0 - q_0\|_{L^1}.
\end{equation}

Finally, we apply this result to the specific PSC operation. Consider the partition $t_0 < t_1 < \dots < t_M$. On each interval $(t_m, t_{m+1}]$, the control profile $\mathbf{u}(t)$ is fixed to a constant mode $\mathbf{u}_m \in \{0,1\}^n$. The total propagator is the composition of the individual block operators:
\begin{equation}
U_{\mathbf{u}}(t_M,t_0) = U_{\mathbf{u}_{M-1}}(t_M,t_{M-1}) \circ \cdots \circ U_{\mathbf{u}_0}(t_1,t_0).
\end{equation}
Since the minorization constant $\varepsilon_0$ (and thus the decay parameters $K, \phi$) depends only on the global bounds of the system parameters and is independent of the specific vector mode active at any instant, the geometric decay established above holds for the switched flow. Substituting the total elapsed time $t = t_M - t_0$ into the estimate, we conclude:
\begin{equation}
\begin{split}
  \|U_\mathbf{u}(t_M,t_0)p_0 - U_{\mathbf{u}}(t_M,t_0)q_0\|_{L^1} \\\le K\,e^{-\phi (t_M - t_0)}\,\|p_0-q_0\|_{L^1}.
  \end{split}
\end{equation}
\hfill$\Box$

\section{Application examples}
\label{sec:examples}

Three representative examples with increasing dimensionality are presented. In each case, the PIDE model~\eqref{eq:PIDE_general} is solved numerically using a semi-Lagrangian scheme, with time reported in dimensionless units $\tau = \gamma_x t$. The non-expansivity and geometric contractivity properties are validated across all cases by applying the PSC inducer profile from a reference scenario, without modification, to two additional initial distributions. Pairwise $L^1$ distances between the resulting controlled solutions are tracked over time. Since all trajectories share the same PSC profile, any observed convergence is attributable solely to the system’s intrinsic contractive dynamics.

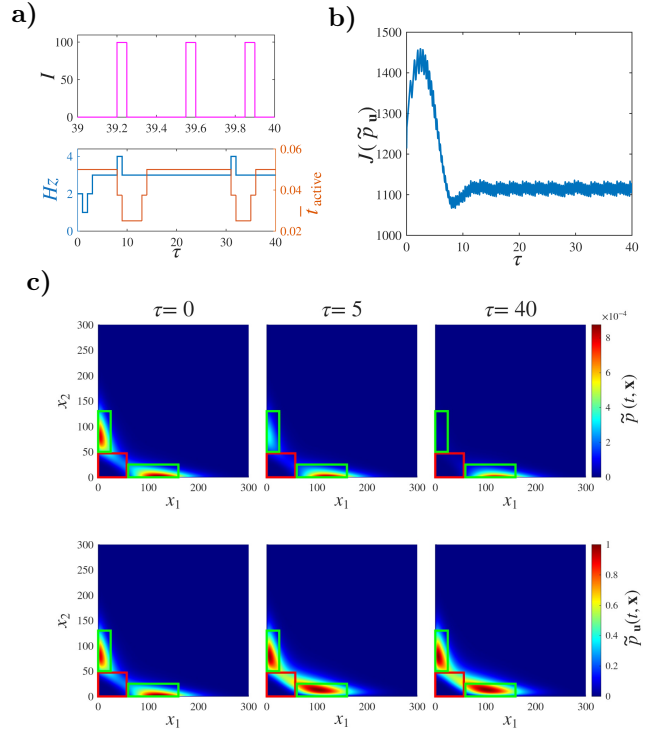
\begin{figure}[!htbp]

    \centering

        \begin{tikzpicture}
        
        \node[anchor=north west, inner sep=0] (imgA) at (-0.2, 6.9) {
            \includegraphics[width=0.21\textwidth]{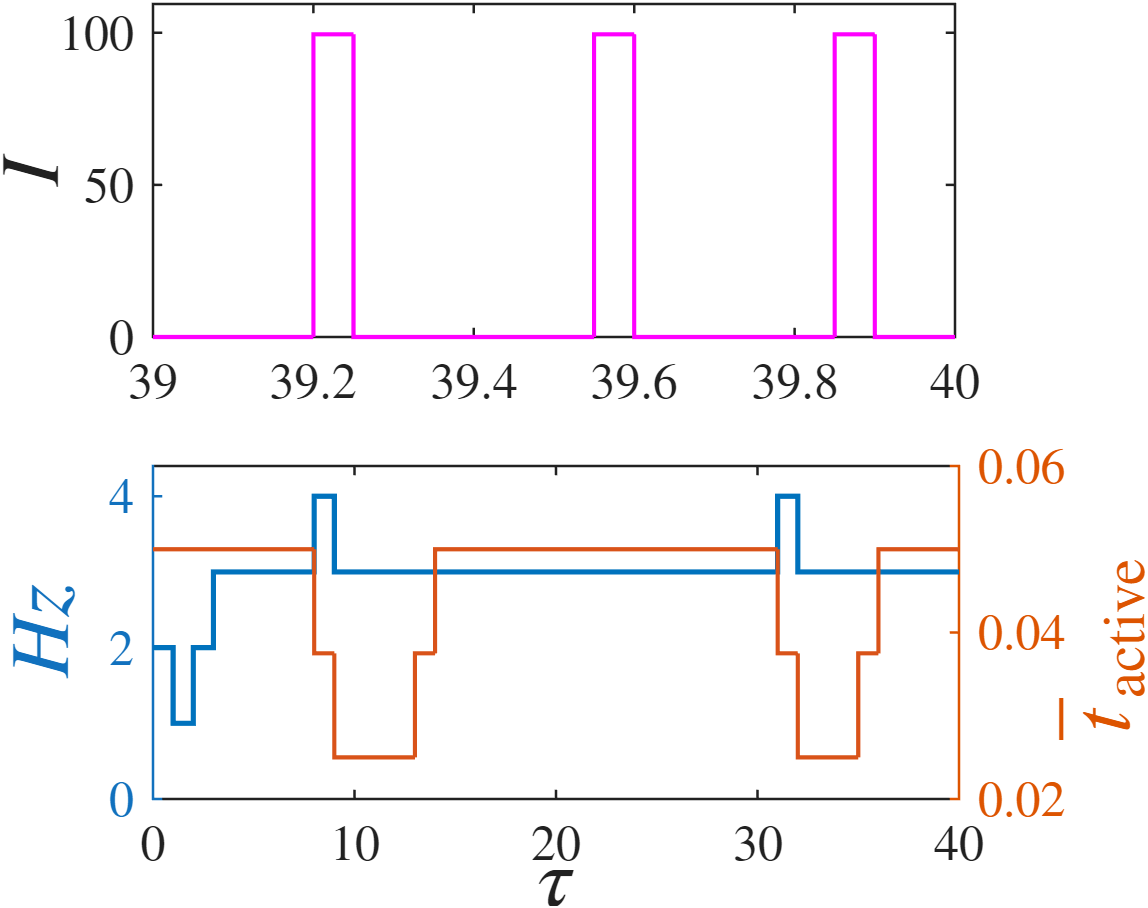}
        };
        \node[anchor=south west, font=\bfseries] at ([xshift=-0.5cm, yshift=-0.05cm]imgA.north west) {a)};

        \node[anchor=north west, inner sep=0] (imgB) at (0.225\textwidth, 7) {
            \includegraphics[width=0.21\textwidth]{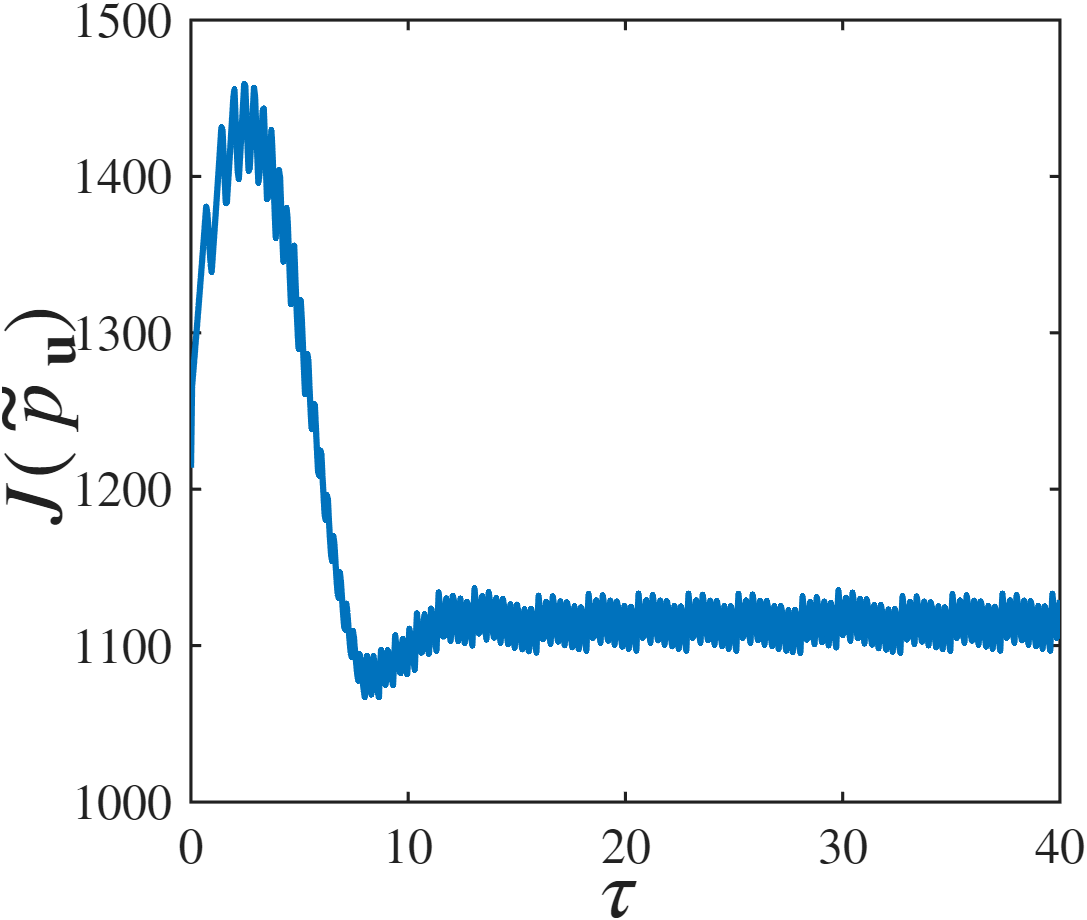}
        };
        \node[anchor=south west, font=\bfseries] at ([xshift=-0.5cm, yshift=-0.2cm]imgB.north west) {b)};

        \node[anchor=north west, inner sep=0] (imgC1) at (0, 3.5) {
            \includegraphics[width=0.45\textwidth]{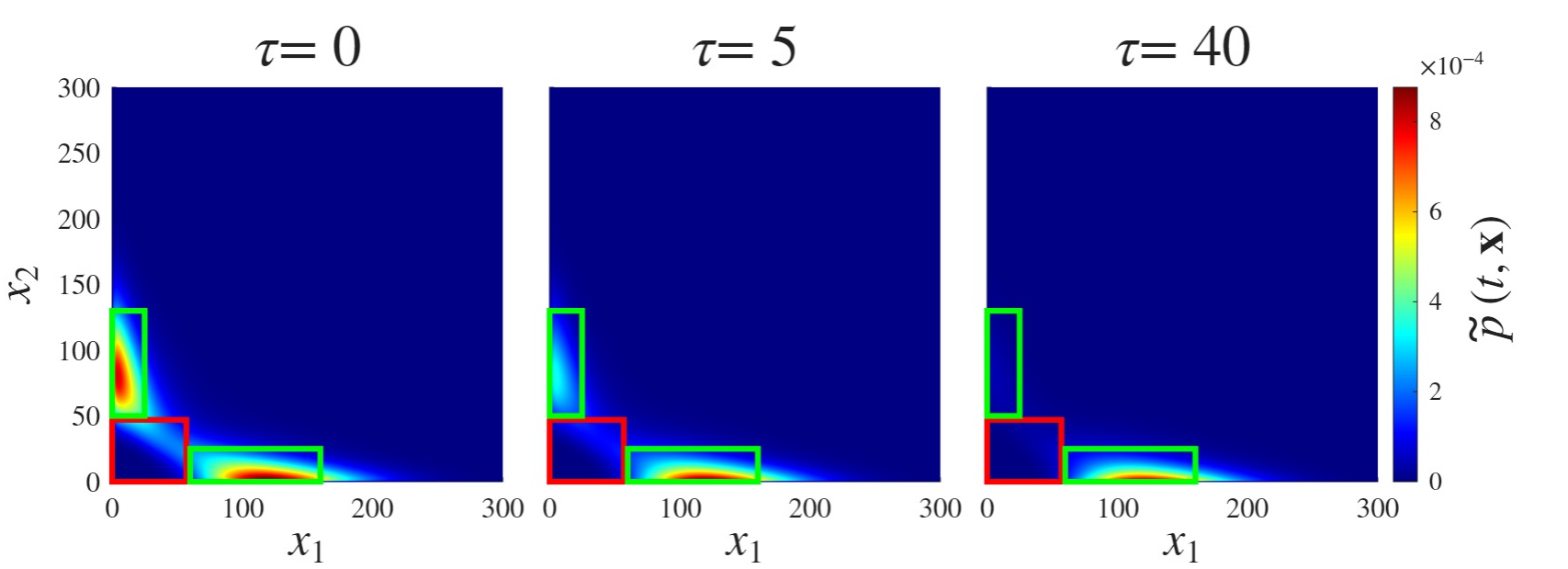}
        };
        
        \node[anchor=north west, inner sep=0] (imgC2) at (0, 0.6) {
            \includegraphics[width=0.45\textwidth]{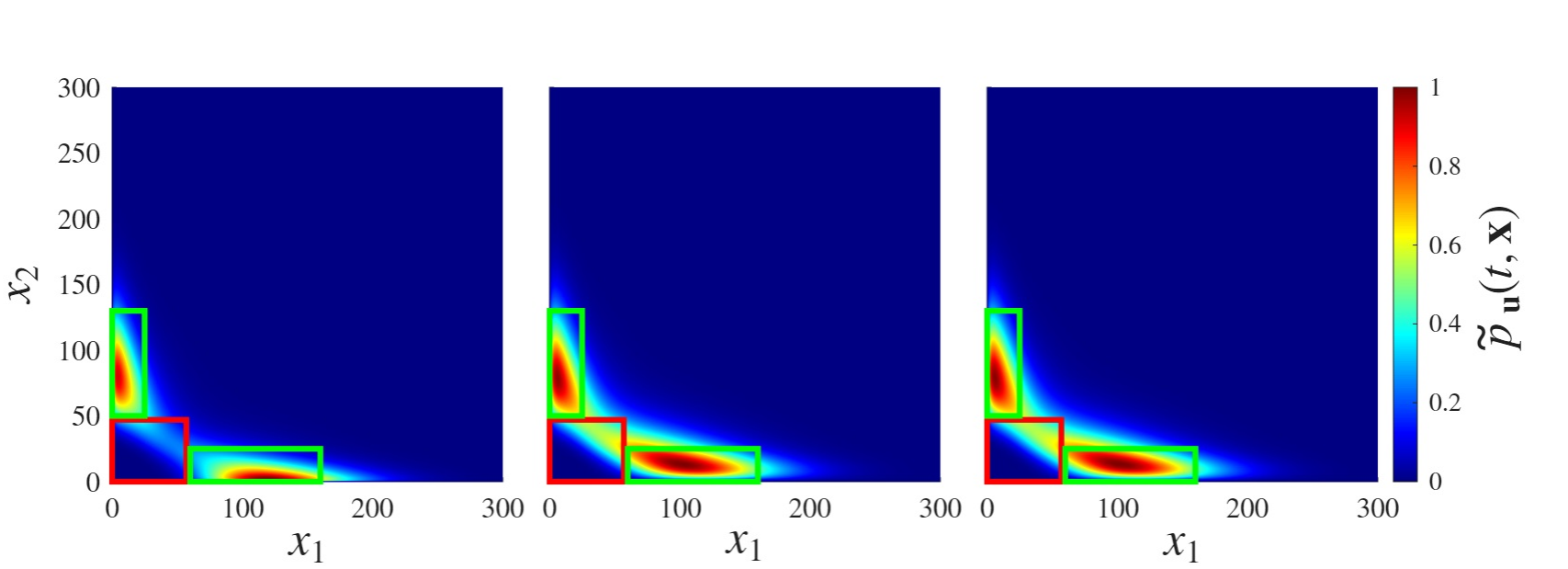}
        };

        \node[anchor=south west, font=\bfseries] at ([xshift=-0.5cm, yshift=-0.2cm]imgC1.north west) {c)};

    \end{tikzpicture}

   \caption{
(a) Input signal. Top: final interval.
Bottom: full dynamics showing activation frequency (blue) and
mean activation time (orange).
(b) Evolution of the cost functional.
(c) PDF for system
uncontrolled (top row) and PSC-controlled (bottom row). Green and
red squares are positive and negative cost contributions.}
\label{fig:Results1}
\end{figure}

\textbf{Case Study I: Bimodal State Preservation.}
First, an asymmetric synthetic genetic toggle switch \citep{Gardner2000} with two mutually repressing proteins $x_1$ and $x_2$ is considered , described by the PIDE~\eqref{eq:PIDE_general} with parameters  $k_m^1 = 11$, $k_m^2 = 9$; 
    $k_x^1 = 100$, $k_x^2 = 80$; 
    $\gamma_m^1 = \gamma_m^2 = 8.4$; 
    $\gamma_x^1 = \gamma_x^2 = 1$, $K_{12} = 30$, $K_{21} = 32$, $H_{12} = H_{21} = 4$,
$\theta_{I_2} = 0.1$, $\mu_{I_2} = 2$, $\varepsilon_1 = \varepsilon_2 =
0.1$. The initial condition is concentrated in a region of high $x_2$ and low $x_1$. The PIDE is solved over
$\Omega = [0, 300] \times [0, 300]$, In the absence of control, the system converges to a stationary distribution in which $x_1$ dominates (Fig.~\ref{fig:Results1} c).
%
%
During the transient phase, the distribution is bimodal at $\tau = 10$, with modes corresponding to $x_1$- and $x_2$-dominated regimes. The control objective is to preserve this bimodality  by counteracting the drift towards $x_1$-dominance.
A single control input $u(t) = I_2(t)$ is introduced in the repression of $x_2$ by $x_1$. The spatial and temporal discretization steps are $\Delta x = 1$ and $\Delta t_k = 0.005$. The cost functional rewards probability mass in two target modal regions $\Omega_1$ and $\Omega_2$ while penalizing the intermediate region $\Omega_c$:
\begin{equation*}
    J(\widetilde{p}_{\mathbf{u}}(t,\mathbf{x}))= \int_{\Omega_1 \cup \Omega_2} \widetilde{p}_{\mathbf{u}}(t,\mathbf{x})\,\mathrm{d}\mathbf{x} -2 \int_{\Omega_c} \widetilde{p}_{\mathbf{u}}(t,\mathbf{x})\,\mathrm{d}\mathbf{x},
\end{equation*}

where $\widetilde{p}_{\mathbf{u}}(t,\mathbf{x}) =
p_{\mathbf{u}}(t,\mathbf{x}) / \max_{\mathbf{x} \in \Omega}
p_{\mathbf{u}}(t,\mathbf{x})$ denotes the distribution normalized by its maximum value (see Fig.~\ref{fig:Results1}c). The inducer saturation level is set to $\kappa_2 = 99.5$,
corresponding to $\alpha = 0.01$ in~\eqref{eq:kappa_def}.
Algorithm~\ref{alg:optimal_inductor_selection} is applied with actuation window $w = 10$,
$\Delta t_m = 10\,\Delta t_k$. The results are shown in Fig.~\ref{fig:Results1}. The uncontrolled system converges to a unimodal distribution concentrated at $x_1$-dominance, whereas the controlled system preserves bimodality. The closed-loop simulation time is $21.35$~s.

\begin{figure}[!htbp]
    \centering

    \makebox[\textwidth][l]{\textbf{a)}}\\[-0.3ex]
    \includegraphics[width=0.4\textwidth]{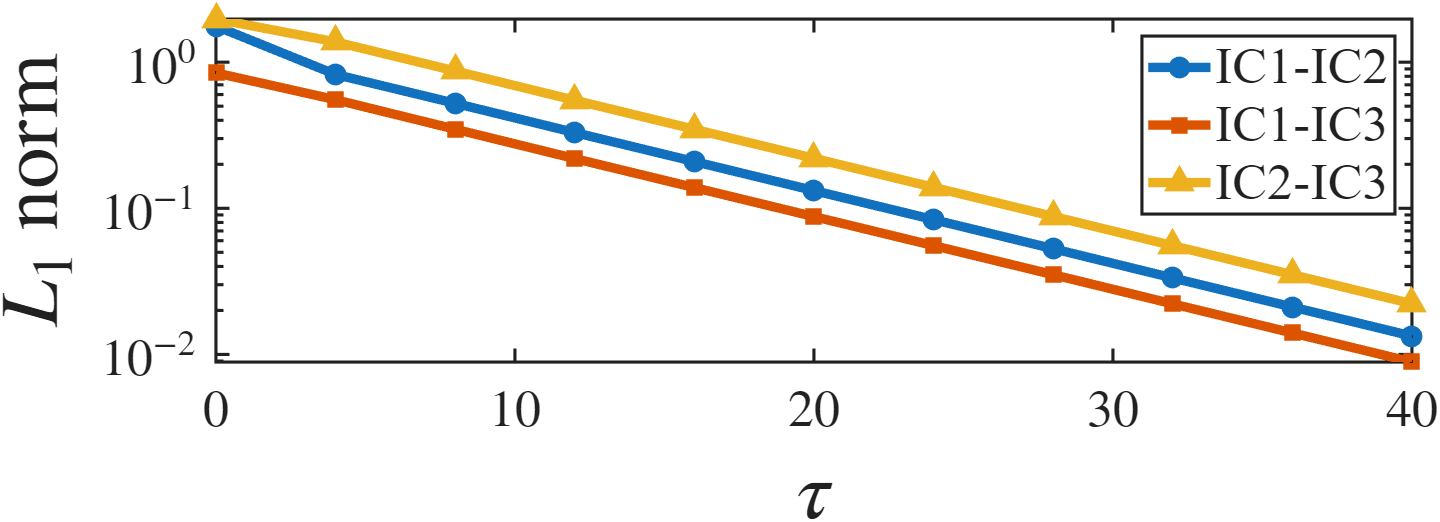}

    \makebox[\textwidth][l]{\textbf{b)}}\\[-0.3ex]
    \includegraphics[width=0.45\textwidth]{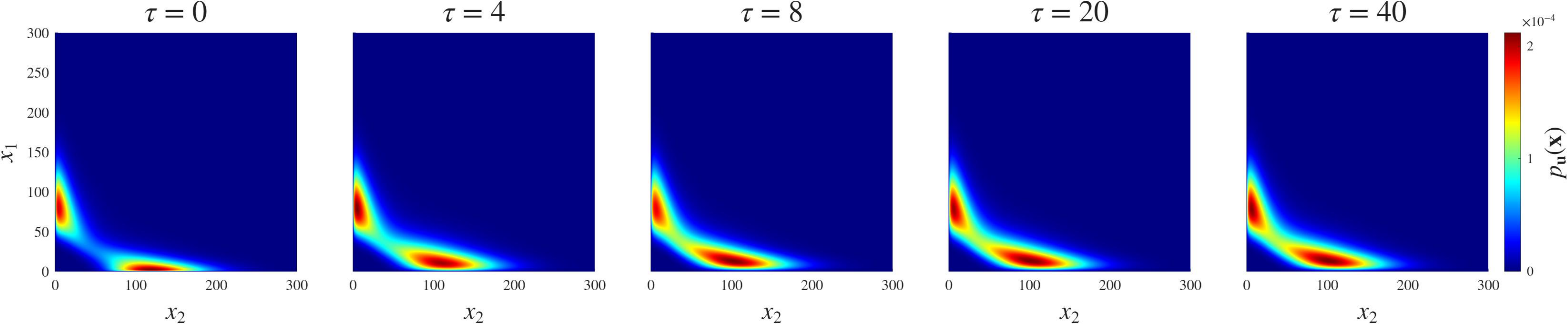}\\
    \includegraphics[width=0.45\textwidth]{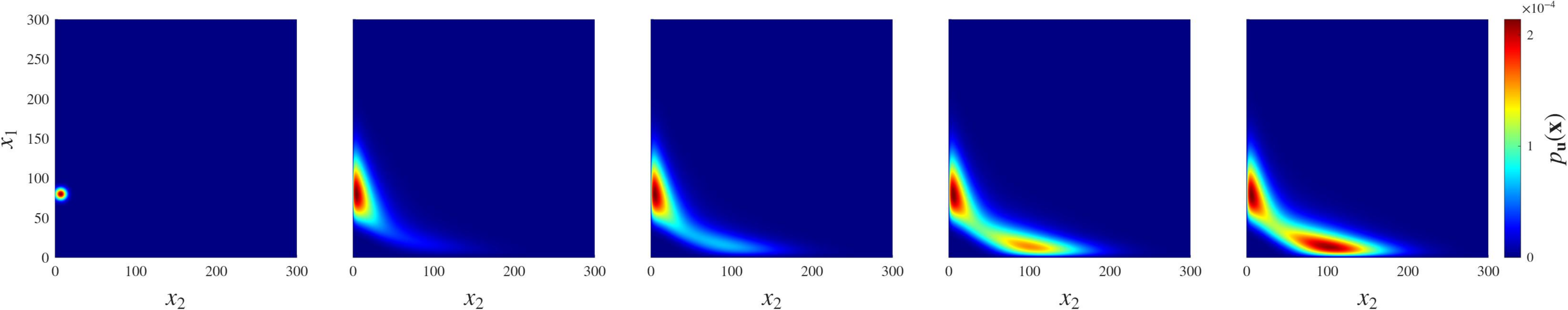}\\
    \includegraphics[width=0.45\textwidth]{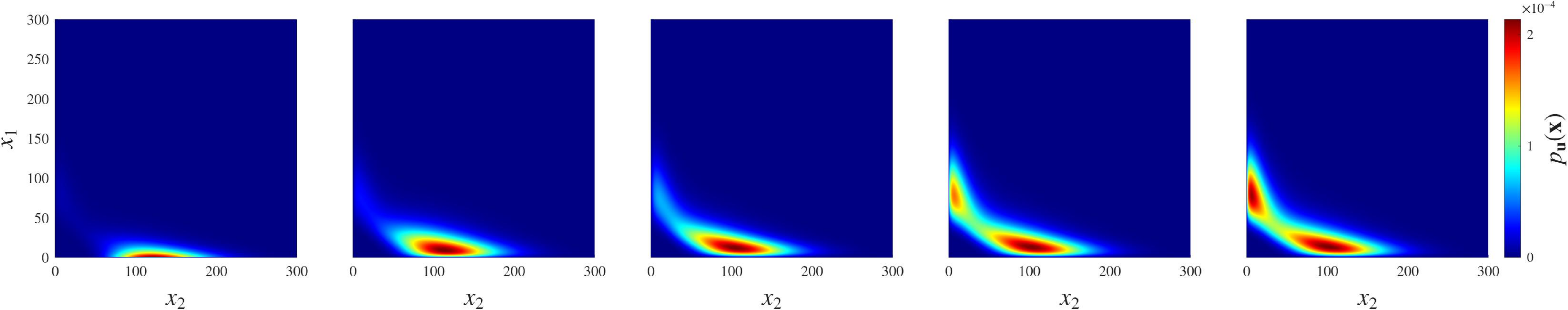}


    \caption{Contractivity analysis for Case Study I. 
    (a) Pairwise $L^1$ distances (logarithmic scale) as functions of time,
    showing monotonic decrease.
    (b) Temporal snapshots of the probability distribution for each initial
    condition (top: IC1; middle: IC2; bottom: IC3).
    }
    \label{fig:contractivty_ts1}
\end{figure}

The robustness of the control is evidenced in Fig.~\ref{fig:contractivty_ts1}a by the precipitous monotonic drop in pairwise $L^1$ distances. Fig.~\ref{fig:contractivty_ts1}b illustrates the distribution profiles at various time points for different initial conditions, where the system demonstrates clear convergence to a common distribution, as anticipated by the contractivity framework developed in Section \ref{sec:contractivity}.


\begin{figure}[!htbp]
    \centering
\begin{tikzpicture}
    
    \node[anchor=north west, inner sep=0] (imgA) at (-0.2, 6.85) {
        \includegraphics[width=0.21\textwidth]{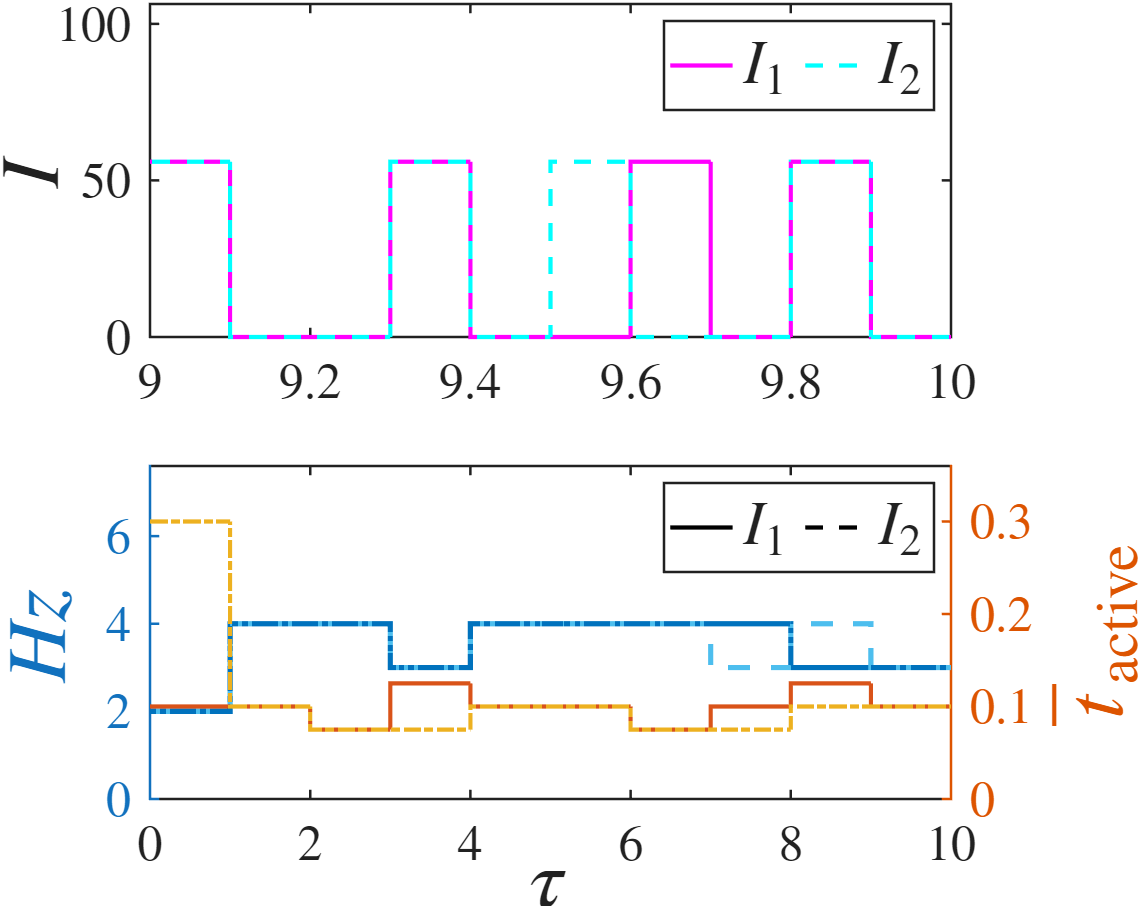}
    };
    \node[anchor=south west, font=\bfseries] at (-0.5cm, 7) {a)};

    \node[anchor=north west, inner sep=0] (imgB) at (0.225\textwidth, 7) {
        \includegraphics[width=0.21\textwidth]{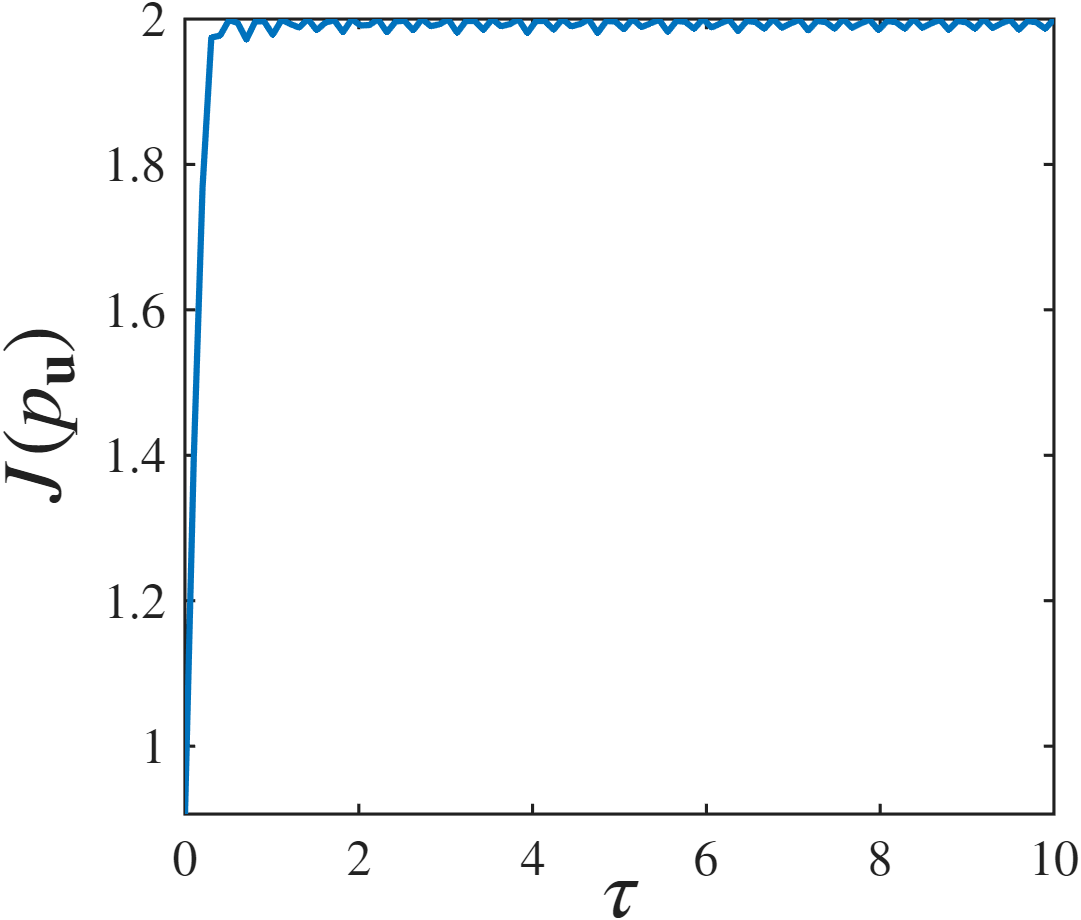}
    };
    \node[anchor=south west, font=\bfseries] at (0.24\textwidth-0.5cm, 7) {b)};

    
    \node[anchor=north west, inner sep=0] (imgC1) at (0, 3) {
        \includegraphics[width=0.21\textwidth]{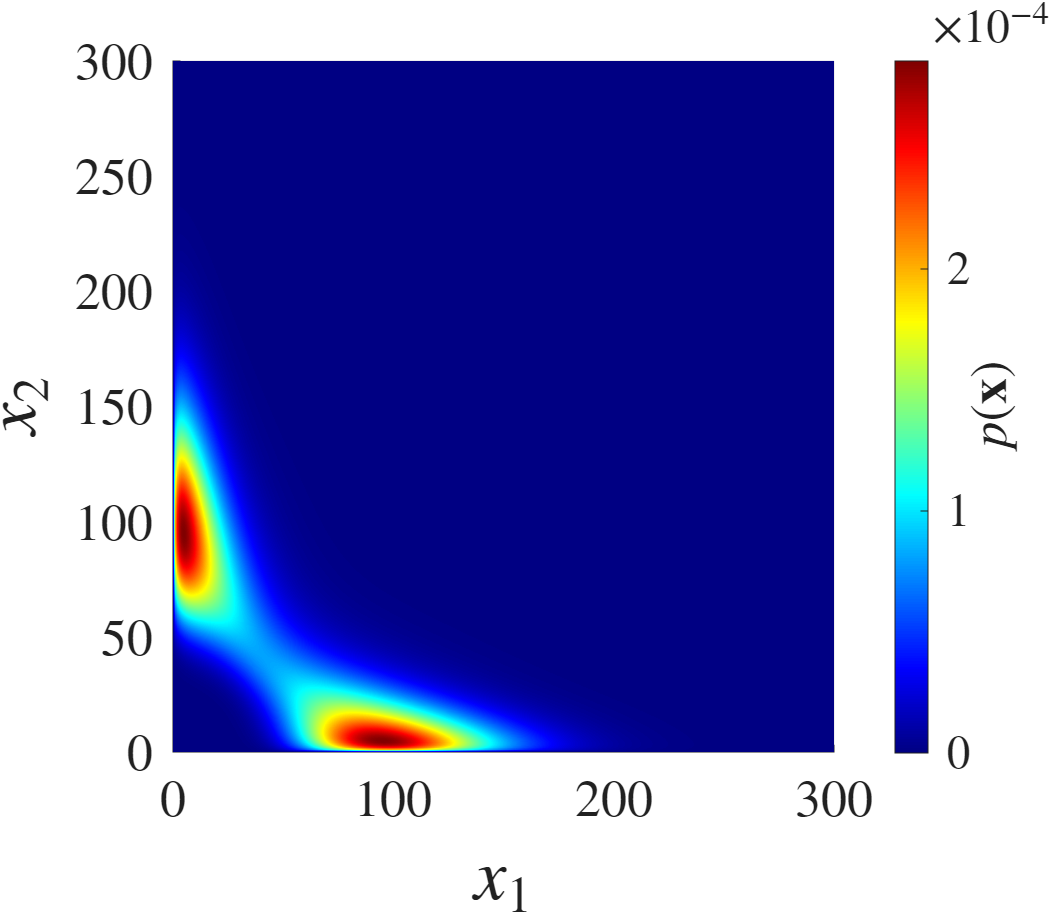}
    };
    \node[anchor=south west, font=\bfseries] at (-0.5cm, 3) {c)};

    \node[anchor=north west, inner sep=0] (imgC2) at (0.23\textwidth, 3) {
        \includegraphics[width=0.21\textwidth]{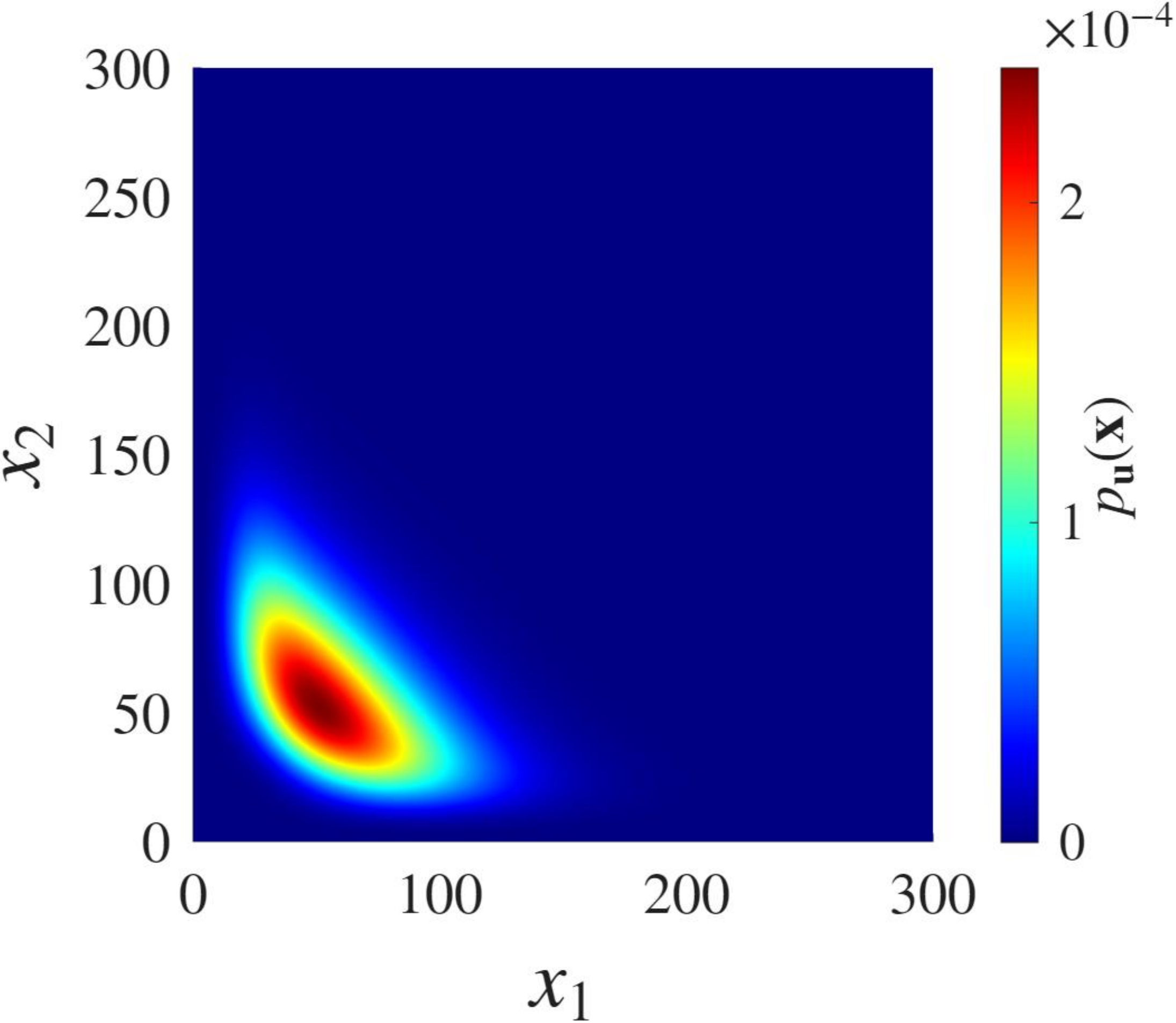}
    };
    \node[anchor=south west, font=\bfseries] at (0.24\textwidth-0.5cm, 3) {d)};

    \node[anchor=north west, inner sep=0] (imgC3) at (-0.2, -0.75) {
        \includegraphics[width=0.2\textwidth]{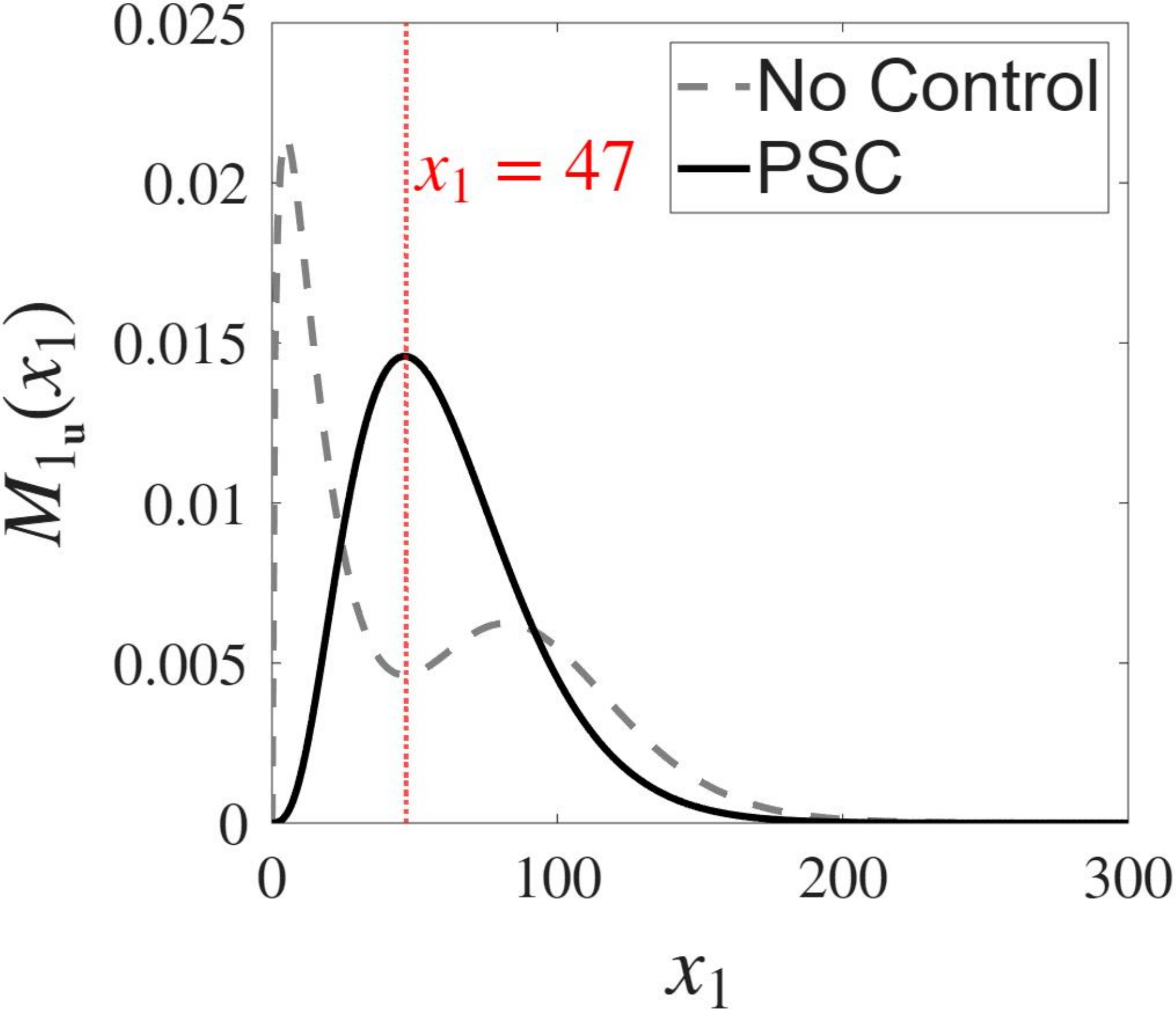}
    };
    \node[anchor=south west, font=\bfseries] at (-0.5cm, -0.75) {e)};

    \node[anchor=north west, inner sep=0] (imgC4) at (0.22\textwidth, -0.75) {
        \includegraphics[width=0.2\textwidth]{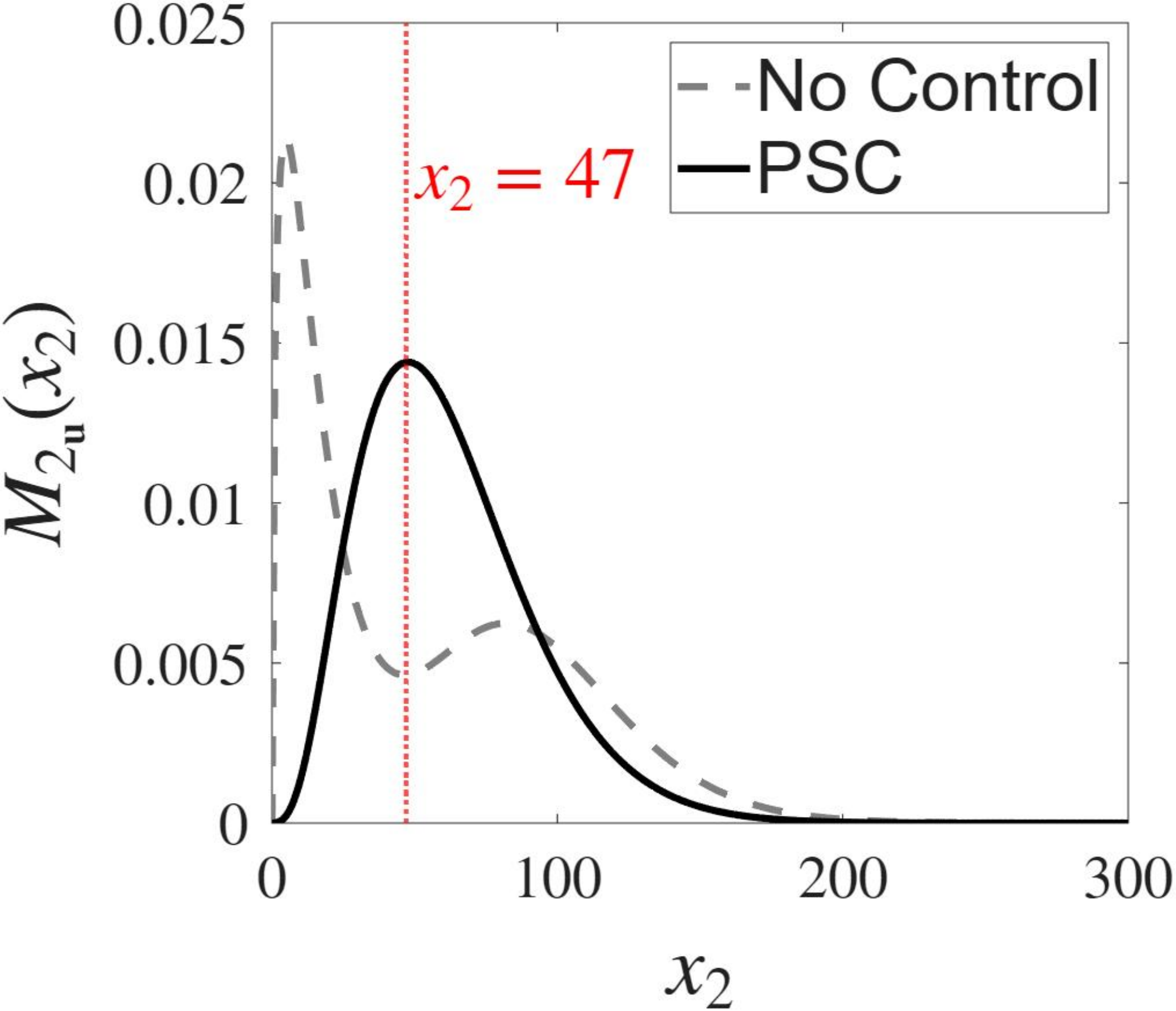}
    };

\end{tikzpicture}

    \caption{
(a) Input signals. Top: final interval. Bottom:
activation frequency (blue) and mean activation time (orange) for $I_1$ (solid) and $I_2$ (dashed).
(b) Evolution of the cost functional.
(c) Uncontrolled stationary distribution. 
(d) Joint distribution concentrated at the previously low-probability region.
(e) Marginal distributions $M_{1_{\mathbf{u}}}$ and $M_{2_{\mathbf{u}}}$ (solid:controlled, dashed: uncontrolled).}
    \label{fig:Results2}
\end{figure}
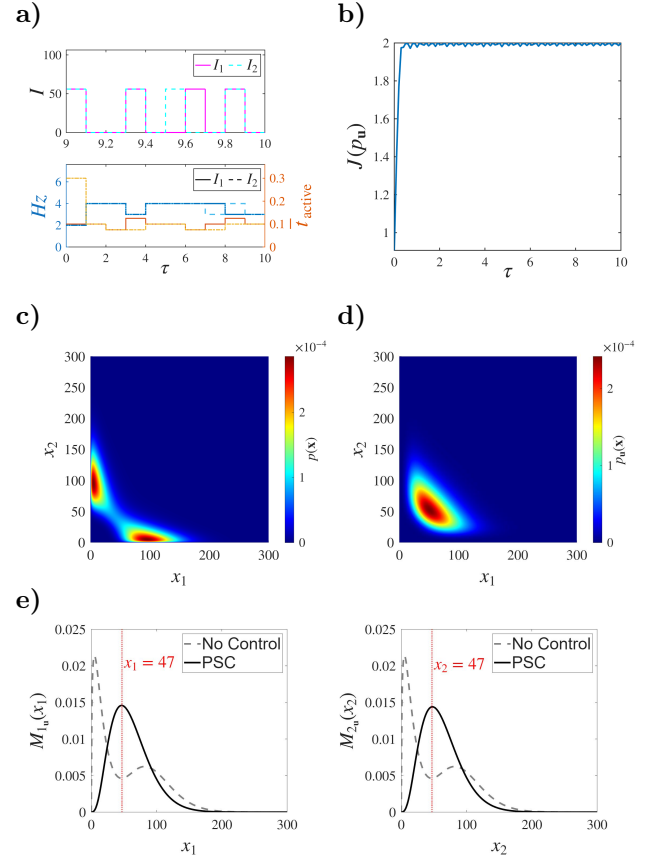

\textbf{Case Study II: Stabilization of a Low-Probability Region
in Two Dimensions.}
Secondly, a symmetric toggle switch is considered, consisting of two mutually repressing genes with identical parameters:  $k_m^1 = k_m^2 = 10$,
    $k_x^1 = k_x^2 = 100$, $\gamma_m^1 = \gamma_m^2 = 10$,
    $\gamma_x^1 = \gamma_x^2 = 1$, $K_{ij} = 40$, $H_{ij} = 4$, $\theta_{I_i} = 0.1$,
$\mu_{I_i} = 2$, $\varepsilon_i = 0.1$, $\Delta t_k = 0.005$, and
$\Delta x_i = 1$. The uncontrolled system admits a bimodal stationary distribution with two dominant modes, corresponding to high expression of either protein, separated by a low-probability region (Fig.~\ref{fig:Results2}c).

The control objective is to balance the switch,  a benchmark for gene network
control \citep{Lugagne2017, Brancato2023, Guarino2020Balancing,
Fernandez2022, Fernandez2025, Vaghy2024}, which in the stochastic regime boils down to concentrate the probability mass in the
intermediate region. The target state
$\mathbf{x}^* = (x_1^*, x_2^*)$ is defined such that each component
$x_i^*$ corresponds to a local minimum of the uncontrolled marginal
distribution $M_i(x_i)$. The control input is $\mathbf{u}(t) = [I_1(t), I_2(t)]$. The
PIDE~\eqref{eq:PIDE_general} is solved over $\Omega = [0, 300]^2$. The marginal distributions are obtained by integrating the joint density over the complementary variable: $
    M_{i_{\mathbf{u}}}(t, x_i) =
    \int_{0}^{300} p_{\mathbf{u}}(t, \mathbf{x})\,\mathrm{d}x_j$,
and normalized by its maximum value at each time instant:
$
    \widetilde{M}_{i_{\mathbf{u}}}(t, x_i) =
    M_{i_{\mathbf{u}}}(t, x_i)/\max_{x_i} M_{i_{\mathbf{u}}}(t, x_i)$.
The cost functional is defined as
\begin{equation*}
    J(p_{\mathbf{u}}(t, \mathbf{x})) =
    \widetilde{M}_{1_{\mathbf{u}}}(t, x_1^*) +
    \widetilde{M}_{2_{\mathbf{u}}}(t, x_2^*),
\end{equation*}
with $J \in (0, 2]$;  $J = 2$ is attained when each marginal
achieves its maximum at the corresponding target coordinate $x_i^*$. The inducer saturation levels are $\kappa_i = 56$, corresponding to
$\alpha_i = 0.01$ in~\eqref{eq:kappa_def}. Algorithm~\ref{alg:optimal_inductor_selection} is applied over the
configuration set $S \in \{0,1\}^{2 \times 2}$. The actuation window is $w = 20$, giving $\Delta t_m = 20\,\Delta t_k$.
The results are reported in Fig.~\ref{fig:Results2}. The PSC induces a
symmetric switching pattern for $I_1$ and $I_2$ (Fig.~\ref{fig:Results2}a). The cost functional converges to $J = 2$ (Fig.~\ref{fig:Results2}b), indicating that both marginal distributions achieve their maximum at the respective target coordinates. 
The controlled stationary distribution is concentrated at $\mathbf{x}^*$ (Fig.~\ref{fig:Results2}d,e). The closed-loop simulation time is $9.92$~s.

\begin{figure}[!htbp]
    \centering

    \makebox[\textwidth][l]{\textbf{a)}}\\[-0.3ex]
    \includegraphics[width=0.4\textwidth]{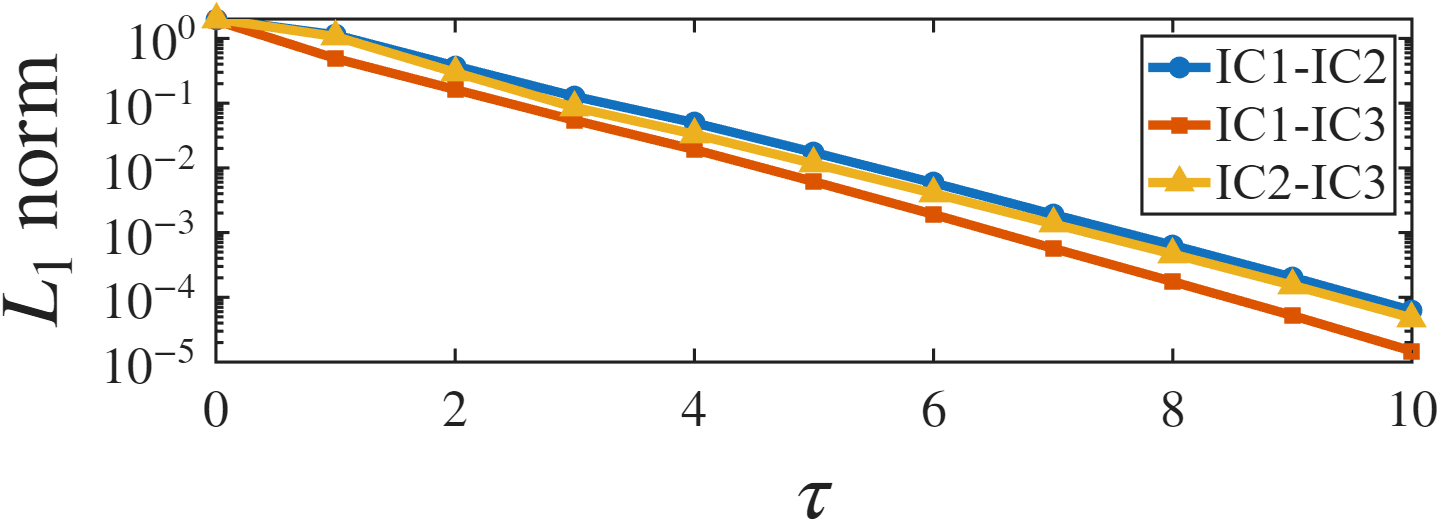}

    \makebox[\textwidth][l]{\textbf{b)}}\\[-0.3ex]
    \includegraphics[width=0.45\textwidth]{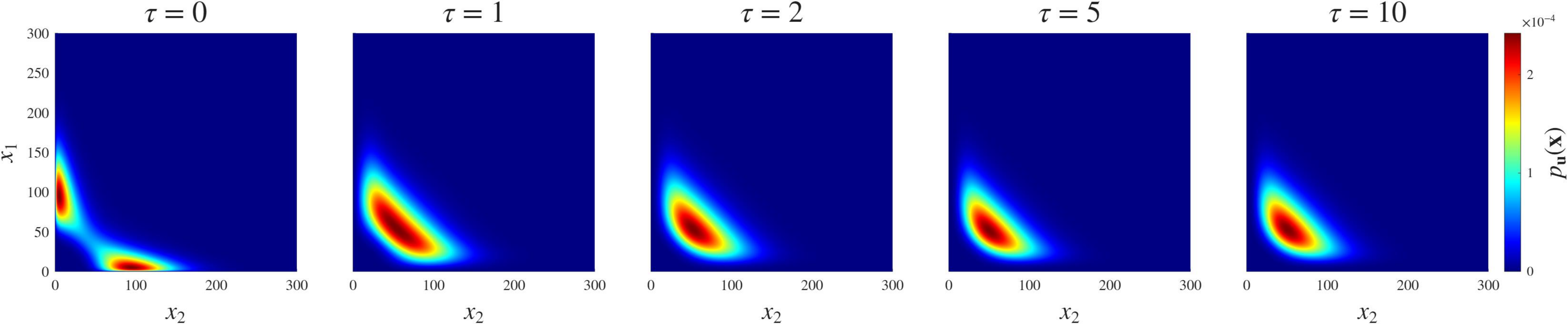}\\
    \includegraphics[width=0.45\textwidth]{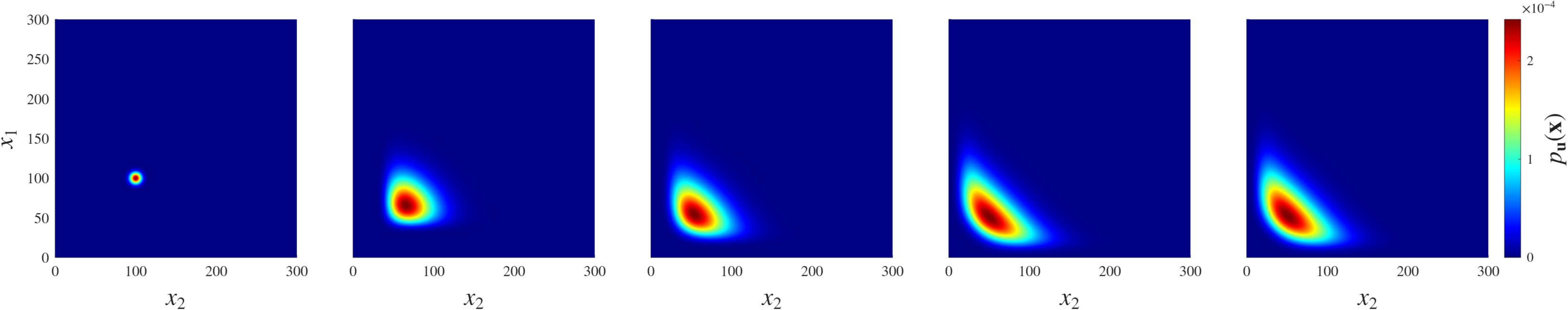}\\
    \includegraphics[width=0.45\textwidth]{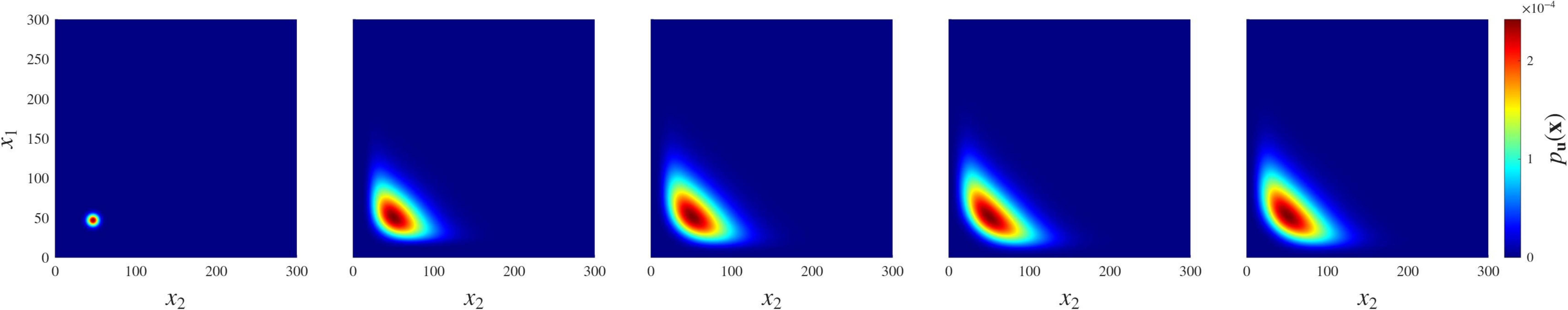}


    \caption{Contractivity analysis for Case Study II. 
    (a) Pairwise $L^1$ distances (logarithmic scale) as functions of time,
    showing monotonic decrease.
    (b) Temporal snapshots of the probability distribution for each initial
    condition (top: IC1; middle: IC2; bottom: IC3).
    }
    \label{fig:contractivty_ts2}
\end{figure}

Fig.~\ref{fig:contractivty_ts2}a shows the pairwise $L^1$ distances, which decrease monotonically. Temporal snapshots of the distributions for various initial conditions are presented in Fig.~\ref{fig:contractivty_ts2}b, showing convergence to a unique state as predicted.

\begin{figure}[!htbp]
    \centering
    \includegraphics[width=0.15\textwidth]{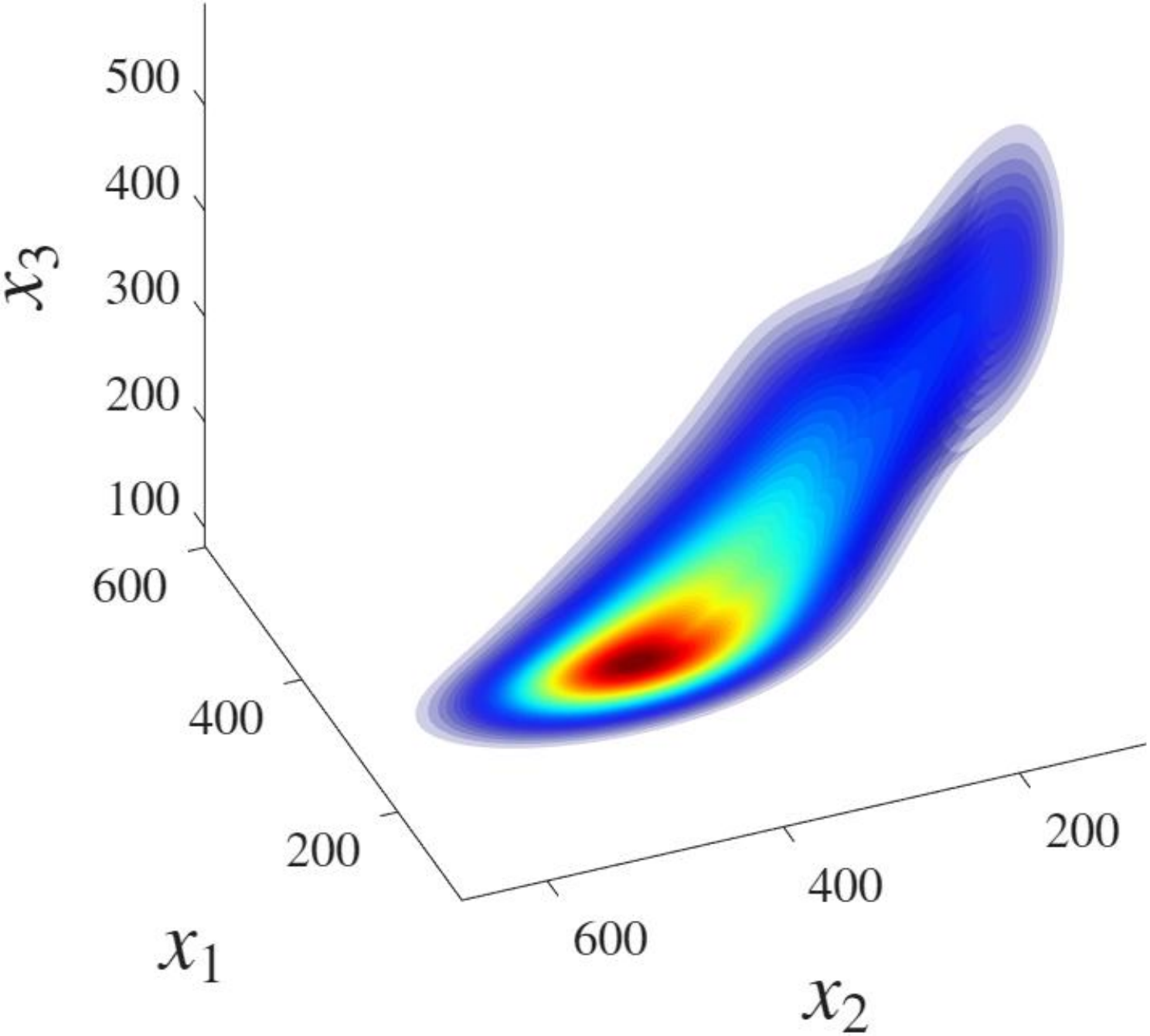}
    \includegraphics[width=0.15\textwidth]{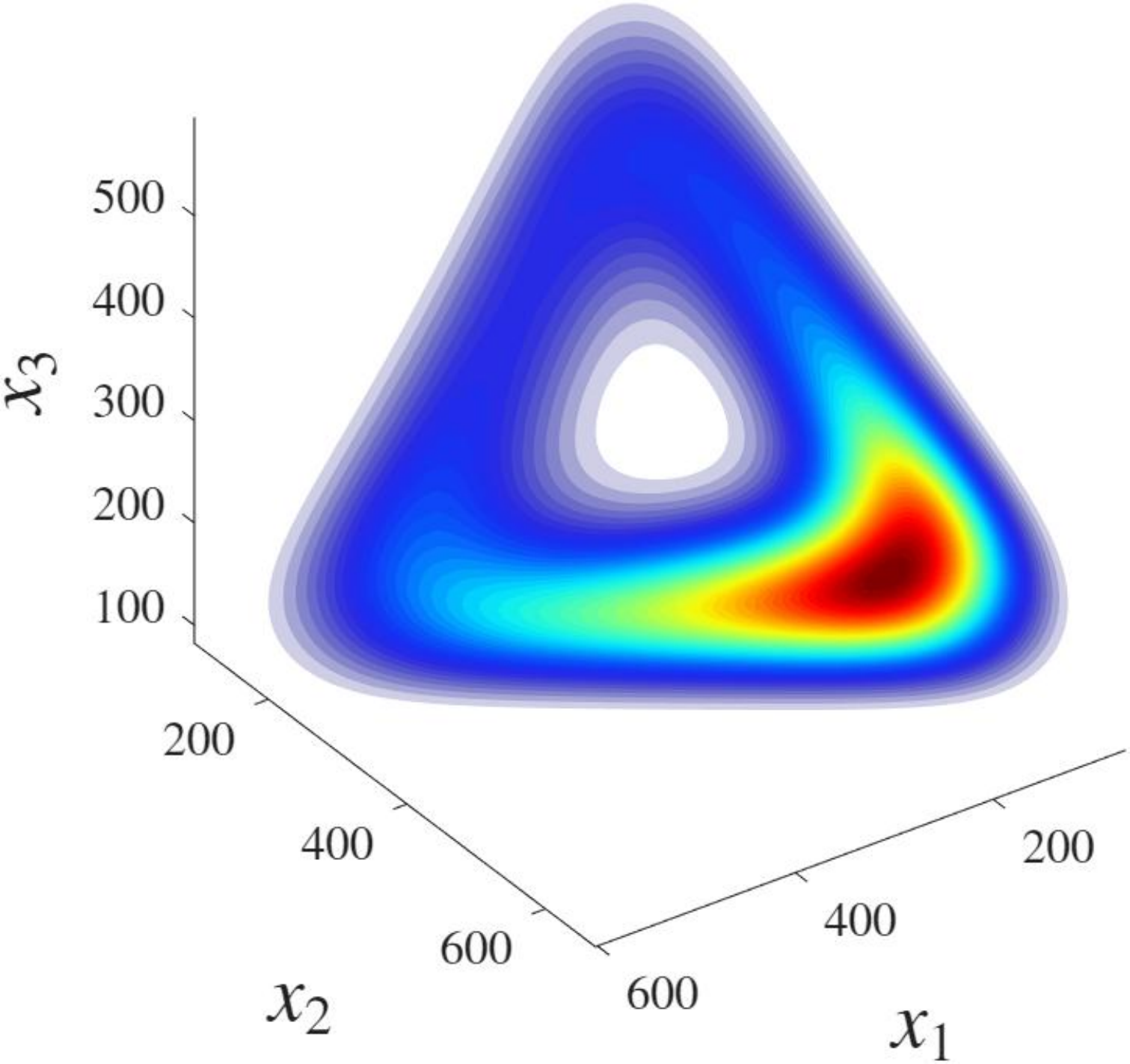}\hspace{7pt}
    \includegraphics[width=0.15\textwidth]{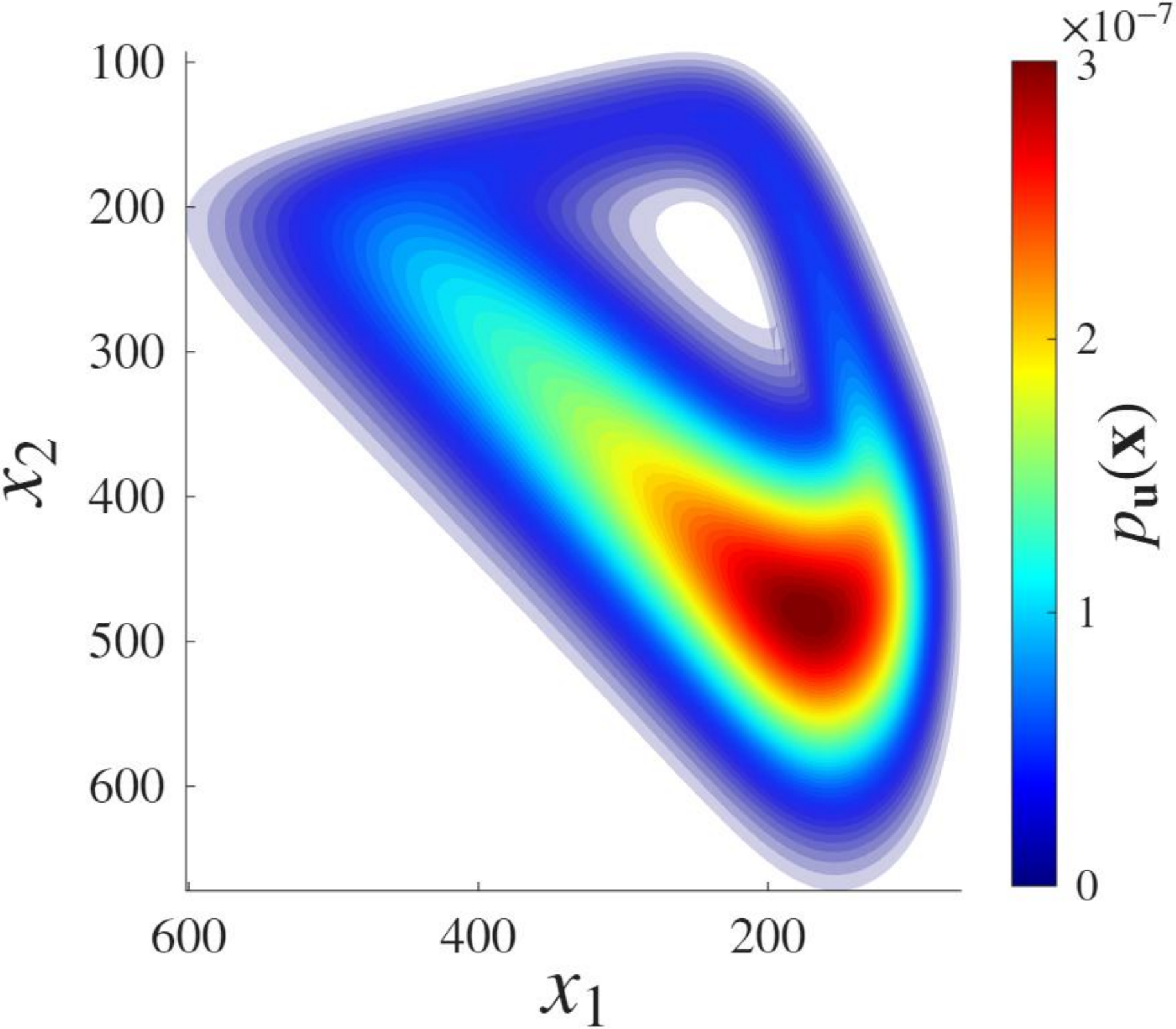}
    \caption{Stationary probability distribution of the uncontrolled
    three-protein genetic oscillator, shown from lateral, frontal, and top views. 
    Parameters: $k_m^1 = 125$, $k_m^2 = 100$, $k_m^3 = 115$;
    $k_x^1 = 90$, $k_x^2 = 110$, $k_x^3 = 100$;
    $\gamma_m^i = 17.6822$ and $\gamma_x^i = 1$ for $i \in \{1,2,3\}$.}
    \label{fig:case3_oscillator}
\end{figure}

\textbf{Case Study III: Stabilization of a Low-Probability Region
in Three Dimensions.}
A stochastic three-protein genetic oscillator is considered, consisting of three species $x_1$, $x_2$, $x_3$ coupled in a cyclic inhibitory network: $x_1 \dashv x_3$, $x_3 \dashv x_2$, $x_2 \dashv x_1$. The uncontrolled system admits a stationary distribution supported on a ring-shaped manifold in $(x_1, x_2, x_3)$ space, with low probability density at the centre (Fig.~\ref{fig:case3_oscillator}).

\begin{figure}[!htbp]
\centering
    \begin{tikzpicture}
        \node[anchor=north west, inner sep=0] (imgA) at (-0.2, 9.9) {
            \includegraphics[width=0.21\textwidth]{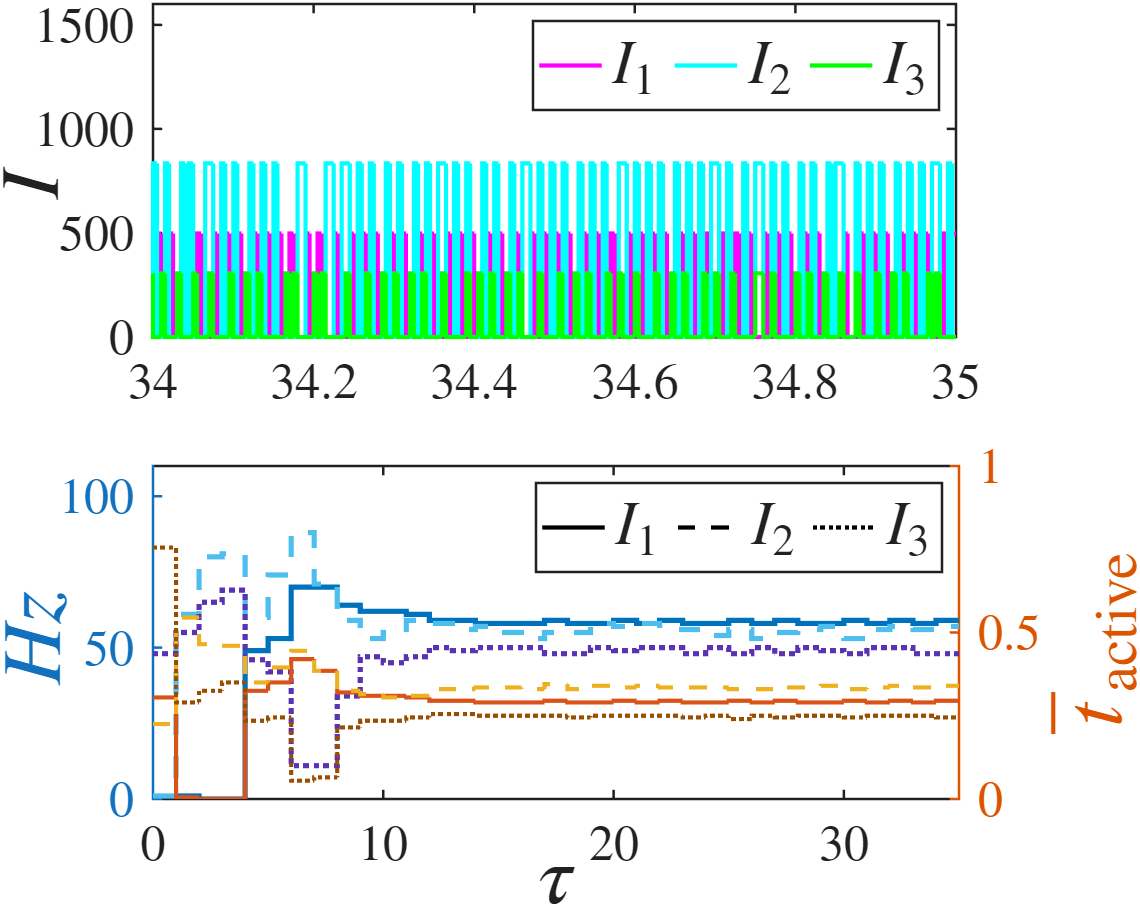}
        };
        \node[anchor=south west, font=\bfseries] at (-0.5cm, 10) {a)};

        \node[anchor=north west, inner sep=0] (imgB) at (0.225\textwidth, 10) {
            \includegraphics[width=0.21\textwidth]{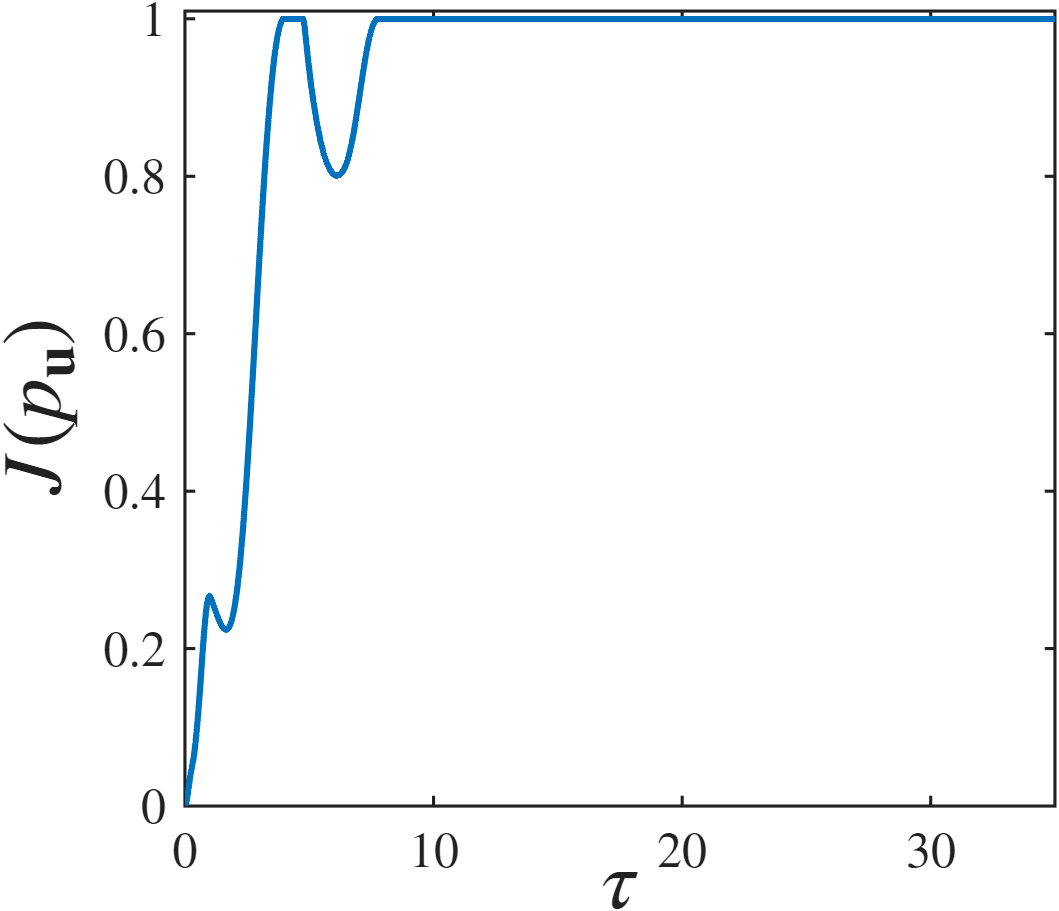}
        };
        \node[anchor=south west, font=\bfseries] at (0.24\textwidth-0.5cm, 10) {b)};

        \node[anchor=north west, inner sep=0] (imgC) at (-0.2, 5.9) {
            \includegraphics[width=0.215\textwidth]{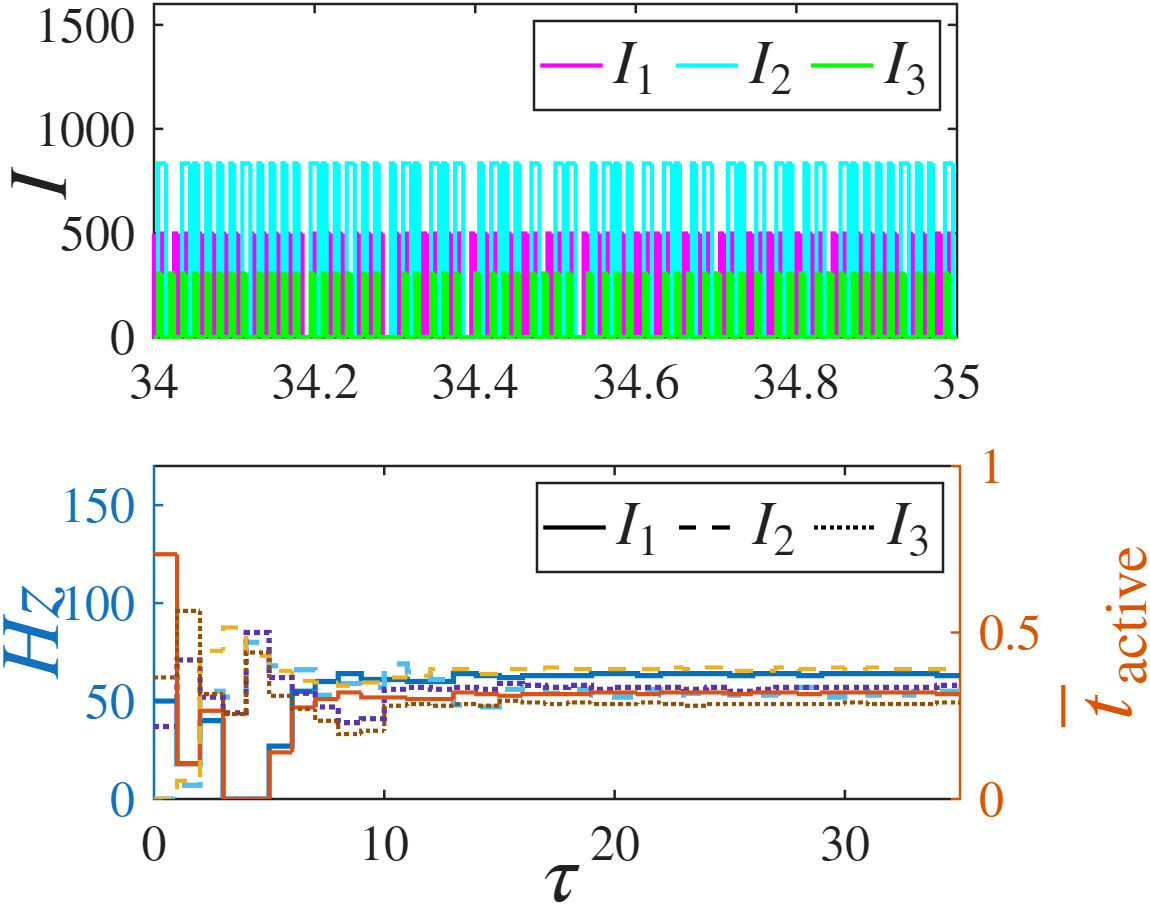}
        };
        \node[anchor=south west, font=\bfseries] at (-0.5cm, 6) {c)};

        \node[anchor=north west, inner sep=0] (imgD) at (0.225\textwidth, 6) {
            \includegraphics[width=0.215\textwidth]{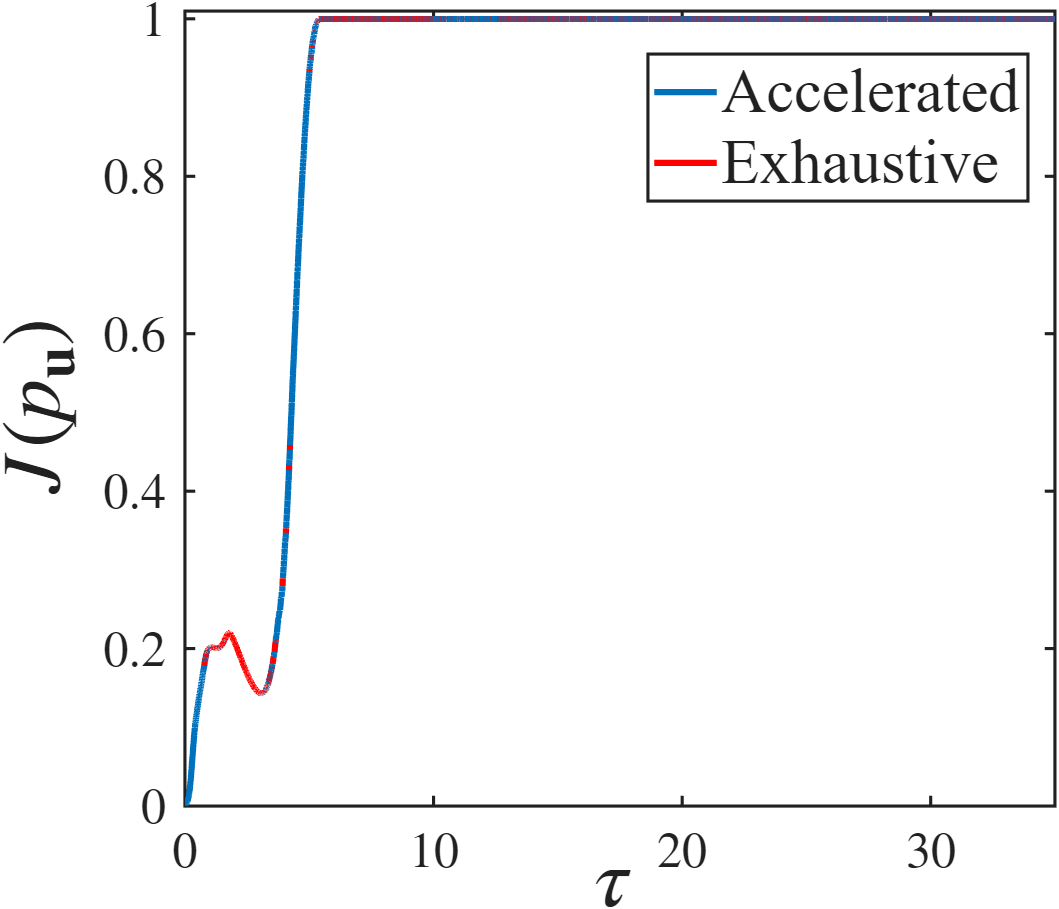}
        };
        \node[anchor=south west, font=\bfseries] at (0.24\textwidth-0.5cm, 6) {d)};

        \node[anchor=north west, inner sep=0] (imgE1) at (-0.2, 2) {
            \includegraphics[width=0.15\textwidth]{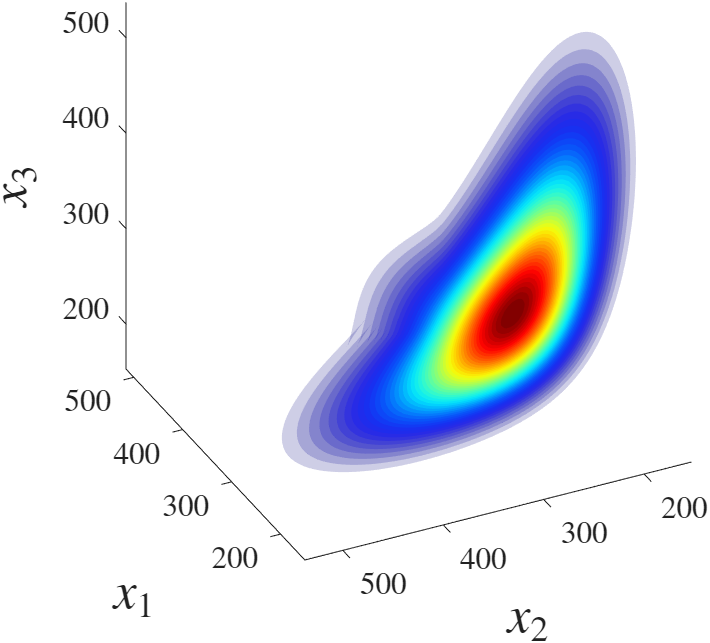}
        };
        \node[anchor=north west, inner sep=0] (imgE2) at (0.14\textwidth, 2) {
            \includegraphics[width=0.15\textwidth]{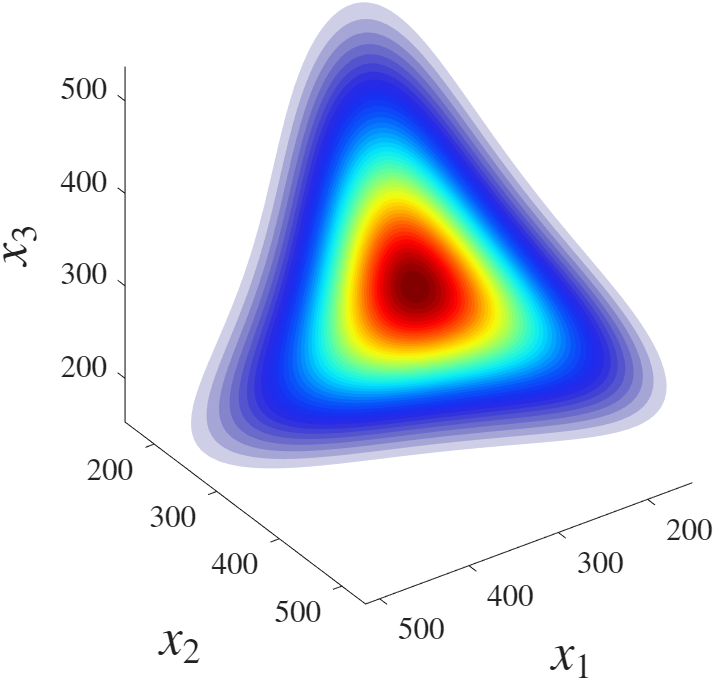}
        };
        \node[anchor=north west, inner sep=0] (imgE3) at (0.3\textwidth, 2) {
            \includegraphics[width=0.15\textwidth]{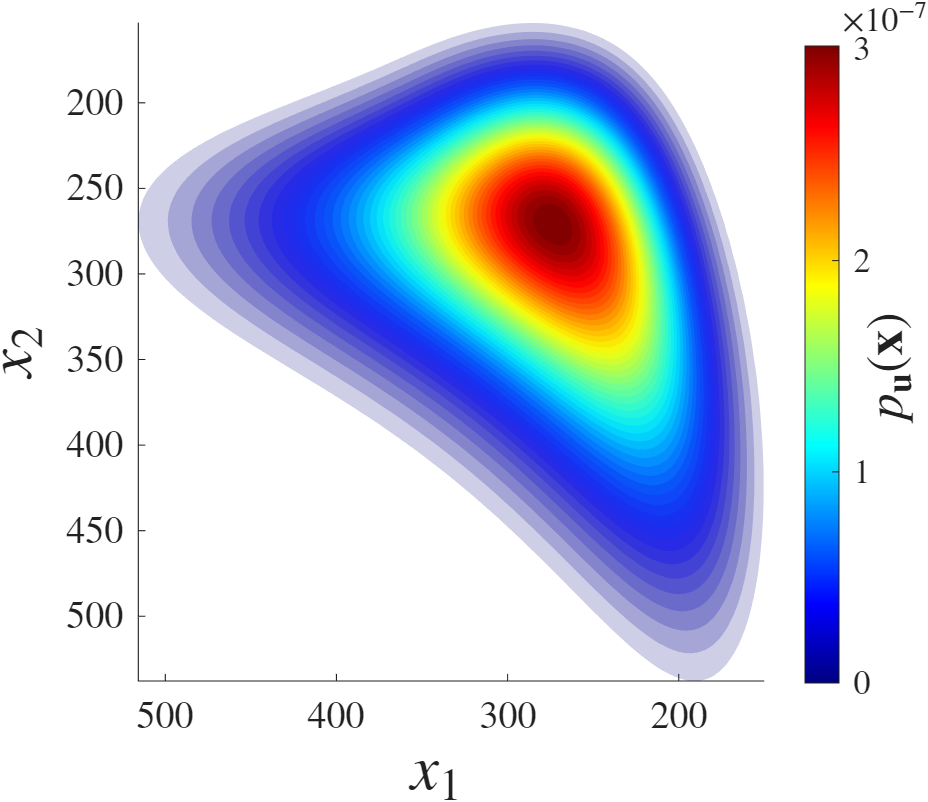}
        };
        \node[anchor=south west, font=\bfseries] at (-0.5cm, 2) {e)};

        \node[anchor=north west, inner sep=0] (imgF1) at (-0.2, -1.25) {
            \includegraphics[width=0.15\textwidth]{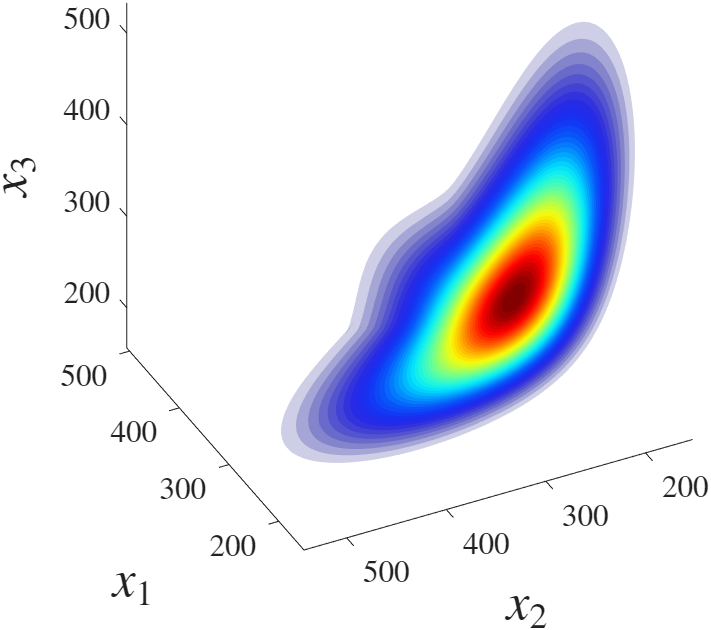}
        };
        \node[anchor=north west, inner sep=0] (imgF2) at (0.14\textwidth, -1.25) {
            \includegraphics[width=0.15\textwidth]{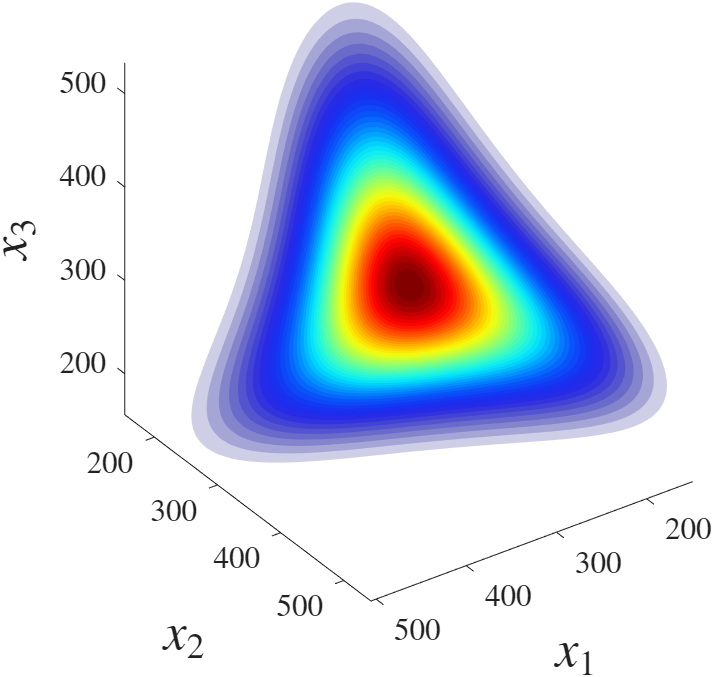}
        };
        \node[anchor=north west, inner sep=0] (imgF3) at (0.3\textwidth, -1.25) {
            \includegraphics[width=0.15\textwidth]{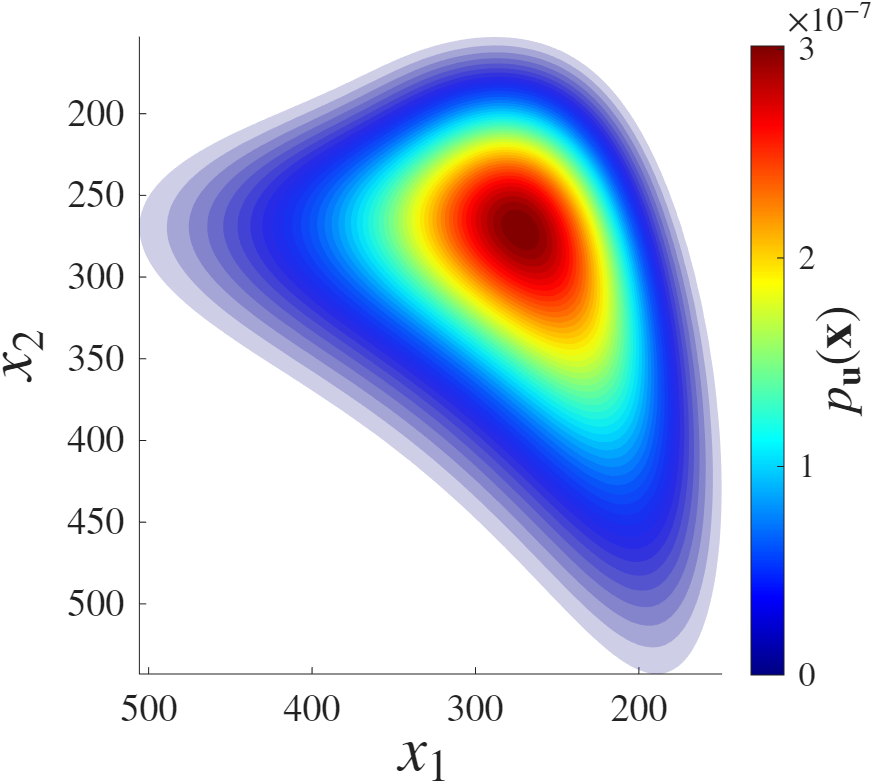}
        };
        \node[anchor=south west, font=\bfseries] at (-0.5cm, -1.25) {f)};

    \end{tikzpicture}

\caption{Results for exhaustive and accelerated PSC.
(a,c) Input signals. Top:  final interval. Bottom: activation frequency (blue) and mean activation time (orange) for $I_1$ (solid), $I_2$ (dashed), and $I_3$ (dotted). (b,d) Evolution of the cost functional. In (d), accepted by the NN (red) and fallback decisions (blue). (e,f) Controlled joint distributions for exhaustive and accelerated PSC, respectively,  from different perspectives.}
    \label{fig:Results3}
\end{figure}
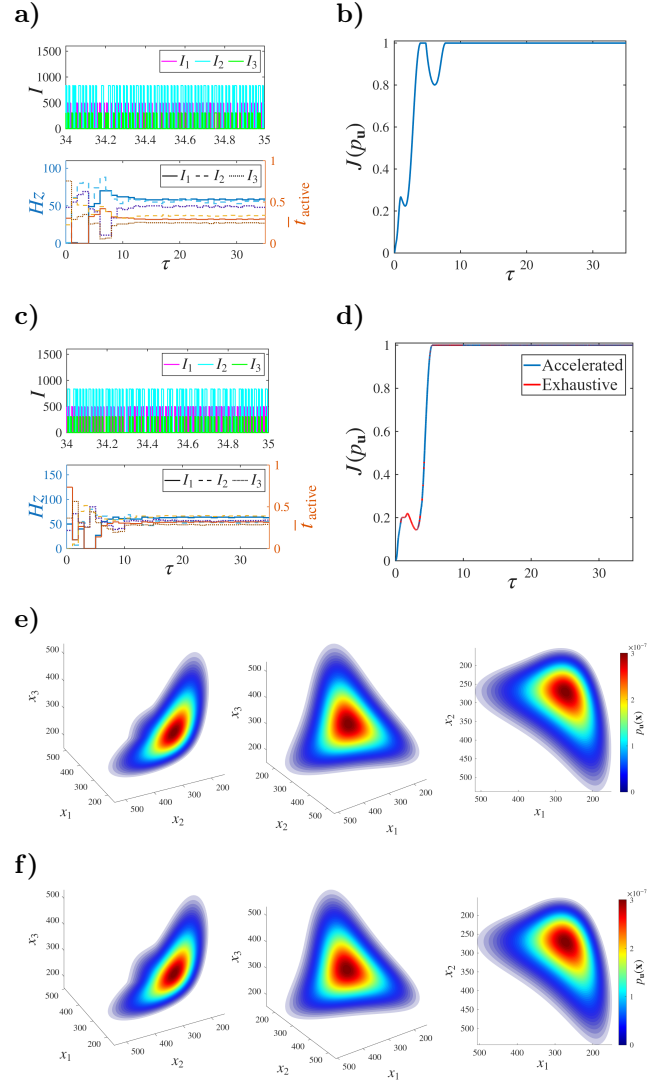

The control objective is to stabilize the system at the low-probability
centre of the stationary distribution. The control input is
$\mathbf{u}(t) = [I_1(t), I_2(t), I_3(t)]$, where each inducer $I_i$ modulates the repression of protein $x_i$ in~\eqref{eq:PIDE_general}. The regulatory interactions follow a cyclic structure with $j = i + 1\ (\mathrm{mod}\ 3)$. The parameters are $K_{ij} = 200$,  $H_{12} = 8$, $H_{23} = 9$, $H_{31} = 7$,  $\theta_{I_1} = 0.08$, $\theta_{I_2} = 0.06$, $\theta_{I_3} = 0.11$, with $\mu_{I_i} = 2$ and  $\varepsilon_i = 0.15$. The domain is $\Omega = [0, 1000]^3$, discretized with $\Delta x_i = 4$ and $\Delta t_k = 0.005$. The cost functional is defined as 
\begin{equation}
    J(p_{\mathbf{u}}(t, \mathbf{x})) =
    \frac{p_{\mathbf{u}}(t, \mathbf{x}^*)}{\max_{\mathbf{x}}\,
    p_{\mathbf{u}}(t, \mathbf{x})}
\end{equation}   
where $\mathbf{x}^*$ denotes the geometric centre of the oscillator. The functional satisfies $J \in (0, 1]$, with $J = 1$ attained when the
probability maximum coincides with $\mathbf{x}^*$. The inducer saturation levels are $\boldsymbol{\kappa} = [497.5,\ 834.3,\ 305.9]$, corresponding to $\alpha = 0.01$ in~\eqref{eq:kappa_def}. The configuration set covers all $2^3 = 8$ binary combinations, $S \in \{0,1\}^3$. The actuation window is set to $w = 1$, so that $\Delta t_m = \Delta t_k$ and the control action is re-evaluated at every integration step. The exhaustive and accelerated PSC (Algorithm~\ref{alg:psc_accelerated}) performances are shown in Figs.~\ref{fig:Results3}b,c,e,f) and compared in terms of the iterations needed to attain $J = 1$ among other metrics summarized in Table~\ref{tab:psc_comparison}. The accelerated scheme reduces the
number of PIDE evaluations by 64\% and the execution time by 57\%. 

\begin{table}[H]
\caption{Metrics for exhaustive and accelerated PSC }
\centering
\begin{tabular}{lcc}
\hline
\textbf{Metric} & \textbf{PSC} & \textbf{Accelerated PSC} \\
\hline
Iterations  & 1545                   & 1076                   \\
Elapsed time  [s]        & $5.932 \times 10^3$   & $2.570 \times 10^3$   \\
NN acceptances        & N/A                   & 596                  \\
PIDE evaluations            & 12360                  & 4436                  \\
\hline
\end{tabular}
\label{tab:psc_comparison}
\end{table}

\begin{figure}[!htbp]
    \centering

    \makebox[\textwidth][l]{\textbf{a)}}\\[-0.3ex]
    \includegraphics[width=0.4\textwidth]{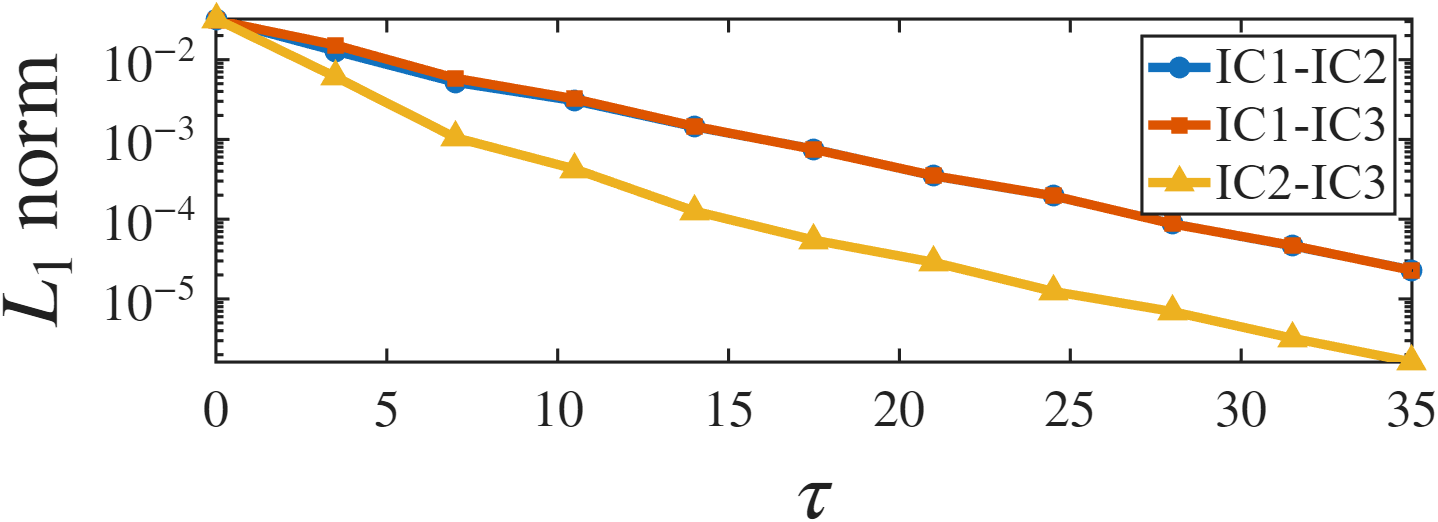}

    \makebox[\textwidth][l]{\textbf{b)}}\\[-0.3ex]
    \includegraphics[width=0.45\textwidth]{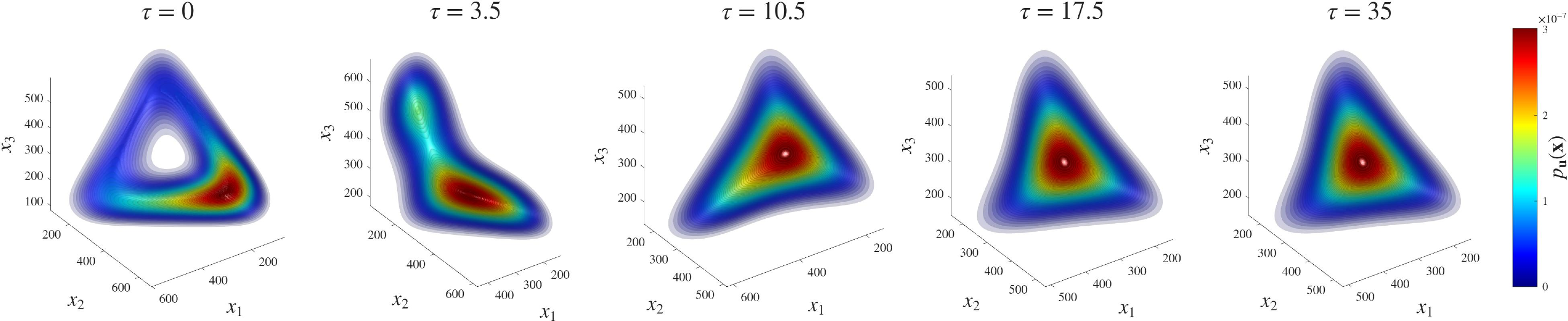}\\
    \includegraphics[width=0.45\textwidth]{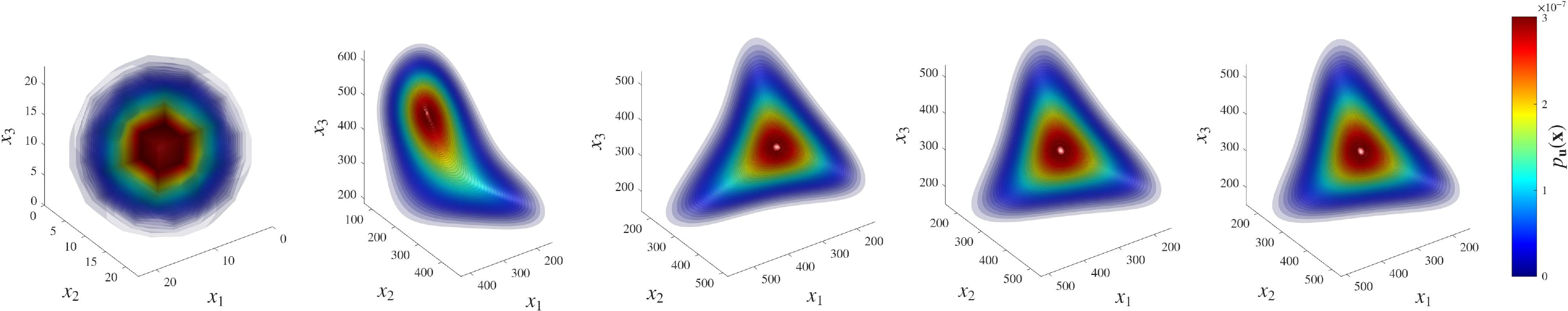}\\
    \includegraphics[width=0.45\textwidth]{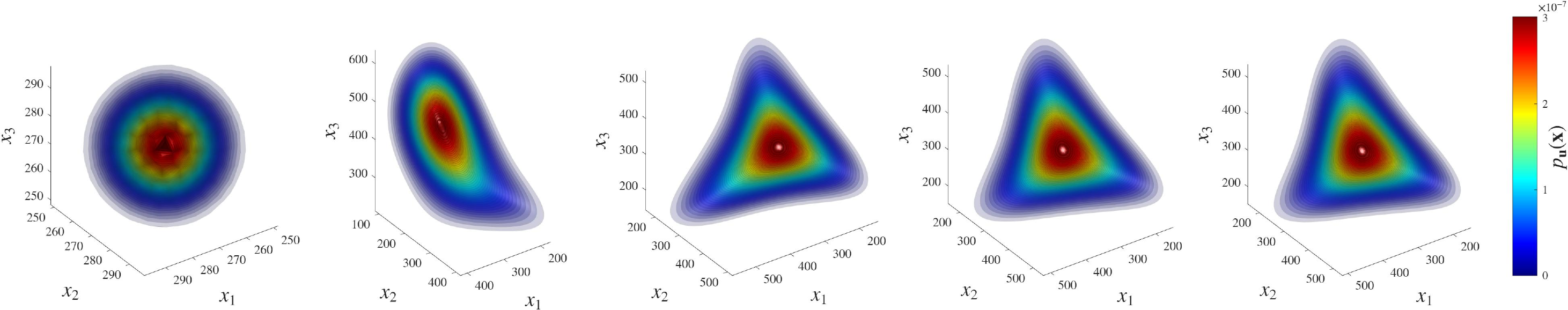}


    \caption{Contractivity analysis for Case Study III. 
    (a) Pairwise $L^1$ distances (logarithmic scale) as functions of time,
    showing monotonic decrease.
    (b) Temporal snapshots of the probability distribution for each initial
    condition (top: IC1; middle: IC2; bottom: IC3).
    }
    \label{fig:contractivty_osc}
\end{figure}

The pairwise $L^1$ distances shown in Fig.~\ref{fig:contractivty_osc}a decrease monotonically over time. Fig.~\ref{fig:contractivty_osc}b displays temporal snapshots of the distributions for three initial conditions, all of which converge to the same steady-state profile. This behavior remains consistent with the theoretical predictions.



\section{Conclusions}
The proposed Predictive Switching Control (PSC) algorithm provides a model-based approach for controlling distributed stochastic dynamics while maintaining computational tractability and enabling direct shaping of full probability distributions. Unlike moment-based or trajectory-level methods, PSC operates over a discrete input set, avoiding gain tuning and reducing the complexity of continuous model predictive control, while outperforming standard ON–OFF strategies through global input selection at each step. A hybrid implementation incorporating a neural-network-based proposal mechanism further enhances scalability, significantly reducing model evaluations without compromising performance. From a theoretical perspective, the proposed control framework admits a contraction-based analysis, guaranteeing $L^1$ contractivity and yielding geometric convergence under strictly positive leakage terms. This provides robustness with respect to initial conditions and enables the use of a single control policy across diverse initial distributions. Overall, PSC offers a structured, stable, and computationally efficient approach for controlling stochastic gene regulatory networks, particularly in settings where dynamic complexity requires going beyond low-order statistical descriptions.

\section{Acknowledgements}
CF and IOM acknowledge ERC Consolidator Grant CellWise
ERC-2024-CoG-101170783. GSz acknowledges the support the Hungarian National Research, Development, and Innovation Office (NKFIH) through the grant K-145934.
The work of MP was supported in part by the Spanish Ministry of Science and Innovation under Grant PID2022-141058OB-I00 and Grant PID2023-146275NB-C21, in part by MI-CIU/AEI/10.13039/501100011033, and in part by ERDF/EU. 

\appendix

\section{Neural acceleration module for PSC}
\label{app:nn_psc}

\textbf{Action selection formulation.}
At each switching instant $t_m$, PSC selects one element from the finite admissible set
\begin{equation*}
    \mathcal{U} = \{\boldsymbol{\kappa} \odot \mathbf{S}_r \,:\, r = 1, \dots, 2^n\},
\end{equation*}
where $\boldsymbol{\kappa}$ contains the inducer saturation levels and $\mathbf{S}_r \in \{0,1\}^n$ is the $r$-th binary switching pattern. The neural module approximates the mapping
\begin{equation*}
    \mathbf{z}(t_m) \mapsto \hat{\mathbf{s}} \in [0,1]^n,
\end{equation*}
with $\mathbf{z}(t_m)$ a compact feature vector describing the current probabilistic state and recent control history. The output $\hat{\mathbf{s}}$ is rounded to obtain the binary vector $\mathbf{S}_{r^*} = \mathrm{round}(\hat{\mathbf{s}})$, hence the candidate input $\mathbf{u}^{(r^*)} = \boldsymbol{\kappa} \odot \mathbf{S}_{r^*}$. If the PIDE simulation verifies that $J_{r^*} \ge J_{t_m}$, the exhaustive search is skipped; otherwise, the controller falls back to the standard PSC evaluation (Algorithm~\ref{alg:psc_accelerated}).

\textbf{Input representation.}
The network input is defined as
\begin{equation*}
    \mathbf{z}(t_m) =
    \bigl(
    \mathbf{s}_{m-1},\;
    \mathbf{x}_{\mathrm{mode}},\;
    p^*,\;
    d^*,\;
    D_{\mathrm{KL}}(p_{\mathbf{u}}(t_m,\mathbf{x}) \| p^*)
    \bigr),
    \label{eq:nn_features_appendix}
\end{equation*}
comprising: the previous control action $\mathbf{s}_{m-1}$; the modal
location $\mathbf{x}_{\mathrm{mode}}$ of the current distribution; the
probability mass at the target $p^*$; the geometric distance to the target
$d^*$; and the Kullback--Leibler divergence
$D_{\mathrm{KL}}(p_{\mathbf{u}}(t_m,\mathbf{x}) \| p^*)$ measuring the
distributional mismatch with respect to the target. his representation avoids processing the full PDF and provides a compact description of the control state.

\textbf{Dataset construction and partitioning.}
The training dataset comprises $5 \times 10^3$ samples generated offline
from exhaustive PSC simulations over a diverse set of trajectories, target
configurations, and symmetry conditions. Each sample pairs the feature
vector $\mathbf{z}(t_m)$ with the optimal binary action selected by the
exhaustive policy. A hold-out test set comprising 15\% of the total data is separated prior
to any preprocessing or training. This subset is reserved exclusively for
final performance evaluation, providing an unbiased estimate of
generalization. 

\textbf{Preprocessing.}
Input features are standardized using Z-score normalization. The mean $\mu_j$ and standard deviation $\sigma_j$ are computed from the training split and then applied unchanged to validation and test data. Targets are not scaled, since the output is projected onto $\{0,1\}^n$ by rounding the network output.

\textbf{Network architecture.}
The neural accelerator is a shallow feedforward network with two hidden
layers of sizes $[20,\,10]$. Hidden layers use hyperbolic tangent
(\texttt{tansig}) activation functions. The output layer employs a
symmetric saturating linear transfer function (\texttt{satlins}), which
facilitates gradient stability during Levenberg--Marquardt optimization.
Continuous outputs are projected onto $[0,1]$ via a min-max operator prior
to rounding. The architecture comprises 483 trainable parameters against approximately
4250 training samples, yielding a samples-to-parameters ratio exceeding
$8\!:\!1$. The projected output is constrained to [0,1] by construction.  The $[20,\,10]$ configuration was adopted and evaluated by 5-fold cross-validation.

\textbf{Training Procedure and Cross-Validation.}
Training is performed offline using the Levenberg--Marquardt algorithm,
minimizing the mean squared error (MSE). Architectural stability is
assessed via 5-fold cross-validation over the training partition.
Performance is monitored using two metrics: Exact Match (perfect vector
prediction) and Bit Accuracy (per-component accuracy, related to the
Hamming distance between predicted and optimal vectors). Cross-validation yields an Exact Match of $55.8\%$ ($\pm 1.6\%$). The
final production network is retrained on the full training partition. 

\textbf{Final Evaluation.}
Evaluated on the hold-out test set, the final model achieves an Exact Match
of $54.9\%$ and a Bit Accuracy of $81.5\%$. The negligible gap between
cross-validation and test performance confirms the absence of overfitting. An Exact Match of $\sim\!55\%$ indicates that the network identifies the
global optimum in more than half of the switching instances. The Bit
Accuracy of $81.5\%$ further reveals that, when the exact optimum is not
predicted, the proposed action typically differs from it by a single
component. In the context of the PSC acceptance criterion, near-optimal
proposals that differ from the optimum by a single inducer state frequently
satisfy the descent condition $J_{r^*} \ge J_{t_m}$ and are accepted. The
effective bypass rate of the exhaustive search therefore exceeds the
$55\%$ Exact Match baseline, as confirmed by the results reported in
Table~\ref{tab:psc_comparison}. The fallback mechanism in
Algorithm~\ref{alg:psc_accelerated} ensures that any proposal failing the
acceptance criterion is discarded in favour of the exhaustive search,
preserving the performance guarantees of the base controller.
For reproducibility, 
the datasets, training and cross-validation scripts, and optimized network
weights are publicly available at: https://github.com/ChristianFdz9/psc-neural-accelerator.git

\section{Auxiliary results}
\label{standard_theorems}

\noindent\textbf{Tonelli's Theorem.}
Let $(X, \mathcal{A}, \mu)$ and $(Y, \mathcal{B}, \nu)$ be $\sigma$-finite
measure spaces, and let $f: X \times Y \to [0, \infty]$ be a non-negative
measurable function. Then:
\begin{equation*}
\begin{split}
    &\int_X \left( \int_Y f(x,y)\,\mathrm{d}\nu(y) \right) \mathrm{d}\mu(x)\\
    &= \int_Y \left( \int_X f(x,y)\,\mathrm{d}\mu(x) \right) \mathrm{d}\nu(y)
    = \int_{X \times Y} f\,\mathrm{d}(\mu \times \nu).
    \end{split}
\end{equation*}

\noindent\textbf{Fubini's Theorem.}
Let $(X, \mathcal{A}, \mu)$ and $(Y, \mathcal{B}, \nu)$ be $\sigma$-finite
measure spaces. If $f: X \times Y \to \mathbb{R}$ is measurable and
$\int_{X \times Y} |f(x,y)|\,\mathrm{d}(\mu \times \nu) < \infty$, then:
\begin{equation*}
\begin{split}
    &\int_X \left( \int_Y f(x,y)\,\mathrm{d}\nu(y) \right) \mathrm{d}\mu(x)\\
    &\qquad = \int_Y \left( \int_X f(x,y)\,\mathrm{d}\mu(x) \right) \mathrm{d}\nu(y).
    \end{split}
\end{equation*}

\noindent \textbf{Bounded Perturbation Theorem \citep{Pazy83, EngelNagel}.} 
Let $X$ be a Banach space. If $A: D(A) \subset X \to X$ is the infinitesimal generator of a strongly continuous semigroup $S(t)$ on $X$, and $B \in \mathcal{B}(X)$ is a bounded linear operator, then the perturbed operator $L = A + B$ with domain $D(L) = D(A)$ generates a strongly continuous semigroup $T(t)$ on $X$. Furthermore, $T(t)$ can be represented by the Dyson-Phillips series:
\begin{equation*}
\begin{split}
&T(t) = \sum_{k=0}^{\infty} T_k(t), \\ &T_0(t) = S(t), \\ &T_{k+1}(t)x = \int_0^t S(t-s) B T_k(s)x \, \mathrm{d}s, \, \forall x \in X.
\end{split}
\end{equation*}

\noindent \textbf{Miyadera-Voigt Criterion \citep{Voigt77}).} 
Let $X = L^1(\Omega, \mu)$. Let $A$ generate a positive strongly continuous semigroup $S(t)$ on $X$, and let $B \in \mathcal{B}(X)$. If there exist constants $\hat{t}_0 > 0$ and $q \in (0,1)$ such that for all non-negative $p \in X$,
\begin{equation*}
\int_0^{\hat{t}_0} \| B S(s) p \|_{L^1} \, \mathrm{d}s \le q \|p\|_{L^1},
\end{equation*}
then the generated semigroup $T(t)$ associated with $A+B$ is a positive strongly continuous semigroup. Additionally, if $S(t)$ preserves total mass and $\int_\Omega (Bq)(x) \, \mathrm{d}x = 0$ for all $q \in X$, then $T(t)$ preserves mass for all non-negative initial conditions.

\noindent\textbf{Harris--Meyn--Tweedie geometric ergodicity \citep{Harris1956, MeynTweedie93}.}
Let $(U_t)_{t\ge0}$ be a mass-preserving, $\psi$-irreducible and aperiodic Markov semigroup on $L^1(\mathbb R_+^n)$, and let $L$ denote its infinitesimal generator. 

Assume that there exist a measurable Lyapunov function $V:\mathbb R_+^n\to[1,\infty)$, constants $a>0$ and $b<\infty$, a petite set $C\subset\mathbb R_+^n$, a time $T>0$, a probability measure $\nu$, and a constant $\alpha\in(0,1]$ such that
\[
LV(\mathbf x)\le -aV(\mathbf x)+b,
\qquad \forall \mathbf x\in\mathbb R_+^n,
\tag*{(Lyapunov drift)}
\]
and
\[
U_T(\mathbf x,\cdot)\ge \alpha\,\nu(\cdot),
\qquad \forall \mathbf x\in C.
\tag*{(uniform minorization)}
\]

Then there exists a unique invariant probability measure $\pi$ for $(U_t)_{t\ge0}$ with $\pi(V)<\infty$, and there exist constants $K\geq 1$ and $\phi>0$ such that for all initial densities $p_0$ with finite $V$-norm,
\[
\|U_t p_0 - \pi\|_{V}
\le K\,e^{-\phi t}\,\|p_0-\pi\|_{V},
\qquad \forall t\ge 0.
\]
\section{Mathematical spaces and functional notation}
\label{spaces_not}

In this work, we define $\mathbb{R}_+^n$ as the non-negative orthant of $\mathbb{R}^n$. The following notation is used:

\begin{itemize}[leftmargin=*]
    \item $L^1(\mathbb{R}_+^n)$: Space of measurable functions such that $\int_{\mathbb{R}_+^n} |f(x)| \, \mathrm{d}x < \infty$.
    
    \item $L^1_+(\mathbb{R}_+^n)$: Subspace of $L^1$ where $f(x) \ge 0$ a.e.
    
    \item $C([0, \infty); L^1)$: Functions $p(t, \cdot)$ continuous in the $L^1$ norm for $t \in [0, \infty)$.

    \item $C_c^\infty(\mathbb{R}_+^n)$: Infinitely differentiable functions with compact support.
    
    \item $\mathcal{B}(L^1)$: Set of bounded linear operators on $L^1(\mathbb{R}_+^n)$.
    
    \item $D(L^\dagger)$: Domain of the operator defined as $\{ f \in L^1 : L^\dagger f \in L^1 \}$.

    \item $\mathrm{Leb}(\cdot)$: Lebesgue measure on $\mathbb{R}^n$.
    
    \item $\Omega \subset \mathbb{R}_+^n$: Bounded numerical domain defined as the product of intervals $\prod_{i=1}^n [x_{i,0}, x_{i,f}]$.

    \item $S \subset \{0, 1\}^n$: Set of admissible control vectors.
\end{itemize}

\bibliography{references}

@article{Canizo2018,
   author = {José A. Cañizo and José A. Carrillo and Manuel Pájaro},
   issue = {1-2},
   journal = {Journal of Mathematical Biology},
   month = {1},
   pages = {373--411},
   publisher = {Springer Verlag},
   title = {{Exponential Equilibration of Genetic Circuits Using Entropy Methods}},
   volume = {78},
   year = {2019},
}

@article{Pajaro17,
  author       = {P{\'a}jaro, Manuel and Alonso, Antonio A. and Otero-Muras, Irene and V{\'a}zquez, Carlos},
  title        = {Stochastic modeling and numerical simulation of gene regulatory networks with protein bursting},
  journal      = {Journal of Theoretical Biology},
  volume       = {421},
  pages        = {51--70},
  year         = {2017},

}

@article{Pajaro18,
  author       = {P{\'a}jaro, Manuel and Otero-Muras, Irene and V{\'a}zquez, Carlos and Alonso, Antonio A.},
  title        = {{SELANSI: SemiLagrangian Numerical Simulation of Gene Regulatory Networks}},
  journal      = {Bioinformatics},
  volume       = {34},
  number       = {5},
  pages        = {893--895},
  year         = {2018},
  
}

@article{Vaghy2024,
  author       = {V{\'a}ghy, Mih{\'a}ly A. and Otero-Muras, Irene and P{\'a}jaro, Manuel and Szederk{\'e}nyi, G{\'a}bor},
  title        = {A Kinetic Finite Volume Discretization of the Multidimensional {PIDE} Model for Gene Regulatory Networks},
  journal      = {Bulletin of Mathematical Biology},
  volume       = {86},
  number       = {2},
  pages        = {22},
  year         = {2024},

}

@book{Pazy83,
  author       = {Pazy, Amnon},
  title        = {Semigroups of Linear Operators and Applications to Partial Differential Equations},
  publisher    = {Springer},
  address      = {New York},
  year         = {1983},
  series       = {Applied Mathematical Sciences},
  volume       = {44},

}

@article{Voigt77,
  author    = {J{\"u}rgen Voigt},
  title     = {On the perturbation theory for strongly continuous semigroups},
  journal   = {Mathematische Annalen},
  year      = {1977},
  volume    = {229},
  number    = {2},
  pages     = {163--171},
}

@book{EngelNagel,
  author       = {Engel, Klaus-Jochen and Nagel, Rainer},
  title        = {One-Parameter Semigroups for Linear Evolution Equations},
  publisher    = {Springer},
  address      = {New York},
  year         = {2000},
  series       = {Graduate Texts in Mathematics},
  volume       = {194},
  
}

@book{MeynTweedie93,
  author       = {Meyn, Sean P. and Tweedie, Richard L.},
  title        = {Markov Chains and Stochastic Stability},
  publisher    = {Springer},
  address      = {London},
  year         = {1993},

}

@book{Davis93,
  author       = {Davis, M. H. A.},
  title        = {Markov Models and Optimization},
  publisher    = {Chapman \& Hall},
  address      = {London},
  year         = {1993},
  series       = {Monographs on Statistics and Applied Probability},
  volume       = {49},
 
}

@article{Benaim2015,
  author       = {Bena{\"i}m, Michel and Le Borgne, St{\'e}phane and Malrieu, Florent and Zitt, Pierre-Andr{\'e}},
  title        = {Qualitative properties of certain piecewise deterministic {Markov} processes},
  journal      = {Annales de l’Institut Henri Poincar{\'e} B, Probability and Statistics},
  volume       = {51},
  number       = {3},
  pages        = {1040--1075},
  year         = {2015},

}

@article{Harris1956,
  author       = {Harris, Theodore E.},
  title        = {The existence of stationary measures for certain {Markov} processes},
  journal      = {Proceedings of the Third Berkeley Symposium on Mathematical Statistics and Probability},
  volume       = {2},
  pages        = {113--124},
  year         = {1956},
  publisher    = {University of California Press}
}

@article{Miyadera66,
  author    = {Isao Miyadera},
  title     = {On perturbation theory for semi-groups of operators},
  journal   = {Tohoku Mathematical Journal},
  year      = {1966},
  volume    = {18},
  number    = {3},
  pages     = {299--310},
}

@article{Faquir2025,
  author = {Faquir, Hamza and P\'ajaro, Manuel and Otero-Muras, Irene},
  title = {A computational framework for optimal and Model Predictive Control of stochastic gene regulatory networks},
  journal = {IEEE Transactions on Computational Biology and Bioinformatics},
  year = {2025},
  
}

@article{Lugagne2017,
  author = {Lugagne, Jean-Baptiste and Sosa Carrillo, Sebasti\'an and Kirch, Melanie and K\"ohler, Agnes and Batt, Gregory and Hersen, Pascal},
  title = {Balancing a genetic toggle switch by real-time feedback control and periodic forcing},
  journal = {Nature Communications},
  volume = {8},
  pages = {1671},
  year = {2017},

}

@article{Gardner2000,
  author = {Gardner, Tim S. and Cantor, Charles R. and Collins, James J.},
  title = {Construction of a genetic toggle switch in Escherichia coli},
  journal = {Nature},
  volume = {403},
  pages = {339--342},
  year = {2000},
 
}

@article{Lugagne2024,
  author = {Lugagne, Jean-Baptiste and Blassick, Caroline M. and Dunlop, Mary J.},
  title = {Deep model predictive control of gene expression in thousands of single cells},
  journal = {Nature Communications},
  volume = {15},
  pages = {2148},
  year = {2024},

}

@article{Oduola2017,
  author = {Oduola, Wuraola O. and Li, Xin and Duan, Cheng and Qian, Lin and Wu, Fan and Dougherty, Edward R.},
  title = {Time-Based Switching Control of Genetic Regulatory Networks: Toward Sequential Drug Intake for Cancer Therapy},
  journal = {Cancer Informatics},
  volume = {16},
  pages = {1176935117706888},
  year = {2017},
  
}

@article{Fernandez2025,
  author = {Fern\'andez, Christian and Faquir, Hamza and Vaghy, Mih\'aly A. and P\'ajaro, Manuel and Szederk\'enyi, G\'abor and Otero-Muras, Irene},
  title = {{PIDE} models for efficient control of stochastic gene regulatory circuits},
  journal = {IFAC-PapersOnLine},
  volume = {59},
  number = {19},
  pages = {621--626},
  year = {2025},
  
}

@article{Lohmiller1998,
  author = {Lohmiller, Winfried and Slotine, Jean-Jacques E.},
  title = {On contraction analysis for non-linear systems},
  journal = {Automatisierungstechnik},
  volume = {46},
  number = {11},
  pages = {10--21},
  year = {1998}
}

@article{Nielsen20,
  author = {Nielsen, R. F. and Gernaey, K. V. and Mansouri, S. S.},
  title = {Hybrid machine learning assisted modelling framework for particle processes},
  journal = {Computers \& Chemical Engineering},
  volume = {140},
  pages = {106916},
  year = {2020},
  
}

@misc{Zhang25_arxiv,
  author = {Zhang, X. and Huang, K. W. and Vo, D.-N. and Han, M. and Decardi-Nelson, B. and Yin, X.},
  title = {Machine learning-based hybrid dynamic modeling and economic predictive control of carbon capture process for ship decarbonization},
  howpublished = {arXiv preprint arXiv:2502.05833},
  year = {2025},

}

@article{Chen22,
  title        = {Large Scale Model Predictive Control with Neural Networks and Primal Active Sets},
  author       = {Chen, Steven W. and Wang, Tianyu and Atanasov, Nikolay and Kumar, Vijay and Morari, Manfred},
  journal      = {Automatica},
  volume       = {135},
  pages        = {109947},
  year         = {2022},
  
}

@article{Menolascina2011,
  author = {Menolascina, Filippo and Di Bernardo, Mario and Di Bernardo, Diego},
  title = {Analysis, design and implementation of a novel scheme for in-vivo control of synthetic gene regulatory networks},
  journal = {Automatica},
  volume = {47},
  number = {6},
  pages = {1265--1270},
  year = {2011},
 
}

@article{Herrera2024_OSS,
  title        = {An Optimal Switching Sequence Model Predictive Control Scheme for the {3L-NPC} Converter with Output {LC} Filter},
  author       = {Herrera, Felipe and Mora, Andrés and Cárdenas, Roberto and Díaz, Matías and Rodríguez, José and Rivera, Marco},
  journal      = {Processes},
  volume       = {12},
  number       = {2},
  pages        = {348},
  year         = {2024},
}

@article{FCSMPC2024,
  title        = {Finite control set model predictive control with limit cycle stability guarantees},
  author       = {Xu, Duo and Lazar, Mircea},
  journal      = {Automatica},
  volume       = {181},
  year         = {2025},
  pages        = {112507},

}

@article{Mosca2005,
  author = {Mosca, Edoardo},
  title = {Predictive switching supervisory control of persistently disturbed input-saturated plants},
  journal = {Automatica},
  volume = {41},
  number = {1},
  pages = {55--67},
  year = {2005},

}

@inproceedings{Brancato2023,
  author = {Brancato, Sara Maria and De Lellis, Francesco and Salzano, Davide and Russo, Giovanni and di Bernardo, Mario},
  title = {External control of a genetic toggle switch via Reinforcement Learning},
  booktitle = {Proceedings of the 2023 European Control Conference (ECC)},
  year = {2023},
  address = {Bucharest, Romania},
  pages = {1--6},

}

@article{Guarino2020Balancing,
  author = {Guarino, Agostino and Fiore, Davide and Salzano, Davide and di Bernardo, Mario},
  title = {Balancing Cell Populations Endowed with a Synthetic Toggle Switch via Adaptive Pulsatile Feedback Control},
  journal = {ACS Synthetic Biology},
  volume = {9},
  number = {4},
  pages = {793--803},
  year = {2020},

}

@article{Fernandez2022,
  author = {Fern\'andez, Christian and Faquir, Hamza and P\'ajaro, Manuel and Otero-Muras, Irene},
  title = {Feedback control of stochastic gene switches using {PIDE} models},
  journal = {IFAC-PapersOnLine},
  volume = {55},
  number = {18},
  pages = {62--67},
  year = {2022},
  
}

@Book{FB-CTDS,
  author =    {F. Bullo},
  title =     {Contraction Theory for Dynamical Systems},
  year =      2026,
  edition =   {{1.3}},
  publisher = {Kindle Direct Publishing},
  ISBN =      {979-8836646806},

}

@article{SAKURAI2022110647,
title = {Interval analysis of worst-case stationary moments for stochastic chemical reactions with uncertain parameters},
journal = {Automatica},
volume = {146},
pages = {110647},
year = {2022},
author = {Yuta Sakurai and Yutaka Hori}
}

@article{FIORE2016279,
title = {Contraction analysis of switched systems via regularization},
journal = {Automatica},
volume = {73},
pages = {279-288},
year = {2016},
author = {Davide Fiore and S. John Hogan and Mario {di Bernardo}}
}

@article{SUN2025112101,
title = {Event-triggered PDF shape control of non-Gaussian stochastic system},
journal = {Automatica},
author = {Xiaoyang Sun and Ping Zhou},
volume = {173},
pages = {112101},
year = {2025}}

@INPROCEEDINGS{11312970,
  author={Zand, Armin M. and Gupta, Ankit and Khammash, Mustafa},
  booktitle={2025 IEEE 64th Conference on Decision and Control (CDC)}, 
  title={Control with Practical Guarantees of Stationary Variance in Stochastic Chemical Reaction Networks}, 
  year={2025},
  volume={},
  number={},
  pages={2911-2916}
  }

\end{document}